\documentclass[10pt,a4paper]{article}

\usepackage{amsmath}
\usepackage{multirow}
\usepackage{slashbox}
\usepackage{float}
\usepackage{graphicx}
\usepackage[colorlinks=true]{hyperref}
\usepackage{amsfonts}
\usepackage{color}
\usepackage{epsfig}
\usepackage{amssymb}
\usepackage{srcltx}
\def\shft{{{\hbox{\tiny\rm drg}}}}
\def\Conv{{\hbox{\rm Conv}}}
\def\PP{{{P}}}
\def\sG{{\hbox{\tiny\rm sG}}}
\def\maj{{\hbox{\rm\tiny maj}}}
\def\mm{{\hbox{\rm\tiny mm}}}
\def\risk{{\hbox{\rm Risk}}}

\def\rank{{\hbox{\rm rank}}}

\def\inter{\mathop{\hbox{\rm int}}}
\def\myepsilonstar{{\varepsilon_\star}}

\def\Risk{\hbox{\rm risk}}

\def\Erf{{\mathop{\hbox{\small\rm Erf}}}}

\def\bR{{\mathbf{R}}}

\def\Opt{{\hbox{\rm Opt}}}

\def\L{{\cal L}}

\def\SS{{\overline{\cal{S}}}}

\def\T{{\cal T}}
\def\E{{\cal E}}
\def\P{{\cal P}}

\def\R{{\cal R}}
\def\G{{\cal G}}
\def\H{{\cal H}}

\def\N{{\cal N}}

\def\S{{\cal S}}

\def\bE{{\mathbf{E}}}
\def\Var{{\mathbf{Var}}}

\def\Prob{\hbox{\rm Prob}}
\def\sign{\hbox{\rm sign}}
\def\Erf{\hbox{\rm Erf}}
\newtheorem{corollary}{Corollary}[section]
\newtheorem{lemma}{Lemma}[section]
\newtheorem{proposition}{Proposition}[section]
\newtheorem{remark}{Remark}[section]
\newtheorem{definition}{Definition}[section]
\newtheorem{assump}{Assumption}

\def\qed{\hfill$\Box$}
\def\Card{\mathop{\hbox{\rm Card}}}
\def\e{{\rm e}}

\def\Erf{{\mathop{\hbox{\small\rm Erf}}}}

\newcommand{\be}{\begin{eqnarray}}
\newcommand{\ee}[1]{\label{#1}\end{eqnarray}}
\newcommand{\nn}{\nonumber \\}

\newcommand{\four}{ \mbox{\small$\frac{1}{4}$}}

\newcommand{\rf}[1]{~(\ref{#1})}
\newcommand{\half}{ \mbox{\small$\frac{1}{2}$}}

\newcommand{\bse}{\begin{eqnarray*}}
\newcommand{\ese}{\end{eqnarray*}}

\definecolor{MyDarkBlue}{rgb}{0,0.08,0.45}
\definecolor{MyViolet}{rgb}{0.45,0.08,0.95}
\definecolor{MyBrown}{rgb}{0.45,0.08,0}

\def\nui{{\hbox{\rm\tiny nui}}}
\def\inp{{\hbox{\rm\tiny inp}}}

\newcommand{\hide}[1]{{}}

\title{Hypothesis Testing via Euclidean Separation}
\author{Vincent Guigues
 \thanks{School of Applied Mathematics FGV/EMAp, 22 250-900 Rio de Janeiro, Brazil
 {\tt vincent.guigues@fgv.br}}
 \and
Anatoli Juditsky
\thanks{LJK, Universit\'e Grenoble Alpes, 700 Avenue Centrale 38041 Domaine Universitaire
de Saint-Martin-d'H\`{e}res, France,
{\tt anatoli.juditsky@imag.fr}}
\and Arkadi Nemirovski
\thanks{Georgia Institute of Technology, Atlanta, Georgia
30332, USA, {\tt nemirovs@isye.gatech.edu}}}
\date{}
\begin{document}
\maketitle
\begin{abstract}
We discuss an ``operational'' approach to testing convex composite hypotheses when the underlying distributions are heavy-tailed. It relies upon
Euclidean separation of convex sets and can be seen as an extension of the approach to testing by convex optimization developed in \cite{GJN2015,AffineDetectors}. In particular, we show how one
can construct quasi-optimal testing procedures for families of distributions which are majorated, in a certain precise sense, by a sub-spherical symmetric one and study the
relationship between tests based on Euclidean separation and ``potential-based tests.'' We apply the promoted methodology in the problem
of sequential detection and illustrate its practical implementation in an application to sequential detection of changes in the input of a dynamic system.
\end{abstract}



\section{Introduction}
The following important observation, attributed to H. Chernoff \cite{chernoff1952} (see also \cite{Burnashev1979,Burnashev1982}), was  the starting point of our research.
\begin{quotation}
\noindent Let $X_1$ and $X_2$ be two nonempty closed and convex sets, one of them being bounded, in $\bR^n$. Suppose that, given a noisy observation
\be
\omega=x+\xi
\ee{jeq0} of the unknown signal $x\in X_1\cup X_2$, where $\xi\sim \N(0,I_n)$ -- the standard $n$-dimensional Gaussian vector, one wants to decide on the hypotheses $H_1:\;x\in X_1$ vs. $H_2:\;x\in X_2$. Then, assuming that $X_1$ and $X_2$ do not intersect (the decision problem is clearly unsolvable otherwise), optimal tests (with respect to different definitions of maximal risks) can be obtained using the following simple construction:\\
1) solve the following (convex optimization) Euclidean separation problem
\be
\Opt=\min_{x^1\in X_1,\,x^2\in X_2} \half \|x^1-x^2\|_2
\ee{eq:l21}
where $\|\cdot\|_2$ is the Euclidean distance.
\\2) Given an optimal solution $(x_*^1,x_*^2)$ to \rf{eq:l21}, compute
\[
h_* =\frac{ x_*^1 - x_*^2 }{\|x^1_* - x^2_* \|_2}\mbox{  and }c_* =\half h_*^T ( x_*^1 + x_*^2 ).
\]
Then the test  $\T_*$ which accepts $H_1$ if
$h_*^T\omega-c_* \geq 0$ and accepts $H_2$ otherwise minimizes the maximal risk of testing -- the maximal over $x\in X_1\cup X_2$ probability of rejecting the true hypothesis -- over the class of all (deterministic or randomized) tests.
\footnote{A test which accepts $H_1$ if $h_*^T\omega-c_* \geq \alpha$ and accepts $H_2$ otherwise, with a properly chosen $\alpha$, is also minimax optimal in several other settings, e.g., Neyman-Pearson problem, etc.} Furthermore, the risk of $\T_*$ is easily computable: the maximal probability of wrongly rejecting the true hypothesis is
\[
\frac{1}{\sqrt{2 \pi }}\int_{\Opt}^{\infty} \e^{-t^2/2}dt=\Erf(\Opt),
\]
where $\Erf(\cdot)$ is the standard error function.
\end{quotation}

This simple observation had important theoretical consequences (see \cite{Ingster2002,Lehmann2006}). Surprisingly, its ``practical implications'' have been largely overseen. Indeed, when the problem  \rf{eq:l21} can be solved efficiently,
\footnote{what is the case when sets $X_1$ and $X_2$ allow for a ``computationally efficient description.'' We refer the reader to \cite{BN2001} for precise definitions and details on efficient implementability. For the time being, it is sufficient to assume
that  \rf{eq:l21} can be solved using {\tt CVX} \cite{cvx2014}.} one can assemble pair-wise tests into multiple-testing procedures to build provably (nearly) optimal tests for a wide class of Gaussian decision problems (e.g., various detection
problems \cite{Burnashev1979,Burnashev1982,Seq2015,alltogether}).   Then  in \cite{GJN2015} the corresponding framework was extended to ``good,''  in a certain precise sense, parametric families of distributions, which include, aside from the Gaussian family, families of Poisson and discrete distributions, the ``common denominator'' of these developments being the fact that for these
families the {\em near-optimal}   (plain optimal, in the case of Gaussian family) tests can be built upon using {\em affine detectors} which can be found by convex optimization. Later, {\em affine and quadratic detectors} were studied in a more general setting in \cite{AffineDetectors}.
\par
In this paper we extend  \cite{GJN2015}, \cite{AffineDetectors} studying the application of tests based on ``Euclidean separation'' to the problems where the distribution of the observation noise $\xi$ has ``heavy tails.'' In particular, in Sections \ref{sec:spher} and \ref{sec:esep} we discuss the problem of testing  convex hypotheses for
{\em sub-spherical families} of distributions, which include Gaussian and Gaussian mixture distributions such as multivariate Student \cite{kotz2004multivariate} and multivariate Laplace \cite{eltoft2006multivariate} distributions. We study the relationship of
sub-spherical families with detector-based tests, and show how ``good'' detectors can be built for sub-spherical families, as well as for some other (e.g., sub-Gaussian) families of distributions in Section \ref{sec:potent0}. Then in Section \ref{sectchp} we explain how tests based on Euclidean separation of pairs of convex hypotheses can be used to construct sequential change detection procedures. Finally, in Section \ref{changepoint} we present a numerical illustration of the proposed techniques: we implement sequential decision rules, developed in Sections  \ref{sectschemeI} and \ref{sectschemeII} for a toy problem of detecting changes in the input of a dynamical system -- changes in the trend of a simple time series.
\section{Basic theory}\label{sectEuclsep}
\subsection{Pairwise hypothesis testing: Situation and goal}
\label{sectsitg}
\paragraph{The basic problem} we intend to consider {\sl in a nutshell} is as follows: we are given observation
\begin{equation}\label{jeq1}
\omega=x+\xi,
\end{equation}
where $x\in\bR^n$ is an unknown signal, and $\xi$ is a random noise with probability density $p(\cdot)$, taken with respect to the Lebesgue measure, known to belong to some given family $\P$.
Our basic goal is, given two nonempty closed convex sets $X_\chi\subset\bR^n$, $\chi=1,2,$ with one of these sets bounded, to decide, via observation (\ref{jeq1}), on the hypotheses $H_\chi$, $\chi=1,2$, with $H_\chi$ stating that the signal $x$ underlying observation belongs to $X_\chi$.
In other words, we are to decide upon the families $\P_1$ and $\P_2$ of distributions, where $\P_\chi$ is the family of distributions of random vectors $x+\xi$, $\xi\sim p$ with $p\in\P$ and $x\in X_\chi$. For the sake of brevity, we shall simplify the description of the hypotheses $H_\chi$ to ``$x\in X_\chi$.''
\paragraph{Stationary repeated observations.}
We ``embed'' the just defined inference problem in the family of inference problems $(\S_K)$, $K=1,2,...$, where $(\S_K)$ is the problem of deciding whether $x\in X_1$ or $x\in X_2$ via sample $\omega^K=(\omega_1,...,\omega_K)$ of $K$ independent observations
\begin{equation}\label{jeq2}
\omega_k=x+\xi_k,\,k=1,...,K,
\end{equation}
with independent noises $\xi_k\sim p(\cdot)\in\P$; we refer to observations (\ref{jeq2}) with independent across $k$ noises $\xi_k\sim p\in \P$ as to {\sl stationary $K$-repeated observations}.
\paragraph{Semi-stationary repeated observations.}
Inference problem $(\S_K)$ can be viewed as a special case of a more general inference problem $(\SS_K)$ as follows.
Suppose that $\{X_1^k,X_2^k:1\leq k\leq K\}$ is a collection of nonempty convex and closed sets such that at least one of the set in every pair $(X^k_{1},X^k_{2})$ is bounded. Let also $\P^1,...,\P^K$ be $K$ given families of probability densities with respect to the Lebesgue measure. Suppose that the observation $\omega^K=(\omega_1,...,\omega_K)$ is given by
\begin{equation}\label{jeqSS2}
\omega_k=x_k+\xi_k,\,k=1,...,K,
\end{equation}
where $\{x_k\}_{k=1}^K$ is a deterministic sequence, $\xi_k\sim p_k$ are independent across $k$ noises, and $\{p_k\in \P^k\}_{k=1}^K$ is a deterministic sequence.\footnote{One can easily verify (cf. \cite[Section 3.1.2]{GJN2015})
that the constructions and results to follow remain intact when the assumptions on observations (\ref{jeqSS2}) are weakened to the assumption that $x_k$ and $\xi_k$ are random, and the conditional
distribution of $\xi_k$, given $x_1,...,x_k$ and $\xi_1,..,\xi_{k-1}$,
always belongs to $\P^k$.} In the sequel, we refer to observations (\ref{jeqSS2}) satisfying the just imposed restrictions on
$\{x_k\}_{k=1}^K$ and $\{\xi_k\}_{k=1}^K$, as to {\sl semi-stationary $K$-repeated observations.} Our objective in problem $(\SS_K)$ is to decide, via the observations \rf{jeqSS2}, on the hypotheses $H_\chi$, $\chi=1,2$, with $H_\chi$ stating that $x_k\in X^k_\chi$ for all $k=1,2,...,K$.
\paragraph{A simple test} for $(\S_K)$ or $(\SS_K)$ is, by definition, a function $\T_K(\omega^K)$, $\omega^K=(\omega_1,...,\omega_K)$, taking values $\{1,2\}$, with $\T_K(\omega^K)=\chi$ interpreted as ``given observation $\omega^K$, the test accepts $H_\chi$ and rejects the alternative.'' We define the {\em partial risks} $\risk_{\chi\S}(\T_K|\P,X_1,X_2)$, $\chi=1,2$
of a test
$\T_K$ on the inference problem $(\S_K)$ as
\[
\begin{array}{rcl}
\risk_{1\S}(\T_K|\P,X_1,X_2)&=&\sup_{x\in X_1}\sup_{p(\cdot)\in \P}\Prob_{x,p}\{\T_K(\omega^K)=2\},\\
\risk_{2\S}(\T_K|\P,X_1,X_2)&=&\sup_{x\in X_2}\sup_{p(\cdot)\in \P}\Prob_{x,p}\{\T_K(\omega^K)=1\},\\
\end{array}
\]
where $\Prob_{x,p}$ stands for the probability with respect to the distribution of observations \rf{jeq2}.
In other words, $\risk_{\chi\S}(\T_K|\P,X_1,X_2)$, is the worst-case probability for $\T_K$ to reject $H_\chi$ when the hypothesis is true.
The partial risks  $\risk_{\chi\SS}(\T_K|[\P^k,X_1^k,X_2^k]_{k=1}^K)$, $\chi=1,2$, of a test $\T_K$ on the inference problem $(\SS_K)$ are defined similarly, but now the supremum is taken with respect to {\em all} deterministic sequences $\{x_k\in X^k_\chi\}_{k=1}^K$ and $\{p_k\in\P^k\}_{k=1}^K$ participating in (\ref{jeqSS2}). Finally, the risks of $\T_K$ on $(\S_K)$ and $(\SS_K)$ are defined as
$$
{
\begin{array}{rcl}
\risk_\S(\T_K|\P,X_1,X_2)&=&\max\limits_{\chi=1,2}\risk_{\chi\S}(\T_K|\P,X_1,X_2),\\
\risk_\SS\left(\T_K|[\P^k,X^k_1,X^k_2]_{k=1}^K\right)&=&\max\limits_{\chi=1,2}\risk_{\chi\SS}\left(\T_K|[\P^k,X^k_1,X^k_2]_{k=1}^K\right).
\end{array}}
$$

\paragraph{Our intention} is to investigate the performance of specific tests stemming from ``Euclidean Separation'' of $X_1$ and $X_2$ {(or $X_1^k$ and $X_2^k$)} to be described in a while.
\subsection{Sub-spherical families of distributions}\label{sec:spher}
We shall be primarily interested in {\sl sub-spherical families $\P$.}
\subsubsection{Sub-spherical families of distributions: definition and basic examples}\label{sec:bexamples}
\begin{definition}\label{defspherfam} {\bf A.}
A \underline{sub-spherical family of distributions} $\P=\P_\gamma^n$ on $\bR^n$ is specified by an even probability density $\gamma(\cdot)$ on the axis such that $\gamma$ is positive in a neighbourhood of the origin. $\P_\gamma^n$ is comprised of all probability densities $p(\cdot)$ on $\bR^n$ such that $p(\cdot)$ is even, and
\begin{equation}\label{zeq1}
\forall (e\in\bR^n,\|e\|_2=1,\,\delta\geq0):\;
\int\limits_{e^T\xi\geq\delta} p(\xi)d\xi\leq P_\gamma(\delta):=\int\limits_\delta^\infty \gamma(s)ds.
\end{equation}
that is, $p$-probability mass of a half-space not containing a neighbourhood of the origin is upper-bounded by $P_\gamma(\delta)$, where $\delta$ is the distance from the origin to the half-space, and
\begin{equation}\label{zeq0}
P_\gamma(r)=\int\limits_r^\infty \gamma(s)ds:\;\bR\to[0,1].
\end{equation}
\par {\bf B.}
A sub-spherical family $\P=\P^n_\gamma$ is called \underline{monotone}, if it contains a \underline{cap} -- a spherically symmetric density
$q(\xi)=f(\|\xi\|_2)$
where $f$ is nonincreasing on the nonnegative axis and such that the induced by $q$ density
of the distribution of $e^T\xi$, $\|e\|_2=1$,
is exactly $\gamma(\cdot)$.
Note that whenever this is the case, $\gamma(\cdot)$ is nonincreasing on the nonnegative ray.
\par{\bf C.}
We call a function $\gamma$ on the real axis \underline{nice}, if $\gamma$ is an even probability density which is continuous  and is nonincreasing on the nonnegative ray. A sub-family $\P^n$ of a sub-spherical family $\P^n_\gamma$ is called \underline{completely monotone},
 if $\gamma$ is nice, and  for every $p(\cdot)\in \P^n$ and every $e\in\bR^n$, $\|e\|_2=1$, the random scalar variable $e^T\xi$, $\xi\sim p$, has probability density
$\gamma_{e,p}(\cdot)$, and this density is nice. Note that due to  $\P^n\subset \P^n_\gamma$, it holds
\begin{equation}\label{aeq1}
\int\limits_\delta^\infty\gamma_{e,p}(s)ds\leq\int\limits_\delta^\infty \gamma(s)ds\,\;\;\forall \delta\geq0.
\end{equation}
\end{definition}
In the sequel, we simplify the notation $\P^n_\gamma$ to $\P_\gamma$ when the value of $n$ is clear from the context.
\par

\paragraph{Example: Gaussian scale mixtures.} Consider the situation where
\begin{equation}\label{zeq4}
\xi\sim \sqrt{Z}\eta,
\end{equation}
where $Z$ is a scalar a.s. positive random variable with given probability distribution $P_Z(t)=\Prob\{Z\leq  t\},t\geq0$ such that {$P_Z(0)=0$}, and $\eta\sim\N(0,\Theta)$ is
a zero mean Gaussian $n$-dimensional random vector independent of $Z$ with unknown a priori positive definite covariance matrix $\Theta$ which is known to be $\preceq I_n$. We refer to $\Theta$ as to {\sl matrix parameter} of the distribution of $\xi$.
Given a unit vector $e$ and $\delta\geq0$, we have
\bse
\begin{array}{rcl}
\Prob\{e^T\xi\geq\delta\}&=&\int\limits_{t>0} \left[\Prob
\{e^T\eta\geq t^{-1/2}\delta\}\right]dP_Z(t)\nn
&=&\int\limits_{t>0} \left[\Prob_{\zeta\sim \N(0,I_n)} \{e^T\Theta^{1/2}\zeta\geq t^{-1/2}\delta\}\right]dP_Z(t)\nn
&=&\int\limits_{t>0} \Erf\left(t^{-1/2}\delta/\sqrt{e^T\Theta e}\right)dP_Z(t)
\leq \int\limits_{t>0}\Erf(t^{-1/2}\delta)dP_Z(t),
\end{array}
\ese
(here $\Erf(\cdot)$ is the standard error function),
implying that
\[
\Prob\{e^T\xi\geq\delta\}=\int\limits_\delta^\infty \gamma_{e,\Theta,P_Z}(s)ds\leq\int\limits_\delta^\infty \gamma_{P_Z}(s)ds\;\;\forall\delta\geq0,
\]
where
\[
\begin{array}{rcl}
\gamma_{e,\Theta,P_Z}(s)&=&\int\limits_{t>0}{1\over\sqrt{2\pi te^T\Theta e}}\exp\left\{-{s^2\over 2te^T\Theta e}\right\}dP_Z(t),\\
\gamma_{P_Z}(s)&=&\int_{t>0}{1\over\sqrt{2\pi t}}\exp\left\{-{s^2\over 2t}\right\}dP_Z(t),\,-\infty<s<\infty.
\end{array}
\]
Clearly, $\gamma_{e,\Theta,P_Z}(s)$ and $\gamma_{P_Z}(s)$ are  positive even probability densities nonincreasing on the nonnegative ray; these densities are continuous, provided
\begin{equation}\label{providedthat}
\int_{t>0}t^{-1/2}dP_Z(t)<\infty.
\end{equation}
Thus, the family of distributions of random vectors (\ref{zeq4}) with $Z$ and $\eta$ as explained above (we refer to these distributions as to {\sl Gaussian scale mixtures}) is contained in the sub-spherical family
$\P_{\gamma_{P_Z}}$. The latter family clearly is monotone,  the cap being the probability density of $\sqrt{Z}\eta$ with independent $Z\sim P_Z$ and
$\eta\sim\N(0,I_n)$. Besides this, in the case of (\ref{providedthat}) the
family $\P_{\gamma_{P_Z}}$ is completely monotone.
\par
A standard example of a Gaussian scale mixture is {given by} {\em $n$-variate  $t$-distributions} $t_n(q,\Theta),\;\Theta\preceq I_n$ (multivariate Student distributions with $q$ degrees of freedom, see \cite{kotz2004multivariate} and references therein). Here $t_n(q,\Theta)$ is, by definition, the distribution of the random vector
  $\xi=\sqrt{Z}\eta$ with $Z=q/\zeta$, where $\zeta$ is the independent of $\eta\sim\N(0,\Theta)$  random variable following $\chi^2$-distribution with $q$ degrees of freedom. One can easily see that all one-dimensional projections $e^T\xi$, $\|e\|_2=1$, of $\xi\sim t_n(q,I_n)$ are random variables with univariate $t_q$-distribution, implying that the multidimensional densities in question form a completely monotone sub-family of
  the sub-spherical family $\P_{\gamma_\S}$ where $\gamma_\S$ is the density of Student's $t_q$ distribution with $q$ degrees of freedom.
 \par
 Another example of scheme \rf{zeq4}
 is the {\em $n$-variate Laplace distributions} $\L_n(\lambda,\Theta),\;\Theta\preceq I_n$, where $Z$ {is}  exponentially distributed with parameter $\lambda$. In this case all one-dimensional projections $e^T\xi$, $\|e\|_2=1$, of $\xi\sim \L_n(\lambda,I_n)$, obey the Laplace distribution with parameter $\lambda$, whence the distributions in question form a completely monotone sub-family of
 the sub-spherical family $\P_{\gamma_\L}$, where $\gamma_{\L}$ is the Laplace density
 \begin{equation}\label{ueq1}
 \gamma_\L(s)=(2\lambda)^{-1}\e^{-|s|/\lambda}, \;s\in \bR.
 \end{equation}
 Finally, with $Z$ taking value 1 with probability 1, scheme \rf{zeq4} describes Gaussian distributions with zero mean and covariance matrices $\preceq I_n$; all these distributions form a completely monotone sub-famuily of
 the sub-spherical family $\P_{\gamma_\G}$, where $\gamma_\G$ is the standard univariate Gaussian density:
 $$
 \gamma_\G(s)={1\over\sqrt{2\pi}}{\rm e}^{-s^2/2}.
 $$
\subsubsection{``Calculus'' of sub-spherical families of distributions}
Sub-spherical families of distributions and their completely monotone subfamilies admit a kind of ``calculus'' with the basic rules which follow.
\par
The following two facts are immediate:
\begin{proposition}\label{prop22} Let $\P_\gamma^n$ be a sub-spherical family of distributions, and let $x\mapsto Qx:\bR^n\to\bR^m$ be an onto mapping satisfying $QQ^T\preceq I_m$. Whenever $p(\cdot)\in\P^n_\gamma$, the distribution of the random vector $Q\xi$, $\xi\sim p$, belongs to $\P^m_\gamma$. Moreover, if $QQ^T=I_m$ and $\P^n_\gamma$ has a cap $q$, then $\P^m_\gamma$ has a cap as well; this cap is the density of the random vector $Q\xi$, $\xi\sim q$.
\end{proposition}
\begin{proposition}\label{prop222} A sub-spherical family of distributions is closed with respect to taking convex combinations of its members. Besides this, the union $\P_{\gamma_1}\cup \P_{\gamma_2}$ of two sub-spherical families of distributions is contained in the
sub-spherical family $\P_\gamma$ with
$$
\gamma(-s)=\gamma(s)=-{d\over ds}\max\limits_{i=1,2} \int\limits_s^\infty \gamma_i(r)dr,\;s\geq 0.
$$
\end{proposition}
Complete monotonicity is preserved by taking sums. The precise statement is as follows:
\begin{proposition}\label{propnew} Let $\mu$ and $\nu$ be nice,
 let  $\P^n$ be a subfamily of the sub-spherical family $\P_\mu^n$, and let $\P^m$
 be a completely monotone subfamily of the sub-spherical family $\P_\nu^m$. Given $r\times n$ matrix $A$, $r\times m$ matrix $B$ and positive definite matrix $\Theta$ such that
\begin{equation}\label{aeq2}
\Theta^2\succeq  AA^T,\,\Theta^2\succeq BB^T,\, AA^T+BB^T\succ0,
\end{equation}
consider random vectors of the form
\begin{equation}\label{aeq3}
\xi=\Theta^{-1}[A\eta+B\zeta],
\end{equation}
where $\eta\sim p(\cdot)\in\P^n$ and $\zeta\sim q(\cdot)\in \P^m$ are independent. Let also
\begin{equation}\label{aeq4}
\gamma(s)=\int\limits_{-\infty}^\infty \mu(s-r)\nu(r)dr.\\
\end{equation}
Then
\par(i) $\gamma(\cdot)$ is nice, and for every $e\in\bR^r$, $\|e\|_2=1$,
and every $\delta\geq0$ one has
\begin{equation}\label{aeq5}
\Prob\{e^T\xi\geq\delta\}\leq\int\limits_\delta^\infty \gamma(s)ds.
\end{equation}
Besides this, the scalar random variable $e^T\xi$ possesses symmetric density, which combines with {\rm (\ref{aeq5})} to imply that the distribution of $\xi$ belongs to the sub-spherical family $\P_\gamma^r$.
\par
(ii)
If, in addition to the above assumptions,  $\P^n$ is completely monotone, then the family of distributions of random variables {\rm (\ref{aeq3})} induced by $\eta\sim p\in\P^n$ and $\zeta\sim q\in \P^m$ is a completely monotone subfamily of the sub-spherical family $\P^r_\gamma$.
\end{proposition}
For a proof, see Section \ref{propnewproof}.
\par
As an immediate consequence of Proposition \ref{propnew}, we get the following
\begin{corollary}\label{cornew} For $1\leq i\leq I<\infty$, let $\mu_i$  be nice functions on the axis, and let $\P^{n_i}$ be subfamilies of the sub-spherical families $\P_{\mu_i}^{n_i}$ such that at least $I-1$ of these subfamilies are completely monotone. Given $r\times n_i$ matrices $A_i$ such that $A_iA_i^T+A_jA_j^T\succ 0$ whenever $i\neq j$, let $\Theta\succ0$ be such that
$$
\Theta^2\succeq A_iA_i^T,\,1\leq i\leq I.
$$
Let $\P^r$ be the family of probability distributions of random vectors of the form
$$
\xi=\Theta^{-1}\sum_{i=1}^IA_i\eta_i
$$
where $\eta_1,...,\eta_I$ are independent of each other and such that $\eta_i\sim p_i$ for some $p_i\in\P^{n_i}$. Then $\P^r$ is contained in the sub-spherical family $\P^r_\gamma$, where
$$
\gamma=\mu_1\star\mu_2\star...\star\mu_I
$$
is nice. If all $\P^{n_i}$ are completely monotone, then so is $\P^r$.
\end{corollary}

\subsection{Euclidean separation and associated tests}\label{sec:esep}
In the sequel, we fix the entities $\P$, $X_1$, $X_2$, $H_1$, $H_2$  introduced in the beginning of Section \ref{sectsitg} and assume that $X_1\cap X_2=\emptyset$ (otherwise no test 
can decide on $H_1$ vs. $H_2$ with risk $<1/2$).

\subsubsection{Single-observation Euclidean separation test}\label{sec:soeuclidean}
 Consider the optimization problem
\begin{equation}\label{jeq3}
\Opt=\min_{x^1\in X_1,x^2\in X_2} \|x^1-x^2\|_2,
\end{equation}
and let $x^1_*$, $x^2_*$ form an optimal solution to the problem (since both $X_1$, $X_2$ are nonempty, closed and convex, and one of the sets is bounded, an optimal solution does exist). We set
\begin{equation}\label{jeq4}
s_*(\omega)=h_*^T\omega-c_*,\,\,h_*={[x^1_*-x^2_*]\over \|x_*^1-x_*^2\|_2},\,
c_*=\half\left[\min_{x\in X_1}h_*^Tx+\max_{x\in X_2} h_*^Tx\right]=\half h_*^T(x^1_*+x^2_*).
\end{equation}
Note that while (\ref{jeq3}) may have many optimal solutions, the vector $h_*$ and the real $c_*$ are uniquely defined by $X_\chi$, $\chi=1,2$. The affine function $s_*(\cdot)$ possesses the following properties:
\begin{equation}\label{jeq5}
\begin{array}{rcl}
x\in X_1&\Rightarrow& s_*(x)\geq s_*(x^1_*)=\delta,\,\;\delta:={1\over 2}\|x_*^1-x_*^2\|_2,\\
x\in X_2&\Rightarrow&s_*(x)\leq s_*(x^2_*)=-\delta.
\end{array}
\end{equation}
\par We can associate with $h_*$ and $c_*$ the {\sl Euclidean separation test} $\T_1$ which, given observation $\omega$  (\ref{jeq1}) with $x$ known to belong to $X_1\cup X_2$, accepts the hypothesis $H_1:x\in X_1$ and rejects the hypothesis $H_2:x\in X_2$ when $s_*(\omega)\geq0$, and accepts $H_2$ and rejects $H_1$ otherwise.\par
Let us make the following immediate observation:
\begin{proposition}\label{prop1} In the situation of Section \ref{sectsitg} and in the notation from {\rm (\ref{jeq4}), (\ref{jeq5})}, let $\P=\P_\gamma$ be a sub-spherical family of distributions, and let $\alpha_1\geq0$ and $\alpha_2\geq0$ be such that $\alpha_1+\alpha_2\leq 2\delta$.
Whenever $x\in X_1$ and $p\in\P$, the $p$-probability of  the event $\{\xi:s_*(x+\xi)<\half(\alpha_2-\alpha_1)\}$ is at most $P_\gamma(\alpha_1)$,  and when $x\in X_2$, the $p$-probability of the event $\{\xi:s_*(x+\xi)\geq\half(\alpha_2-\alpha_1)\}$ is at most $P_\gamma(\alpha_2)$. As a result,
risks of the test $\T_1$ satisfy
\[\begin{array}{l}
\risk_{1\S}(\T_1|\P,X_1,X_2)\leq P_\gamma(\alpha_1),
\\
\risk_{2\S}(\T_1|\P,X_1,X_2)\leq P_\gamma(\alpha_2).
\end{array}
\] In particular, when $\alpha_1=\alpha_2=\delta$, the risk $\risk_\S(\T_1|\P,X_1,X_2)$ of the test is at most
\begin{equation}\label{weq14}
\myepsilonstar=\myepsilonstar(\delta|\gamma):=P_\gamma(\delta)=\int\limits_\delta^\infty \gamma(s)ds.
\end{equation}
\end{proposition}
For a proof, see Section \ref{prop1proof}.
 \begin{remark}\label{rem1} In the situation of Proposition \ref{prop1}, let the sub-spherical family $\P_\gamma$ be monotone. Then the test   $\T_1$ described in the proposition has the minimal risk among all single-observation tests deciding on $H_1$ vs. $H_2$.
  \end{remark}
  Indeed, denoting by $q(\cdot)$ the cap of $\P_\gamma$, and by $p_1(\cdot)$ and $p_2(\cdot)$ the densities of random vectors $x^1_*+\xi$, $x^2_*+\xi$, $\xi\sim q(\cdot)$, we clearly have
$$
\int\min[p_1(\omega),p_2(\omega)]d\omega=2\myepsilonstar,
$$
implying by the Neyman-Pearson lemma that the risk of any test deciding on two simple hypotheses $x=x_1^*$, $x=x_2^*$  via a single observation (\ref{jeq1}) is at least $\myepsilonstar$.
\subsubsection{Majority tests based on Euclidean separation}\label{majtest}
Let $K$ be a positive integer, and let $\P=\P_\gamma$ be a sub-spherical family of distributions. The  {\sl $K$-observation majority test} $\T^\maj_K$  for problem $(\S_K)$ works as follows: given observations $\omega_k$, $k=1,...,K$, see \rf{jeq2}, the test accepts $H_1$ and rejects $H_2$ when $s_*(\omega_k)\geq0$ for at least $K/2$ values of $k$, and accepts $H_2$ and rejects $H_1$ otherwise. From Proposition \ref{prop1} it follows that the risk  of $\T^\maj_K$ satisfies the bound
\begin{equation}\label{jeq20}
\risk_{\S}(\T^\maj_K|\P,X_1,X_2)\leq \sum_{K\geq k\geq K/2} \left({K\atop k}\right)\myepsilonstar^k(1-\myepsilonstar)^{K-k}.
\end{equation}
Let us assume that in the problem $(\SS_K)$ the families $\P^1,...,\P^K$ in the definition of $(\SS_K)$ are sub-spherical families of distributions, and that the sets $X_1^k$ and $X_2^k$ do not intersect for $k=1,...,K$.
The majority test $\T^\maj_K$ can be easily modified to become applicable to problem $(\SS_K)$.
Let the affine function $s_k(\cdot)$ and positive real $\delta_k$ be the entities $s_*(\cdot),\delta$ associated, via  \rf{jeq3} -- \rf{jeq5}, with $\P^k$ in the role of $\P$ and $X^k_\chi$ in the role of $X_\chi$, $\chi=1,2$.  The  {\sl $K$-observation majority test} $\T^\maj_K$  for the problem $(\SS_K)$ given observations $\omega_k$, $k=1,...,K$, see (\ref{jeqSS2}), accepts $H_1$ and rejects $H_2$ when $s_k(\omega_k)\geq0$ for at least $K/2$ values of $k$, and accepts $H_2$ and rejects $H_1$ otherwise.
Let now $\gamma_k$ be the density underlying the sub-spherical family $\P^k$. When applying Proposition \ref{prop1} to $\P=\P^k$, $X_1=X_1^k$, and $X_2=X_2^k$, we conclude that the risk $\risk_\S(\T_1|\P^k,X^k_1,X^k_2)$ of the test $\T_1$ in the problem of deciding, {\em given a single observation $\omega_k$}, upon the hypotheses $H_1^k:\,x_k\in X_1^k$ vs
 $H_2^k:\,x_k\in X_2^k$,  does not exceed
$\epsilon_{k}(\delta_k)=\int\limits_{\delta_k}^\infty\gamma_k(s)ds<\half.$ We conclude that the risk of the majority test $\T^\maj_K$ in the problem $(\SS_K)$ satisfies the bound \footnotemark[1]
\footnotetext[1]{We use the following fact (if absolutely evident facts indeed exist, this is one of them): given a vector $p=[p_1;...;p_n]$ with entries $p_i\in [0,1]$, and a positive integer $k\leq n$, let $P_{\geq k|n}(p)$ be the probability to get $\geq k$ heads in $n$ independent flips of a coin, with probability $p_i$ to get a head in the $i$-th trial. Then $P_{\geq k|n}(p)$ is a nondecreasing function of $p$. Whatever ``evident,'' this fact needs a proof, and here it is.
Let $B=\{\omega\in\bR^n:\;0\leq \omega_i\leq 1,\,i\leq n\}$ be the $n$-dimensional cube equipped with the Lebesgue measure $\mu$. We have $P_{\geq k|n}(p)=\mu( B_p)$, where $B_p=\{\omega\in B:\,\Card\{i:\omega_i\leq p_i\}\geq k\}$. When $0\leq p_i\leq p_i^\prime\leq1$, $1\leq i\leq n$, we clearly have $B_p\subset B_{p'}$, and consequently $P_{\geq k|n}(p)=\mu( B_p)\leq \mu(B_{p'})=P_{\geq k|n}(p')$, as claimed.  }
\be
\risk_\SS\left(\T^\maj_K|[\P^k,X^k_1,X^k_2]_{k=1}^K\right)\leq \sum_{K/2\leq k\leq K} p_{k|K},
\ee{jeq20SS}
where $p_{k|K}=p_{k|K}(\epsilon_{1}(\delta_1),...,\epsilon_{K}(\delta_K))$ is the probability of  $k$ successes in the first $K$ non-stationary independent Bernoulli trials with the probability $\epsilon_k(\delta_k)$ of success in  $k$-th trial. Observe that $p_{j|k}$, $0\leq j\leq k$, $1\leq k\leq K$, satisfy the recursion:
\be
p_{j|k}=(1-\epsilon_k)p_{j|k-1}+\epsilon_k p_{j-1|k-1},\;\mbox{with $p_{0|0}=1$},
\ee{pnb}
and where, by convention, $p_{-1|k}=0$, $k=0,...,K-1$. Recursion \rf{pnb} allows to compute all the quantities $p_{j|k}$, $0\leq j\leq k\leq K$, and thus the right hand side of (\ref{jeq20SS}), in $O(K^2)$ arithmetic operations.
\subsubsection{Near-Optimality}\label{sectNearOpt}
We are about  to show that under appropriate assumptions, the majority test built in Section \ref{majtest} is near-optimal. The precise
statement is as follows:
\begin{proposition}\label{OptMajTest} In the situation and in the notation described in the beginning of Section \ref{sec:esep}, assume that sub-spherical family
$\P=\P_\gamma$ and positive reals $\bar{d}$, $\alpha$, $\beta$ are such that
\begin{equation}\label{eqDD}
{\beta \bar d \leq \frac{1}{2}},
\end{equation}
\begin{equation}\label{eqDDb}
\int_0^\delta \gamma(s)ds \geq \beta\delta,\;\;0\leq\delta\leq \bar{d},
\end{equation}
and $\P$ contains a density $q(\cdot)$ such that
\begin{equation}\label{eqDDa}
\int_{\bR^n}\sqrt{q(\xi-e)q(\xi+e)}d\xi \geq \exp\{-\alpha e^Te\}\;\;\forall (e:\|e\|_2\leq \bar{d}).
\end{equation}
Let, further, the sets $X_1$, $X_2$ be  such that $\Opt$ as given by {\rm (\ref{jeq3})} satisfies the relation
\begin{equation}\label{eqDDc}
\delta:=\Opt/2\leq\bar{d}.
\end{equation}
Given tolerance $\epsilon\in(0,1/5)$, the risk of $K$-observation majority test
$\T^\maj_K$ for problem $(\S_K)$ associated with $X_1$, $X_2$, $\P_\gamma$ ensures the relation
\begin{equation}\label{eqDDd}
K\geq K^*:=\left\rfloor{\ln(1/\epsilon)\over2\beta^2\delta^2}\right\lfloor \;\Rightarrow \; \risk_{\S}(\T^\maj_K|\P,X_1,X_2)\leq\epsilon
\end{equation}
(here $\rfloor x\lfloor$ stands for the smallest integer $\geq x\in \bR$).
In addition, for every $K$-observation test $\T_K$ for $(\S_K)$ satisfying $\risk_{\S}(\T_K|\P,X_1,X_2)\leq\epsilon$ it holds
\begin{equation}\label{eqDDg}
K\geq K_*:={\ln\left({1\over 4\epsilon}\right)\over2\alpha\delta^2}.
\end{equation}
As a result, the majority test $\T^\maj_{K^*}$ for $(\S_{K^*})$ has risk at most $\epsilon$ and is near-optimal, in terms of the required number of observations,  among all tests
with risk $\leq\epsilon$: the number $K$ of observations in such test satisfies the relation
$$
K^*/K\leq \theta:=K^*/K_*=O(1){\alpha\over\beta^2}.
$$
\end{proposition}
For proof, see Section \ref{ProofOptMajTest}.
\paragraph{Illustration.} Given $\nu\geq1$, consider the case when $\P=\P_\gamma$ is the sub-spherical monotone family with $n$-variate (spherical) Student's $t_n(\nu,I_n)$-distribution in the role of the cap, so that
\begin{equation}\label{eqDDz}
\gamma(s)=\gamma_\nu(s):={\Gamma\left({\nu+1\over 2}\right)\over \Gamma\left({\nu\over 2}\right)(\pi\nu)^{1/2}}
\left[1+s^2/\nu\right]^{-(\nu+1)/2}.
\end{equation}
It is easily seen (see Section \ref{ProofOptMajTest}) that $\P$ contains the $\N(0,\half I_n)$ density $q(\cdot)$, implying that setting $$
\bar{d}=1,\;\alpha=1,\;\beta=\gamma_1(1)={1\over 2\pi},
$$
one ensures  relations (\ref{eqDD}), (\ref{eqDDb}) and (\ref{eqDDc}). As a result, when $\Opt$ as yielded by (\ref{jeq3}) is $\leq 2$, the non-optimality factor $\theta$ of the majority test $\T^\maj_{K^*}$ as defined in Proposition \ref{OptMajTest} does not exceed $O(1)$.\footnote{In fact, accurate computation shows that when $\epsilon\leq 0.1$ and $\Opt\leq 1000$, the non-optimality factor $\theta$ does not exceed 17 when $\nu\geq1$, and does not exceed 10 when $\nu\geq 4$.}

\subsection{Potential-based tests with Euclidean separation}\label{sec:potent0}
\subsubsection{Potentials and potential-based tests}\label{sec:potent}
\begin{definition}\label{defpotrisk} A \underline{potential} is an odd and nondecreasing Borel real-valued function $\eta(\cdot)$ on the axis.
Given a family $\P$ of probability densities on $\bR^n$, a nonnegative $\delta$  and a potential $\eta(\cdot)$, we define the {\sl $\delta$-risk} $\Risk_\delta(\eta|\P)$ of the potential on $\P$ as the smallest $\epsilon$ such that
\begin{equation}\label{zeq10}
\begin{array}{lrcl}
(a)&\int{\rm e}^{-\eta(\delta+e^T\xi)}p(\xi)d\xi&\leq& \epsilon\;\;\forall (e\in\bR^n:\|e\|_2=1,p\in\P),\\
(b)&\int{\rm e}^{\eta(e^T\xi-\delta)}p(\xi)d\xi&\leq& \epsilon\;\;\forall (e\in\bR^n:\|e\|_2=1,p\in\P).
\end{array}
\end{equation}
\end{definition}
Let us make the following immediate observation:
\begin{proposition}\label{prop0}
For $k=1,...,K$, let $\P^k$ be a family of probability densities on $\bR^n$,  $X^k_1$ and $X^k_2$ be closed nonempty non-intersecting convex sets in $\bR^n$, one of the sets being bounded, and let $h_{k},c_k$ and $\delta_{k}$ be associated   with $\P=\P^k$, $X_\chi=X_\chi^k$, $\chi=1,2$, via  {\rm (\ref{jeq4}) and (\ref{jeq5})}. Given potentials $\eta_1,...,\eta_K$, let us define the Euclidean detector induced by the potential $\eta_k$ and $h_k$, $c_k$ as the function
\begin{equation}\label{req1}
\phi_k(\omega)=\eta_k(h_k^T\omega-c_k):\;\bR^n\to\bR,
\end{equation}
and let
\begin{equation}\label{req2}
\phi^{(K)}(\omega_1,...,\omega_K)=\sum_{k=1}^K\phi_k(\omega_k):\;\bR^{nK}\to\bR,\;\,K=1,2,...
\end{equation}
Finally, let $\{x_k\}_{k=1}^K$ and $\{p_k\in\P^k\}_{k=1}^K$ be deterministic sequences.
Then
\[
\begin{array}{lrcl}
(a)&\int{\rm e}^{-\phi^{(K)}(x_1+\xi_1,...,x_K+\xi_K)}\prod_{k=1}^K[p_k(\xi_k)d\xi_k]&\leq&
\prod_{k=1}^K\Risk_{\delta_k}(\eta_k|\P^k)\;\;\forall (x_k\in X^k_1,\,p_k\in\P^k)_{k=1}^K,\\
(b)&\int{\rm e}^{\phi^{(K)}(x_1+\xi_1,...,x_K+\xi_K)}\prod_{k=1}^K[p_k(\xi_k)d\xi_k]&\leq&
\prod_{k=1}^K\Risk_{\delta_k}(\eta_k|\P^k)\;\;\forall (x_k\in X^k_2,\,p_k\in\P^k)_{k=1}^K.
\end{array}
\]
\end{proposition}
For a proof, see Section \ref{prop0proof}. \par
Under the premise of Proposition \ref{prop0}, consider a test $T^\eta_K$ which, given $K$ observations {\rm (\ref{jeqSS2})} with independent of each other $\xi_k\sim p_k\in\P^k$, accepts the hypothesis $H_1:\{x_k\in X^k_1,k\leq K\}$ whenever $\phi^{(K)}(\omega_1,....,\omega_K)\geq0$, and accepts $H_2:\{x_k\in X_2^k,k\leq K\}$ otherwise.
An immediate corollary of Proposition \ref{prop0} is as follows:
\begin{corollary}\label{cor1} In the notation and under the assumptions from the premise of Proposition \ref{prop0},
combining Proposition \ref{prop0} with the Markov inequality we obtain that
the risk $\risk_\SS(\T^\eta_K|[\P^k,X_1^k, X_2^k]_{k=1}^K)$ of $\T^\eta_K$ does not exceed the quantity
$\prod_{k=1}^K\Risk_{\delta_k}(\eta_k|\P^k)$.
\end{corollary}
\subsubsection{Potentials for sub-spherical families of distributions}\label{sectsfpot}
\begin{definition}\label{defpotreg} {\bf A.} Given $\delta>0$, we call a potential $\eta(s):\bR\to \bR$ $\delta$-regular, if the function
$$H_{\delta\eta}(s)={\rm e}^{-\eta(\delta-s)}+{\rm e}^{-\eta(\delta+s)}$$
is nondecreasing on the ray $s\geq0$.
\par
{\bf B.}
The {\sl $\delta$-index} of the $\delta$-regular potential $\eta$ on a sub-spherical family $\P_\gamma$ is the quantity
\begin{equation}\label{jeq8}
\begin{array}{rcl}
\epsilon_\delta(\eta|\gamma):=\int\limits_{-\infty}^\infty {\rm e}^{-\eta(r)}\gamma(r-\delta)dr
=\int\limits_{-\infty}^\infty{\rm e}^{-\eta(\delta+s)}\gamma(s)ds
=\int\limits_0^\infty H_{\delta\eta}(s)\gamma(s)ds,\\
\end{array}
\end{equation}
where the concluding equality is due to the fact that $\gamma(\cdot)$ is even.
\end{definition}
\begin{proposition}\label{prop2} Let $\P=\P_\gamma$ be a sub-spherical family of probability densities on $\bR^n$, let $\delta\geq0$, and let $\eta$ be a $\delta$-regular potential. Then
\begin{equation}\label{jeq9}
\Risk_\delta(\eta|\P_\gamma)\leq \epsilon_\delta(\eta|\gamma).
\end{equation}
\end{proposition}
For a proof, see Section \ref{prop2proof}.
\paragraph{Example 1: step potential.} Assume that $\P=\P_\gamma$ is a sub-spherical family of distributions, and let $\delta>0$, implying that $\myepsilonstar=\myepsilonstar(\delta|\gamma)$, as given by (\ref{weq14}), belongs to $[0,1/2)$ (recall that $\gamma$ is an even probability density positive in a neighbourhood of the origin). 
We define the {\sl step potential} as\footnote{Hereafter we put $\sign(s)=s/|s|$ if $s \neq 0$ and $\sign(0)=0$.}
\begin{equation}\label{jeq10}
\eta(s)=\half \ln\left(
{1-\myepsilonstar\over \myepsilonstar}
\right)\sign(s).
\end{equation}

Taking into account that $0\leq\myepsilonstar<1/2$, it is immediately seen that $\eta$ is a
$\delta$-regular potential. The $\delta$-index of the step potential satisfies
\begin{equation}\label{eq110}
\epsilon_\delta(\eta|\gamma)=2\sqrt{\myepsilonstar(1-\myepsilonstar)},
\end{equation}
as is shown by the following computation:
$$
\begin{array}{rcl}
\epsilon_\delta(\eta|\gamma)&=&\int\limits_{-\infty}^\infty {\rm e}^{-\eta(s)}\gamma(s-\delta)ds=\sqrt{1-\myepsilonstar\over\myepsilonstar}\int\limits_{-\infty}^{-\delta}\gamma(s)ds+\sqrt{\myepsilonstar\over 1-\myepsilonstar}\int\limits_{-\delta}^\infty \gamma(s)ds\\
&=&\sqrt{1-\myepsilonstar\over\myepsilonstar}\myepsilonstar+\sqrt{\myepsilonstar\over 1-\myepsilonstar}(1-\myepsilonstar)=2\sqrt{\myepsilonstar(1-\myepsilonstar)}.
\end{array}
$$
Note that in the situation of Section \ref{sectsitg} with stationary $K$-repeated observations (\ref{jeq2}), a sub-spherical family $\P=\P_\gamma$ and non-intersecting $X_1$, $X_2$, the test $\T^\eta_K$ associated with the step potential $\eta$ is exactly the majority test $\T^\maj_K$ defined in Section \ref{majtest}. Because the
index $\epsilon_\delta(\eta|\gamma)$ of the step potential is $<1$ due to $\myepsilonstar<1/2$,  its $\delta$-risk $\Risk_\delta(\eta|\P_\gamma)$ is $<1$ as well.
\paragraph{Example 2: ramp potential.} Assume that $\P$ is the sub-spherical family $P_{\gamma_\L}$ with
$\gamma_\L(s)=\frac{1}{2 \lambda} \e^{-|s|/\lambda}$, $s\in \bR$ (see the ``Laplace'' example in Section \ref{sec:spher}).
Given $\delta>0$, let us consider the {\em ramp potential} $\eta_\delta$ which minimizes the risk $\epsilon_\delta(\eta|\gamma_\L)$:
\[
\eta_\delta(s)=\half \ln\left({\gamma_\L(s-\delta)\over \gamma_\L(s+\delta)}\right)=
\left\{
\begin{array}{ll} {s/\lambda} &\mbox{for $|s|\leq \delta$,}\\
{\frac{\delta}{\lambda}}\,\sign(s),&\mbox{for $|s|> \delta$. }\end{array}\right.
\]
One can easily verify that $\eta_\delta$ is $\delta$-regular for any $\delta>0$, and the corresponding
$\delta$-index satisfies
\[
\epsilon_\delta(\eta_\delta|\gamma_\L)=\int\limits_{-\infty}^\infty\sqrt{\gamma_\L(s-\delta)\gamma_\L(s+\delta)}ds=
{\e^{-\delta/\lambda}\Big(\frac{\delta}{\lambda}+1\Big),}
\]
which is $<1$ for all $\delta>0$.
\par Examples above show that in the situation of Section \ref{sectsitg}, assuming that $\P=\P_{\gamma}$ is a sub-spherical family, there exist potentials $\eta$ with
$\delta$-indexes $\epsilon_\delta(\eta|\gamma)<1$
and therefore, using Proposition \ref{prop2}, with
risks $\Risk_\delta(\eta|\P)<1$, provided $\delta>0$.
In this situation, when solving problem $(\SS_K)$ with $X_1^k\equiv X_1,X_2^k\equiv X_2$, $X_1$ and $X_2$ not intersecting, and specifying $\delta>0$ according to (\ref{jeq5}), we can decide on the hypotheses $H_1$ and $H_2$ with any desired risk $\epsilon\in(0,1)$ via  semi-stationary $K$-repeated observations (\ref{jeqSS2}), provided that
\be
K\geq K(\epsilon):=
\left\rfloor
{\ln(\epsilon^{-1})\over \ln(\Risk_\delta(\eta|\P)^{-1})}
\right\lfloor
,
\ee{eq:K_potent}
see Corollary \ref{cor1}.

\begin{remark}
{Under the premise and in the notation of Proposition \ref{OptMajTest}, let
$\eta$ be a potential satisfying $\Risk_\delta(\eta|\P) \leq 2 \sqrt{\myepsilonstar (1-\myepsilonstar)}$ with $\myepsilonstar$ given by {\rm (\ref{weq14})} (i.e.,
in terms of its $\delta$-risk, $\eta$ is not worse than the step potential). By {\rm (\ref{eqDDb})} we have $2 \sqrt{\myepsilonstar (1-\myepsilonstar)}\leq (1-4\beta^2\delta^2)^{1/2}\leq
\exp\{-2\beta^2\delta^2\}$, implying that $K(\epsilon)$ as given by {\rm (\ref{eq:K_potent})} is at most the quantity $K^*$ given by {\rm (\ref{eqDDd})}. As a result, in the situation under consideration the induced by $\eta$ $K$-observation
test $\T^\eta_K$ for problem $(\SS_K)$ shares the  near-optimality properties of the majority test $\T^\maj_K$ stated in Proposition \ref{OptMajTest}.}
\end{remark}

\subsubsection{Potentials for the sub-Gaussian family}\label{secd:pot_subg}
Aside from sub-spherical families, there are other families of probability distributions allowing, in the case of $X_1\cap X_2=\emptyset$, for potentials with risks $<1$ and thus for tests with an arbitrarily low risk, provided stationary or semi-stationary $K$-repeated observations with properly selected $K$ are available. The
 simplest family of this type is the family $\P^n_\sG$ of sub-Gaussian  probability densities $p(\cdot)$, with parameters $0$ and $I_n$, that is, probability densities $p$ on $\bR^n$ such that
\begin{equation}\label{zeq33}
\int {\rm e}^{h^T\xi} p(\xi)d\xi\leq \e^{{1\over 2} h^Th}\;\;\forall h\in\bR^n.
\end{equation}
Assuming that $\P=\P_\sG$ and given $\delta>0$, let us put
\begin{equation}\label{zeq44}
\eta_{\sG,\delta}(s)=\delta s.
\end{equation}
\begin{proposition} \label{prop33} Whenever $\delta\geq0$, we have
\begin{equation}\label{zeq20}
\Risk_\delta(\eta_{\sG,\delta}|\P_\sG)\leq{\rm e}^{-\delta^2/2}.
\end{equation}
\end{proposition}
Proposition \ref{prop33} is a special case of Proposition 3.3 in \cite{AffineDetectors}; to make the presentation self-contained, we provide its proof in Section \ref{sec:proofprop33}.
\par
Note that if $\P$ is the sub-spherical family $\P_{\gamma_\G}$ with
$\gamma_\G(s)= {1\over\sqrt{2\pi}}{\rm e}^{-s^2/2}$, $s\in \bR$
(recall that this family contains, for instance, Gaussian distributions with zero mean and covariance matrix $\preceq I_n$, cf. Section \ref{sec:bexamples})
the potential $\eta_{\sG,\delta}(s)=\half \ln \left({\gamma_\G(s-\delta)\over \gamma_\G(s+\delta)}\right)=\delta s$ minimizes the $\delta$-index
over $\P_{\gamma_\G}$ with
$\epsilon_\delta( \eta_{\sG,\delta} |\gamma_\G)=  {\rm e}^{-\delta^2/2}$.
\subsubsection{``Majority of means'' tests}
Let now $X_1{\subset} \bR^n$ and $X_2 {\subset} \bR^n$ be closed convex nonempty sets, one of the sets being bounded, and such that $X_1\cap X_2=\emptyset$. In this situation, given $\epsilon\in (0,1)$ and
$K$-repeated stationary observations $\omega^K$, the potential-based tests developed in Sections \ref{sec:potent} -- \ref{secd:pot_subg}
allow to decide on $H_1$, $H_2$ with risk $\leq \epsilon$, where $K$ grows logarithmically with $\epsilon^{-1}$.
For instance, in the case of a sub-Gaussian family of distributions, the corresponding test
attains the risk $\epsilon$ provided that
$K\geq K(\epsilon)=O\left({\ln(\epsilon^{-1})/\delta^{2}}\right)$ (cf. \rf{zeq20}), where $2\delta>0$ is the Euclidean distance between $X_1$ and $X_2$, see \rf{jeq3} --\rf{jeq5}. On the other hand, these tests
rely upon the ``Cramer-type''   Definition \ref{defpotrisk} of the risk of the potential, and thus assume control of the exponential moment of $\eta(\xi)$. In this section we present a different testing procedure, which only uses second order characteristics of the potential (mean and variance), and yet
allows to achieve arbitrarily low risks of testing for essentially the same sizes of observation sample.
We describe this modification in the simplest case of stationary repeated observations, generalizations to more general settings (e.g., that of  Proposition \ref{prop0}) being straightforward.
\par
Now, let
$x_*^1,\,x_*^2$ and $h_*$ be associated with $X_1$ and $X_2$ via \rf{jeq3} and \rf{jeq4}, and let $\P$ be a family of probability distributions on $\bR^n$. We suppose that a potential $\eta(\cdot)$ and $c\in \bR$ are such that for some $\varrho>0$ and all $p\in \P,$

\be
\begin{array}{ll}
(a)&\bE_{\xi\sim p}\{\eta(h_*^T(x_*^1 +\xi)+c)\}-\bE_{\xi\sim p}\{\eta(h_*^T(x_*^2 + \xi)+c)\}\geq  \varrho,\\
(b)&\Var_{\xi\sim p}\{\eta(h_*^T(x+\xi)+c)\}
\leq 1
\end{array}
\ee{eq:evar}
for all $x\in X_1\cup X_2$.\footnote{Here and below $\Var$ denotes the ``usual'' variance: for a probability density $p$ on $\bR^n$ and $f:\,\bR^n\to\bR$, $\Var_{\xi\sim p}\{f(\xi)\}=\bE_{\xi\sim p}\{f^2(\xi)\}-[\bE_{\xi\sim p}\{f(\xi)\}]^2$.}
\paragraph{Example.} Let $\P$ be a family of zero-mean distributions on $\bR^n$ with covariance matrix $\preceq\sigma^2 I_n$:
\[
\int \xi p(\xi)d\xi=0,\;\;\int \xi\xi^Tp(\xi)d\xi\preceq I_n,\;\;\forall p\in \P.
\]
For the linear potential $\eta(t)=t$ with $c=0$ we clearly have
\[\begin{array}{l}
\bE_{\xi\sim p}\{\eta(h_*^T(x_*^1 +\xi){+c)}\}-\bE_{\xi\sim p}\{{\eta(}h_*^T(x_*^2 + \xi){+c)}\}=h_*^T(x_*^1 - x_*^2)=2\delta=:\varrho,\;\;\forall p\in \P\\
\Var_{\xi\sim p}\{{\eta(} h_*^T(x+\xi){+c)} \}=\bE_{\xi\sim p}\{(h_*^T\xi)^2\}\leq 1\;\;\forall (x\in X_1\cup X_2, \;p\in \P)
\end{array}
\]
(recall that $\delta=\half \|x^1_* - x^2_* \|_2$, and  $\|h_*\|_2=1$).
\par
Now, given
$K$-repeated stationary observations ${\omega^K}=[x+\xi_1,...,x+{\xi_K}]$, and $\epsilon_1,\epsilon_2\in (0,1)$, we consider the inference problem $(\S_K)$ of deciding via the observation $\omega^K$ whether $x\in X_1$ (hypothesis $H_1$) or $x\in X_2$ (hypothesis $H_2$). Our objective is to build the test $\T_K$ for $\S^K$ such that, uniformly over $x\in X_1\cup X_2$ and $p\in \P$, the probability of wrongly rejecting $H_1$ (accepting $H_2$) is $\leq \epsilon_1$,
and the probability of
wrongly  rejecting $H_2$ (accepting $H_1$) is $\leq \epsilon_2$.
and we want to attain this goal using the smallest possible size $K$ of observation sample $\omega^K$.
\par
For the sake of definiteness, assume that $\epsilon_1\leq \epsilon_2$. We denote
$\kappa=\half\left(1-{\ln (\epsilon_2^{-1})\over \ln (\epsilon_1^{-1})}\right)$ (note that $0\leq \kappa<\half$). Let now
$m=\big\rfloor{4\e(\e^{-\kappa}+1)^2\varrho^{-2}}\big\lfloor$, and let
\[
\psi_{j}(\omega^K)={1\over m}\sum_{i=(j-1)m+1}^{mj} \eta(h_*^T \omega_i + c),\;\;j=1,2,...
\]
Denote
\[
c_*=\bE_{\xi\sim p}\{\eta(h_*^T(x_*^1 + \xi)+c)\}-{\e^{\kappa}\varrho\over 1+\e^{\kappa}} \geq
\bE_{\xi\sim p}\{\eta(h_*^T(x_*^2 + \xi)+c)\}+{\varrho\over 1+\e^{\kappa}}.
\]

The {\em $K$-observation majority of means test} $\T_K^{\mm}$ for $\S_K$ is as follows: given
\be
K:=Jm\geq \rfloor 2\ln(\epsilon_1^{-1})\lfloor\, m
\ee{eq:Kdef}
observations $\omega_k$, $k=1,...,K$, $\T_K^{\mm}$ accepts the hypothesis $H_1$ when $\psi_{j}(\omega^K)\geq c_*$ for at least
$ J/2$ values of $j$, and accepts the hypothesis $H_2$ otherwise.
\begin{proposition}\label{prop133}
In the just described situation, the risks of the test $\T^{\mm}_K$ meet the problem specifications, namely,
\[
\risk_{1 {\S}}(\T_K|\P,X_1,X_2)\leq \epsilon_1,\;\;\risk_{2{\S}}(\T_K|\P,X_1,X_2)\leq \epsilon_2.
\]
\end{proposition}
For a proof, see Section \ref{sec:proofprop133}.

Let us consider the test $\T^\mm_K$ using the linear potential {with $c=0$}. Under  the premise of the proposition, i.e., in the situation where the noise covariance matrix is $\preceq I_n$ and the distance between the sets $X_1$ and $X_2$ is $\geq 2\delta$, the size of the stationary $K$-repeated observation sufficient for the test to satisfy the risk specifications is
\[
\big\rfloor{\e(\e^{-\kappa}+1)^2\delta^{-2}}\big\lfloor\;\rfloor 2\ln(\epsilon_1^{-1})\lfloor=O\big(\delta^{-2}\ln (\epsilon_1^{-1})\big).
\]
Note that when $\epsilon_1=\epsilon_2=\epsilon$, the above number of observations sufficient to decide on $H_1$, $H_2$ with risk $\epsilon$ is within absolute constant factor of the number of observations, as given by (\ref{eq:K_potent}) and (\ref{zeq20}), needed for the same purpose in the case when $\P=\P_\sG$.

\section{Sequential detection via Euclidean separation}\label{sectchp}

\subsection{Motivating example}\label{sec:motivatingexample}

\par We start introducing a motivating example which will help understand the general setting discussed in Sections \ref{dynsitgoal},
\ref{sectschemeI}, \ref{sectschemeII} and which will be used in Section \ref{numericalill}
to test our methodology.

Consider {the simple} time series model
\begin{equation}\label{meq1}
\begin{array}{rcl}
y_k&=&\alpha_k+\zeta_k\\
\alpha_k&=&\alpha_{k-1}+\eta_k+u_k\\
\end{array},\,k=1,2,...,d,
\end{equation}
where
\begin{itemize}
\item  $y_k$ is the observation at time $k$, and $u=[u_1;...;u_{d}]\in\bR^{d}$ is the deterministic input,
\item $\zeta=[\zeta_1;...;\zeta_{d}]$ is zero mean ${d}$-dimensional Gaussian random vector with unknown covariance matrix known to be $\preceq \sigma^2I_{d}$;
\item $\eta=[\eta_1;...;\eta_{d}]$ is independent of $\zeta$ ${d}$-dimensional random vector  obeying multivariate Student distribution with $\nu$ degrees of freedom and matrix parameter $\preceq I_{d}$. Here $\nu$ is a positive integer or $+\infty$, with $\nu=\infty$ interpreted as the fact that $\eta$ is
a zero mean Gaussian random vector with covariance matrix $\preceq I_{d}$.
\end{itemize}
\par We intend to decide from observations $y_1,...,y_{d}$ on {the nuisance} hypothesis $u=0$ vs. signal alternative
\begin{equation}\label{meq2}
u\in\bigcup\limits_{{1\leq i\leq {d},\atop0<\rho< R}}
[U_i^+(\rho)\cup U_i^-(\rho)],
\end{equation}
where $R\in(0,\infty)$ is a parameter, $U_i^-(\rho)=-U_i^+(\rho)$, and $[U_i^+(\rho)\cup U_i^-(\rho)]$ is the set of ``signal inputs of shape $i$ and magnitude $\geq\rho$.'' We consider two cases:
\begin{itemize}
\item {\sl pulse signal inputs:} $U_i^+(\rho)=\{u\in\bR^{d}:u_k=0,k\neq i, \rho \leq u_i \leq R\}$, $1\leq i\leq {d}$;\\
\item {\sl step signal inputs:} $U_i^+(\rho)=\{u\in\bR^{d}:u_k=0,k<i, \rho \leq  u_i=u_{i+1}=...=u_{d} \leq R\}$, $1\leq i\leq {d}$.
\end{itemize}
To make the notation consistent with {the one used} in the general setting described in Section \ref{dynsitgoal}, we set $U_{2i-1}(\rho)=U_i^+(\rho)$ and $U_{2i}(\rho)=U_i^-(\rho)$, thus
getting $N=2{d}$
parametric families of signal inputs and we refer to signal input $u\in U_j(\rho)$ as to signal input of shape $j$ and magnitude $\geq\rho$.\\

\par {\bf Control parameters} are $R>0$ appearing in (\ref{meq2}) and  tolerance $\epsilon\in(0,1/2)$ responsible for the risks of our decision rules.
\\
\par {\bf Our goal}
is to find decision rules $\T_k$ and positive reals $\rho_{kj}\in (0,R]$, $1\leq k\leq {d}$, $1\leq j\leq N$, {such} that
\begin{itemize}
\item $\T_k$ makes a decision given observation $y^k=[y_1;...;y_k]$, and this decision, depending on $y^k$, is either ``signal conclusion,'' or ``nuisance conclusion'' (exactly one of them).
 In the case of signal conclusion at step $k$, the inference procedure is terminated at this step, otherwise we pass to time $k+1$ (when $k<{d}$) or terminate (when $k={d}$).
\item The following risk specifications are met:
\begin{itemize}
\item When the nuisance hypothesis is true (i.e.,  $u=0$), the probability of terminating with signal conclusion  {(false alarm)} somewhere on the time horizon $1,...,{d}$ is $\leq\epsilon$.
\item For every $k\leq {d}$ and every $j\leq N$, if the input $u$ underlying our observation belongs to $U_j(\rho)$ with $\rho_{kj}\leq\rho< R$, the probability of signal conclusion somewhere on the time horizon $1,...,k$ should be at least $1-\epsilon$.
\end{itemize}
\end{itemize}
\par We intend to achieve this goal with as small $\rho_{kj}$ as possible.
\par We now explain how to cast the problem we have just described as the more general detection problem
discussed in Sections \ref{dynsitgoal}, \ref{sectschemeI}, \ref{sectschemeII}.
To this end, we need to build a new observation scheme. More precisely, we have
$$
y^k=\alpha_0[\underbrace{1;...;1}_{k}]+B_k\eta+C_k\zeta+B_ku
$$
for some known matrices $B_k$ and $C_k$. Let $F_k$ be ${d}\times (k-1)$ matrix with  columns  forming an orthonormal basis of the {$(k-1)$-${d}$imensional} subspace of $\bR^k$ comprised of vectors with zero mean. Let us set
\begin{equation}\label{eq45}
z^k:=F_k^Ty^k=D_k\eta+E_k\zeta+D_ku
\end{equation}
where $(k-1)\times {d}$ matrices $D_k$ and $E_k$ are given by $D_k=F_k^TB_k$, $E_k=F_k^TC_k$, and, as is immediately seen,  have rank $k-1$. We treat $z^k$, rather than $y^k$, as an intermediate
observation at time $k$.\par
We then find a positive definite matrix $\Theta_k$ such that
$$
\Theta_k^2\succeq D_kD_k^T\ \&\ \Theta_k^2\succeq E_kE_k^T.
$$

\par\noindent
Finally, the observation $\omega^k$ at time $k$, is the vector (cf. \rf{eq233})
\begin{equation}\label{meq5}
\omega^k=\Theta_k^{-1}z^k=\underbrace{[\Theta_k^{-1}D_k]}_{A_k}u+\underbrace{\Theta_k^{-1}\left[D_k\eta+E_k\zeta\right]}_{\xi_k}.
\end{equation}
\par We have come to observation scheme \eqref{meq5} which can be handled by the change detection procedure we are about to describe
in Sections \ref{dynsitgoal}, \ref{sectschemeI}, \ref{sectschemeII}.
In particular, by Proposition \ref{propnew}, $\xi_k$ has a symmetric density $p_k(\cdot)\in\P_{\gamma}^{k-1}$ with
$
\gamma=\gamma_{\S}\star\gamma_\sigma,
$
where $\gamma_\S$ is the density of the standard  univariate Student's {$t_\nu$} distribution with $\nu$ degrees of freedom, and $\gamma_\sigma$ is the density of $\N(0,\sigma^2)$.

\subsection{Situation and goal}\label{dynsitgoal}
{In this Section our objective is to make decisions about an unknown  vector $u\in \bR^{n_u}$ representing inputs of a linear system in the situation where the information about $u$ is acquired sequentially, and the goal is to decide whether
the input observed so far is a nuisance or is  ``meaningful.'' Our modeling methodology goes back to \cite{alltogether}. Specifically,
suppose that the observation $\omega^k\in \bR^{m_k}$ available at step $k=1,...,K$, is
\be
\omega^k=A_ku+\xi_k,
\ee{eqOSt}
where $u$ is the unknown system input, $\xi_k\in \bR^{m_k}$ are random noises, and
$A_k\in \bR^{m_k\times n_u}$ are known matrices. Throughout this section we suppose that we are given}\\
\def\myQ{Z}
\begin{enumerate}
\item a family $\P$ of distributions on $\bR^{n_\xi}$ and  matrices $\myQ_k\in \bR^{m_k\times n_\xi}$ such that $\myQ_k$ is of rank $m_k$ and
\be
\xi_k=\myQ_k\xi,\,k=1,...,K,
\ee{eq:xidef} where $\xi$ obeys some (possibly unknown) distribution $p_\xi\in \P$;
\item a convex compact set $U_\inp\subset \bR^{n_u}$ of {\sl admissible inputs};
\item a convex compact {\sl nuisance set} $V_\nui\subset \inter U_\inp$ such that $0\in V_\nui$;
\item $N$ convex and closed {\sl ``activation'' sets} $W_j\subset \bR^{n_u}$ such that $0\not\in W_j$ and
\begin{equation}\label{eq2}
w\in W_j,\rho\geq1\Rightarrow \rho w\in W_j;
\end{equation}
\item $N$ nonempty convex compact ``drag sets''  $V^j_\shft$ such that $0\in V^j_\shft$.
\end{enumerate}
\par For the example of Section \ref{sec:motivatingexample}, we have
$\xi=[\eta;\zeta]$, $Z_k[\eta;\zeta]=\Theta_k^{-1}\left[D_k\eta+E_k\zeta\right]$,
$U_\inp = \{u \in \mathbb{R}^d : 0 \leq u_i \leq R \},$ $V_\nui =V^j_\shft =\{0\}$,
and
\begin{itemize}
 \item for {\sl pulse signal inputs:} $W_{2i-1}=\{u\in\bR^{d}:u_k=0,k\neq i, 1 \leq u_i \}$, $W_{2i}=-W_{2i-1}$, $1\leq i\leq {d}$;\\
\item for {\sl step signal inputs:} $W_{2i-1}=\{u\in\bR^{d}:u_k=0,k<i, 1 \leq  u_i=u_{i+1}=...=u_{d}\}$, $W_{2i}=-W_{2i-1}$, $1\leq i\leq {d}$.
\end{itemize}
We call an input $u\in V_\nui$ {\sl a nuisance}, and a vector $z$ of the form $z=\rho w$, with $w\in W_j$ and $\rho>0$, an {\sl activation of shape $j$ and magnitude $\geq\rho$.} Note that if $0<\rho'\leq\rho$ and $z$ is an activation of shape $j$ and magnitude $\geq\rho$, then, as it should be, $z$ is an activation of shape $j$ and magnitude $\geq\rho'$, due to $z=\rho'w'$, $w'=(\rho/\rho')w$, and $w'\in W_j$ due to $0<\rho'\leq\rho$ and (\ref{eq2}).
\par
We call input $u$ a {\sl signal of shape $j$ and magnitude $>\rho$ (or $\geq\rho>0$)}, if $u\in U_\inp$ and
$$
u=v+\rho' w
$$
\def\sig{{\hbox{\tiny\rm sig}}}
with $v\in V_\shft^j$, $w\in W_j$, and $\rho'>\rho$ (resp. $\rho'\geq\rho$), and denote by
$$
U_{j}(\rho)=\{u=v+\rho w: u\in U_\inp,v\in V_\shft^j,w\in W_j\}
$$
the set of all signal inputs of shape $j$ and magnitude $\geq\rho$.
\paragraph{Our goal} is to decide via observations $\omega^k$, $k\leq K$, on the nuisance hypothesis ``the input $u$ to (\ref{eqOSt}) is a nuisance'' (i.e., $u\in V_\nui$) vs. the signal alternative
``the input $u$ to (\ref{eqOSt}) is a signal of some shape and (positive) magnitude.''
Note that in this case, under the signal alternative the input $u$ belongs to a nonconvex set expressed as a union of convex sets.
More precisely, we assume that we are given tolerances
$$
\{\epsilon_k\in(0,1/2),\, 1\leq k\leq K\}, \;\{\epsilon_{kj}\in(0,1/2),1\leq k\leq K,1\leq j\leq N\}
$$
and want to design a sequence of decision rules $\{\T_k:1\leq k\leq K\}$ along with thresholds $\rho_{kj}>0$, $1\leq k\leq K$, $1\leq j\leq N$ with the following properties: for every $k\leq K$,
\begin{itemize}
\item
rule $\T_k$ makes a decision based on observation $\omega^k$, and this decision is either to accept the null hypothesis or to accept the signal one (but not both);
\item
if the input is a nuisance, the probability for $\T_k$ to accept the signal alternative is $\leq \epsilon_k$;
\item
if, for some $j$, the input is a signal of shape $j$ and magnitude $>\rho_{kj}$, the probability for $\T_k$ to accept the nuisance hypothesis is at most $\epsilon_{kj}$.
\end{itemize}
Given tolerances $\epsilon_k$, $\epsilon_{kj}$, we would like to achieve the outlined goal with thresholds $\rho_{kj}$ as small as possible.
\subsection{Tests $\T_k$, Scheme I.}\label{sectschemeI}
\subsubsection{Assumptions on the distribution of noise $\xi$}
Throughout Section \ref{sectschemeI}, we make the following assumption on the family $\P$ of  probability densities of random disturbance $\xi$ in \rf{eq:xidef}:
\begin{assump}\label{ass:1}
For every $k\leq K$ we can point out a parametric family
$\H_k=\{\eta^k_\delta(\cdot):\delta\geq0\}$ of potentials and a continuous nonincreasing function $\R_k(\delta):\,\bR_+\to (0,1]$
such that $\R_k(0)=1$ and
\begin{equation}\label{teq1}
\Risk_\delta(\eta^k_\delta|\P^k)\leq \R_k(\delta)\,\,\forall \delta\geq0,
\end{equation}
where $\P^k$ is the family of probability densities of random vectors $\xi_k=\myQ_k\xi$ with $\xi\sim p_\xi\in\P$.
\end{assump}
Note that Assumption \ref{ass:1} indeed holds in the situations of our primary interest. Specifically, assume that $\myQ_k\myQ_k^T\preceq I_{m_k}$.\footnote{this always can be achieved by appropriate scaling of observations (\ref{eqOSt}).} Then
\begin{itemize}
\item
when $\P=\P_\sG^{n_\xi}$ is the family of sub-Gaussian distributions on $\bR^{n_\xi}$, with parameters $0,I_{n_\xi}$, then
$\P^k$ belongs to the family $\P_\sG^{m_k}$ of sub-Gaussian, with parameters $0,I_{m_k}$, distributions on $\bR^{m_k}$. By Proposition \ref{prop33}, (\ref{teq1}) is ensured by the choice
    \[
    \H_k=\{\eta^k_\delta(s)=\delta s,\delta\geq0\},\;\; \R_k(\delta)={\rm e}^{-\delta^2/2};
\]
\item
when $\P=\P_\gamma$ is a sub-spherical family of distributions on $\bR^{n_\xi}$, \rf{eq:xidef} combines with Proposition \ref{prop22} to ensure that the family $\P^k$ is contained in the sub-spherical family of distributions $\P^{m_k}_\gamma$. Invoking the example of the step potential from Section \ref{sectsfpot}, (\ref{teq1}) is ensured by setting
\[
\begin{array}{c}
    \H_k=\left\{\eta^k_\delta(s)=\half \ln\left(
{1-\myepsilonstar(\delta)\over \myepsilonstar(\delta)}
\right)\sign(s),\delta\geq0\right\}, \;\;\R_k(\delta)=2\sqrt{\myepsilonstar(\delta)(1-\myepsilonstar(\delta))},\\
\myepsilonstar(\delta)=\int\limits_\delta^\infty \gamma(s)ds.\\
\end{array}
\]
\end{itemize}
\subsubsection{Building decision rules}
For every $k\leq K$ and every $j\leq N$, consider the parametric convex optimization problem
$$
\Opt_{kj}(\rho)=\min_{v',v,w}\left\{\half\|A_k(v'-[v+\rho w])\|_2: v'\in V_\nui,v\in V^j_\shft,w\in W_j,v+\rho w\in U_\inp\right\}\eqno{(\PP_{kj}[\rho])}
$$
where the parameter $\rho$ is positive.
\par From our assumptions on $V_\nui$, $U_\inp$, $V^j_\shft$  and $W_j$ it immediately follows that
\begin{enumerate}
\item the set $\Delta_{j}$ of those $\rho>0$ for which $(\PP_{kj}[\rho])$ is feasible, is a half-open segment $(0,R_{j}]$ with $0<R_{j}<\infty$, and
\item when $\rho\in\Delta_{j}$, problem $(\PP_{kj}[\rho])$ is solvable, and $\Opt_{kj}(\rho)$ is a  real-valued continuous nondecreasing function on $\Delta_{j}$
such that $\lim_{\rho\to+0} \Opt_{kj}(\rho)=0$.
\end{enumerate}
Given $k\leq K$, we specify the test $\T_k$ as follows.
\begin{enumerate}
\item We compute the quantities $R_j$.
\footnote{The simplest way to identify $R_j$ is to run bisection in $\rho$ on a large initial range of $\rho$ in order to find the largest $\rho$ for which $(\PP_{kj}[\rho])$ is feasible.}
\item We select somehow a (perhaps, empty) set $J_k\subset\{1,2,...,N\}$ and reals $\rho_{kj}\in(0,R_j]$, $j\in J_k$, in such a way that
\begin{enumerate}
\item $\rho_{kj}=R_j$ when $j\not\in J_k$,
\item we have
\begin{equation}\label{eq13}
\sum_{j\in J_k}\epsilon_{kj}^{-1}\R_k^2(\Opt_{kj}(\rho_{kj}))\leq \epsilon_k.
\end{equation}
\end{enumerate}
Note that (\ref{eq13}) implies that $\Opt_{kj}(\rho_{kj})>0$ for $j\in J_k$,  otherwise the left hand side in (\ref{eq13}) is at least 1 due to $\R_k(0)=1$, and we have assumed that $\epsilon_k<1/2$.
\item
If $J_k=\emptyset$, $\T_k$ accepts the nuisance hypothesis. When $J_k\neq\emptyset$, we act as follows:
\begin{enumerate}
\item for $j\in J_k$, we solve the optimization problem $(\PP_{kj}[\rho_{kj}])$ and denote by $(v_{kj}^\prime;v_{kj},w_{kj})$ an optimal solution to the problem. We set
\begin{equation}\label{eq9}
\begin{array}{c}
\delta_{kj}=\Opt_{kj}(\rho_{kj}),\;\;\;u_{kj}=v_{kj}+\rho_{kj}w_{kj},\;\;\;h_{kj}={A_k[v_{kj}^\prime-u_{kj}]\over \|A_k[v_{kj}^\prime-u_{kj}]\|_2},\\
c_{kj}=\half h_{k j}^T A_k[v_{kj}^\prime+u_{kj}],\;\;\;\phi_{kj}(\omega^k)=\eta^k_{\delta_{kj}}(h_{kj}^T \omega^k-c_{kj}),
\end{array}
\end{equation}
with $\eta^k_\delta(\cdot)$ given by Assumption {\ref{ass:1}},
and define $\alpha_{kj}$ from the relation
\begin{equation}\label{eq909}
{\rm e}^{\alpha_{kj}}=\epsilon_{kj}/\R_k(\delta_{kj}).
\end{equation}
\item Finally, given observation $\omega^k$, the rule $\T_k$ accepts the nuisance hypothesis if $\phi_{kj}(\omega^k)+\alpha_{kj}\geq 0$ for all $j\in J_k$, and accepts the signal hypothesis otherwise.
\end{enumerate}
\end{enumerate}
\subsubsection{Performance analysis}
We now check that the just defined decision rules meet the goal stated in Section \ref{dynsitgoal}.
Let us fix $k\leq K$ and $j\in J_k$. Given feasible input $u$, let $P^k_u$ be the distribution of observation $\omega^k$, the input being $u$. By construction and due to Proposition \ref{prop0}, applied
to observations $\omega^k=A_k u +\xi_k$ and  $X_1= A_k V_\nui$, $X_2=A_kU_j(\rho_{kj})$, we have
\begin{equation}\label{eq11}
\begin{array}{lcll}
\int{\rm e}^{-\phi_{kj}(\omega^k)}P^k_u(d\omega^k)&\leq&\R_k(\delta_{kj}),&\hbox{when $u$ is a nuisance},\\
\int{\rm e}^{\phi_{kj}(\omega^k)}P^k_u(d\omega^k)&\leq&\R_k(\delta_{kj}),&\hbox{when $u$ is a signal of shape $j$ and magnitude $\geq\rho_{kj}$}.
\end{array}
\end{equation}
\par
Assume first that the input $u$ is a nuisance, and let us upper-bound the $P^k_u$-probability of rejecting the nuisance hypothesis at step $k$. This takes place only when $\phi_{kj}(\omega^k)+\alpha_{kj}<0$ for some $j\in J_k$, and,
by the first relation in (\ref{eq11}), for a given $j\in J_k$ the $P^k_u$-probability of this event  is at most ${\rm e}^{-\alpha_{kj}}\R_k(\delta_{kj})=\R_k^2(\delta_{kj})/\epsilon_{kj}$ where the last equality is due to (\ref{eq909}). As a result, the $P^k_u$-probability to reject the nuisance hypothesis at step $k$  is at most
$$
\sum_{j\in J_k}\epsilon_{kj}^{-1}\R_k^2(\delta_{kj})=\sum_{j\in J_k}\epsilon_{kj}^{-1}\R_k^2(\Opt_{kj}(\rho_{kj}))\leq \epsilon_k,
$$
where the concluding inequality is due to (\ref{eq13}).
\par
Now, let the input $u$ be a signal of shape $j$ and magnitude $>\rho_{kj}$. Due to $\rho_{kj}=R_j$ for $j\not\in J_k$, we have $j\in J_k$. Let us upper-bound the $P^k_u$-probability of the nuisance conclusion at step $k$. By construction, the latter occurs only when $\phi_{kj}(\omega^k)+\alpha_{kj}\geq0$, and,
by the second relation in (\ref{eq11}), this can happen with $P^k_u$-probability  at most ${\rm e}^{\alpha_{kj}}\R_k(\delta_{kj})=\epsilon_{kj}$, where the concluding equality is due to (\ref{eq909}). \par
The bottom line is that the probability of false alarm at step $k$ (a signal conclusion when the input is nuisance) is $\leq \epsilon_k$, and the probability of $\rho_{kj}$-miss (the probability to make a nuisance conclusion when the  input is a signal of shape $j$ and magnitude $>\rho_{kj}$) is at most $\epsilon_{kj}$.
\subsection{Decision rules $\T_k$, Scheme II}\label{sectschemeII}
Throughout Section \ref{sectschemeII}, we make the following
\begin{quote}
{\bf Assumption A:}  The
family $\P$ of probability densities of the disturbance $\xi$ in \rf{eq:xidef} is such that for all $k\leq K$, probability densities of noises $\xi_k$ belong to sub-spherical families $\P^{m_k}_\gamma$ with some common $\gamma$.
 \end{quote}
Note this assumption takes place, e.g., when  $\P$ is a sub-spherical family $\P^{n_\xi}_\gamma$ and $\myQ_k\myQ_k^T\preceq I_{m_k}$, $k\leq K$, see Proposition \ref{prop22}.\par
We put (cf. (\ref{weq14}))
\begin{equation}\label{yeq1}
\myepsilonstar(\delta| \gamma)=\int\limits_\delta^\infty\gamma(s)ds,
\end{equation}
\subsubsection{Building the decision rules}
Given $k\leq K$, we specify $\T_k$ as follows.
\begin{enumerate}
\item Exactly as in Scheme I, for every $j\leq N$, we consider the parametric convex optimization problem ($\PP_{kj}[\rho]$)
with positive $\rho$, and compute $R_j$, the largest  $\rho$ for which the problem is feasible.
\item We select somehow a (perhaps, empty) set $J_k\subset\{1,2,...,N\}$ and reals $\rho_{kj}\in(0,R_j]$, $1\leq j\leq N$, and $\alpha_{\chi kj}$, $\chi=1,2$, $j\in J_k$ in such a way that
 $\rho_{kj}=R_j$ when $j\not\in J_k$,
and we have
\begin{equation}\label{eq23}
\begin{array}{lrcl}
(a)~~~~~~~~~~~~& \myepsilonstar(\alpha_{2kj}| \gamma)&\leq &\epsilon_{kj},\;\;\forall j\in J_k,\\
(b)~~~~~~~~~~~~& \sum\limits_{j\in J_k}\myepsilonstar(\alpha_{1kj}| \gamma)&\leq& \epsilon_k,\\
(c)~~~~~~~~~&\alpha_{1kj}+\alpha_{2kj}&\leq& 2\delta_{kj},\;\;\delta_{kj}=\Opt_{kj}(\rho_{kj}),\;\;\forall j\in J_k.\\
\end{array}
\end{equation}
Note that for $j\in J_k$ we have $\alpha_{1kj}>0$ and $\alpha_{2kj}>0$ (by (\ref{eq23}.$a$-$b$) combined with $\epsilon_{kj}<1/2$, $\epsilon_k<1/2$; recall that $\myepsilonstar(s| \gamma)\geq1/2$ when $s\leq0$). As a result,  (\ref{eq23}.$c$) implies that $\delta_{kj}>0$ whenever $j\in J_k$.
\item
If $J_k=\emptyset$, $\T_k$ accepts the nuisance hypothesis. When $J_k\neq\emptyset$, we act as follows:
\begin{enumerate}
\item for $j\in J_k$, we solve the optimization problem $(\PP_{kj}[\rho_{kj}])$ and denote by $(v_{kj}^\prime;v_{kj},w_{kj})$ an optimal solution to the problem. Similarly to (\ref{eq9}), we set
{
\begin{equation}\label{eq29}
\begin{array}{c}
u_{kj}=v_{kj}+\rho_{kj}w_{kj},\;\;\;
h_{kj}={A_k[v_{kj}^\prime-u_{kj}]\over \|A_k[v_{kj}^\prime-u_{kj}]\|_2},\\
c_{kj}=\half h_{k j}^T A_k[v_{kj}^\prime+u_{kj}],\;\;\;\phi_{kj}(\omega^k)=h_{kj}^T \omega^k- c_{kj}.
\end{array}
\end{equation}
}
\item Finally, given observation $\omega^k$, the rule $\T_k$ accepts the nuisance hypothesis if $\phi_{kj}(\omega_k)\geq\half {(\alpha_{2kj}-\alpha_{1kj})}$ for all $j\in J_k$, and accepts the signal hypothesis otherwise.
\end{enumerate}
\end{enumerate}
\subsubsection{Performance analysis}\label{sectperfan}
Let $k\leq K$ be fixed, and let $P^k_u$ be the probability distribution of observation $\omega^k$, the input being $u$. Taking into account (\ref{eq23}.$c$) and applying Proposition \ref{prop1} with $\alpha_\chi=\alpha_{\chi kj}$, $\delta=\delta_{kj}$ and
$X_1=A_k V_\nui$, $X_2=A_kU_j(\rho_{kj})$, we obtain
\begin{equation}\label{eq21}
\begin{array}{llcll}
(a)~~~~&P^k_u\left\{\phi_{kj}(\omega^k)<\half{(\alpha_{2kj}-\alpha_{1kj})}\right\}&\leq& \myepsilonstar(\alpha_{1kj}| \gamma),&
\hbox{if $u\in V_\nui$,}\\
(b)~~~~&P^k_u\left\{\phi_{kj}(\omega^k)\geq\half{(\alpha_{2kj}-\alpha_{1kj})}\right\}&\leq&\myepsilonstar(\alpha_{2kj}| \gamma),&\hbox{if $u\in U_j(\rho_{kj}).$}\\
\end{array}
\end{equation}
\par
Assume first that the input $u$ is a nuisance, and let us upper-bound the $P^k_u$-probability of rejecting the nuisance hypothesis at step $k$. This rejection implies that
$\phi_{kj}(\omega^k)<\half{(\alpha_{2kj}-\alpha_{1kj})}$ for some $j\in J_k$, and by (\ref{eq21}.a) the $P^k_u$-probability of this event for a given $j\in J_k$ is at most
$\myepsilonstar(\alpha_{1kj}| \gamma)$. As a result, the $P^k_u$-probability of signal conclusion at step $k$ when $u$ is a nuisance is at most
$$
\sum_{j\in J_k}\myepsilonstar(\alpha_{1kj}| \gamma)\leq\epsilon_k,
$$
where the  inequality is due to (\ref{eq23}.$b$).
\par
Now let the input $u$ be a signal of shape $j$ and magnitude $>\rho_{kj}$, i.e., $u \in U_j(\rho)$ with $\rho>\rho_{k j}$. Due to $\rho_{kj}=R_j$ for $j\notin J_k$ {this means that} $j\in J_k$. Let us upper-bound the $P^k_u$-probability of a nuisance conclusion at step $k$. The nuisance hypothesis is not rejected only when $\phi_{kj}(\omega^k)\geq\half {(\alpha_{2kj}-\alpha_{1kj})}$, and by {(\ref{eq21}.b)} the {$P^k_u$-probability of the latter event does not exceed}
$\myepsilonstar(\alpha_{2kj}| \gamma)\leq \epsilon_{kj}$ (recall that by definition of sets $U_j(\rho)$, $u \in U_j(\rho)$ with
$\rho>\rho_{k j}$ implies that $u \in U_j( \rho_{k j})$), where the concluding inequality is due to (\ref{eq23}.$a$). \par
The bottom line is that the $P^k_u$-probability of false alarm at step $k$
(rejecting the nuisance hypothesis when it is true) is $\leq \epsilon_k$, and the $P^k_u$-probability of $\rho_{kj}$-miss (making a nuisance conclusion when the  input is a signal of shape $j$ and magnitude $>\rho_{kj}$) is at most $\epsilon_{kj}$.

\section{Application: change detection in linear dynamical system}\label{changepoint}
\subsection{Problem statement}\label{sectchpsit}
We {consider the} change detection problem as follows.
\begin{enumerate}
\item We are given a discrete time linear time invariant system
\begin{equation}\label{eq1}
\begin{array}{rcl}
x_{t}&=&P_tx_{t-1}+Q_tu+R_t\xi,\\
y_t&=&C_tx+D_tu+S_t \xi,\;\;t=1,...,d,
\end{array}
\end{equation}
where
\begin{itemize}
\item
$x=[x_0;x_1;...;x_d ]$, $x_t\in \bR^{n_x}$ is the state trajectory, $u\in \bR^{n_u}$ is the input,
$y_t\in \bR^{n_y}$ is the output at time $t$,
\item
$\xi\in\bR^{n_\xi}$ is a random disturbance with (unknown) distribution $P$,
\item $P_t$,...,$S_t$, $1\leq t\leq d$, are known matrices of appropriate sizes.
\end{itemize}
Given $\tau\leq d$, we set $y^\tau=[y_1;...;y_\tau]$.
We denote
$
E_x=\underbrace{\bR^{n_x}\times...\times \bR^{n_x}}_{d}
$
and similarly for $E_u$.
\item We are given sets $U_\inp\subset E_u$ (admissible inputs), $V_\nui\subset E_u$ (nuisances), $V^j_\shft$
(drags), $W_j\subset E_u$ (activations of shape $j$ and magnitude $\geq1$), $1\leq j\leq N$, meeting the requirements of Section \ref{dynsitgoal}. These sets, exactly as in Section \ref{dynsitgoal},  give rise to the notions of
    admissible and  nuisance inputs to (\ref{eq1}), same as signal inputs of shape $j$ and magnitude $\geq \rho$.
\item The distribution $P$ of disturbance $\xi$ is known to belong to a given family  $\P$ of probability densities on $\bR^{n_\xi}$. We are also given
tolerances $\epsilon_t\in(0,1/2)$, $\epsilon_{tj}\in (0,1/2)$, $1\leq t\leq d$, $1\leq j\leq N$.
\end{enumerate}
We {acquire} observations $y_\tau$ one by one, so that at time $t$ we have at our disposal the observation $y^t=[y_1;...;y_t]$. Our objective is to design  tests $\T_t$ and thresholds $\rho_{tj}>0$, $1\leq t\leq d$, $1\leq j\leq N$,
meeting requirements completely similar to {those} from Section \ref{dynsitgoal}:
\begin{itemize}
\item given observation $y^t$, the test $\T_t$ should make either a nuisance, or a signal conclusion (but not both);
\item if the input to (\ref{eq1}) is a nuisance, the probability of the non-nuisance conclusion at time $t$ should be at most $\epsilon_t$, and if, for some $j$, the input is a signal of shape $j\leq N$ and magnitude $>\rho_{tj}$, the probability of the nuisance conclusion at time $t$ should be at most $\epsilon_{tj}$.\par
\end{itemize}
Both these requirements should be satisfied for every $t\leq d$, and we would like to meet them with as small thresholds $\rho_{tj}$ as possible.
\par We are about to demonstrate that the just outlined change detection problem  can be handled via the techniques developed in Sections  \ref{sectschemeI} and \ref{sectschemeII}. Basically all we need to this end is to convert our observation scheme into the one
considered in Section \ref{sectchp}, and this is what we are about to do next.
\subsection{Building the observation scheme}
Given an input $u$, a noise $\xi$ and $t\leq d$, the observation $y^t$ is not uniquely determined by the input and the noise; it is also affected by the initial state $x_0$ of the system. To {get rid of the influence of the initial condition} we act as follows.
\begin{enumerate}
\item We denote by $F^t$ the linear subspace of $E_y^t=\{[y_1;...;y_t]\in\bR^{tn_y}\}$  comprised of all outputs $y^t=[y_1;...;y_t]$ which {\sl in the noiseless case $\xi=0$} stem from zero input and some
initial state of the system, build an orthonormal basis of the orthogonal complement of $F^t$ in $E_y^t$ and make the vectors of this basis the rows of a matrix, thus arriving at a $\mu_t\times(tn_y)$ matrix $M_t$. Note that $\mu_t$ is a nondecreasing function of $t$.
\item
We {set}
$
z^t=M_t[y_1;...;y_t],
$
where $y_\tau$ is given by (\ref{eq1}).
\par
It may happen that $\mu_t=0$ for some $t$. In this case, our decision rule $\T_t$ by construction accepts the nuisance hypothesis, so that nontrivial decision rules will be associated only with those time instants $t$ for which $\mu_t\geq 1$. {Let $t=\kappa+1$ be the first instant such that the corresponding $\mu_t\geq 1$. Note that time} instants $t$ with $\mu_t\geq 1$ form the {final} segment $\{\kappa+1,\kappa+2,...,\kappa+K=d\}$ of $1,...,d$ (recall that $\mu_t$ is nondecreasing in $t$). We set $m_k=\mu_{\kappa+k}$, $1\leq k\leq K$.
By construction,
\begin{equation}\label{eq3}
z^{\kappa+k}=\overline{A}_{k}u+\overline{B}_k\xi,
\end{equation}
with some $m_k\times n_u$ matrix $\overline{A}_k$ and $m_k\times n_\xi$ matrix $\overline{B}_k$ readily given by our data.
\end{enumerate}
\par
 From now on, we
assume that\footnote{Observe that when $K=0$ our approach results in trivial tests always accepting the nuisance hypothesis -- in this case the input-related component in the observations is fully masked by the influence of initial condition $x_0$.}
 $K>0$ and that $\overline{B}_k$ are of full row rank, i.e. $\rank(\overline{B}_k )=m_k$, $1\leq k\leq K$.
Finally, we select somehow invertible $m_k\times m_k$ matrices $L_k$ and pass from observations (\ref{eq3}) to observations
\begin{equation}\label{eq233}
\begin{array}{c}
\omega^k=A_ku+\xi_k,\hbox{\ where\ }\omega^k=L_k z^{\kappa + k},\,
A_k=L_k\overline{A}_k,\;\;\xi_k=\myQ_k\xi,\;\;\myQ_k=L_k\overline{B}_k.
\end{array}
\end{equation}
For example, we can set $L_k=(\overline{B}_k\overline{B}_k^T)^{-1/2}$, thus ensuring that $\myQ_k\myQ_k^T=I_{m_k}$.
\par Notice that observations (\ref{eq233}) meet the requirements imposed in Section \ref{dynsitgoal}. As a result, we find ourselves in the situation considered in Section \ref{dynsitgoal} and therefore can apply to the change detection problem in question the machinery developed in Sections \ref{sectschemeI} and \ref{sectschemeII}.
\subsection{{Illustration: detecting changes in the trend of a simple time series}}\label{numericalill}

We consider the example of Section \ref{sec:motivatingexample} and apply, with  minor modifications, the construction outlined in Section \ref{sectschemeII}.

{\subsubsection{Constructing decision rules}}\label{sec:constructing}

To attain our  goal (see Section \ref{sec:motivatingexample}) we act as follows. The rule $\T_1$ is trivial -- it always accepts the nuisance hypothesis. To describe how $\T_k$, $k>1$, is built, let us fix $k\in\{2,...,{d}\}$.
\begin{enumerate}
\item {\sl Building observation {$\omega^k$}.}
Setting $\xi=[\eta;\zeta]$ and
$Z_k[\eta;\zeta]=\Theta_k^{-1}\left[D_k\eta+E_k\zeta\right]$, our observations (\ref{meq5}) are as required in (\ref{eqOSt}), and we meet Assumption A.

In our implementation, we use $\Theta_k=\Omega_k^{1/2}$, where $\Omega_k$ is the minimum trace matrix satisfying $\Omega_k\succeq D_kD_k^T$, $\Omega_k\succeq E_kE_k^T$.

\par We now apply the decision rules of Scheme II from Section \ref{sectschemeII} to observation \rf{meq5}. To define parameters
$J_k, \rho_{k j}, \alpha_{1 k j}$ and  $\alpha_{2 k j}$ we proceed as follows.
\item We set
\begin{equation}\label{meq75}
\begin{array}{c}
\widehat{J}_k=\{1,...,2k\};\,\,\widehat{\epsilon}={\epsilon\over {d}({d+1})-2};\,\,
\epsilon_k=2k\widehat{\epsilon};\,\,
\epsilon_{kj}=\epsilon,\,j\in \widehat{J}_k;\,\,
\rho_{kj}=R,\,j\not\in \widehat{J}_k.\\
\end{array}
\end{equation}
For $j\in \widehat{J}_k$, we specify $\alpha_{1kj}$, $\alpha_{2kj}$, $\delta_{kj}$  by the relations
$$
\int\limits_{\alpha_{1kj}}^\infty \gamma(s)ds=\widehat{\epsilon},\,\,\int\limits_{\alpha_{2kj}}^\infty \gamma(s)ds=\epsilon_{kj}=\epsilon,\,\,\delta_{kj}=\half[\alpha_{1kj}+\alpha_{2kj}].
$$
 Note that, by construction, setting $\epsilon_1 =0$ we have
$$
\forall k: \sum_{j\in\widehat{J}_k}\,\int\limits_{\alpha_{1kj}}^\infty\gamma(s)ds\leq\epsilon_k,\;\;\sum_{k=1}^{d}\epsilon_k=\epsilon\; \mbox{and}\;\; \alpha_{1kj}+\alpha_{2kj}=2\delta_{kj},j\in \widehat{J}_k,
$$
cf. (\ref{eq23}).

\item For $j\in \widehat{J}_k$, we consider the convex optimization problem
(cf. $(\PP_{kj}[\rho])$). $\Opt_{kj}(\rho)$ clearly is continuous and nonincreasing in $\rho>0$ and $\lim_{\rho\to+0}\Opt(\rho)=0$. When $\Opt_{kj}(R)\leq\delta_{kj}$ 
we set $\rho_{kj}=R$. Otherwise, 
we find  the smallest $\rho=\rho_{kj}$ such that  $\Opt_{kj}(\rho)\geq\delta_{kj}$; note that in {the latter case} we have  $\rho_{kj}\in(0,R)$ and $\Opt_{kj}(\rho_{kj})=\delta_{kj}$. We have specified $\rho_{kj}$ for all $j\in\widehat{J}_k$.
\item We set $J_k=\{j\in \widehat{J}_k:\rho_{kj}<R\}$ thus ensuring that $\rho_{kj}<R$ if and only if $j\in J_k$.\footnote{Note that, typically, $|J_k|<|\widehat{J}_k|$. Thus, one can easily improve the estimation procedure by better accounting for the ``remaining at step $k$'' part of
false alarm probability. A simple ``dynamic'' management of false alarm probabilities of tests is implemented in the numerical experiments described in the next section.
The detailed construction is presented in the Online complement of the paper available at {\url{http://arxiv.org/abs/1705.07196}}.}
\item
Same as in  Section \ref{sectschemeII}, for $j\in J_k$ we denote by $u_{kj}$ an optimal solution  to problem $(\PP_{kj}[\rho])$
with $\rho=\rho_{kj}$ and set ( cf. (\ref{eq29}))
    $$
    h_{kj}=-{A_ku_{kj}\over\|A_ku_{kj}\|_2}=-{A_ku_{kj}\over 2\delta_{kj}},\;\;c_{kj}=\half h_{kj}^TA_ku_{kj},\;\; \phi_{kj}(\omega)=h_{kj}^T\omega-c_{kj}.
    $$

\item Finally, given observation $\omega^k$, our rule $\T_k$ makes the nuisance conclusion if and only if $\phi_{kj}(\omega_k)\geq\half[\alpha_{2kj}-\alpha_{1kj}]$ for all $j\in J_k$, and makes signal conclusion otherwise, cf. Section \ref{sectschemeII}.
Invoking the results of Section \ref{sectperfan}, the decision rules we have built do satisfy risk specifications of Section \ref{sec:motivatingexample}.
\end{enumerate}

{\subsubsection{Quantifying conservatism: performance indexes}}
Let us pass from intermediate observations (\ref{eq45}) to observations
\begin{equation}\label{meq47}
w^k=[D_kD_k^T]^{-1/2}z^k=\underbrace{Q_k\eta+S_k\zeta}_{\lambda^k}+Q_ku,
\end{equation}
where $Q_k=[D_kD_k^T]^{-1/2}D_k$ satisfies $Q_kQ_k^T=I_{k-1}$ and $S_k=[D_kD_k^T]^{-1/2}E_k$. {Since $Q_k$ has orthonormal rows,} specifying the distribution of $\eta$ to be {multivariate} Student's $t_{{d}}(\nu,I_{d})$ distribution on $\bR^{d}$ with $\nu$ degrees of freedom and unit matrix parameter, the distribution $p(\cdot)$ of the random variable $Q_k\eta$ will be {multivariate} Student's {$t_{k-1}(\nu,I_{k-1})$ distribution}. Now, let $\theta$ be the {smallest nonvanishing} singular value of $S_k$ (or, equivalently, $\theta^2$ is the smallest eigenvalue of $S_kS_k^T$). Clearly, we can specify the covariance matrix $\Sigma$ of zero mean ${d}$-dimensional Gaussian random vector  $\zeta$ to satisfy $\Sigma\preceq \sigma^2I_{d}$ and to be such that the covariance matrix of $S_k\zeta$ is $\theta^2\sigma^2I_{k-1}$.
We conclude that  we can {point out} distributions of $\eta$ and $\zeta$, {satisfying specifications of the model} (\ref{meq1}), {and such} that the random noise $\lambda^k$ in (\ref{meq47}) will be the sum of two independent {zero-mean} random vectors, one with
$(k-1)$-dimensional Student distribution $t_{k-1}(\nu, I_{k-1})$, and {the other} -- Gaussian, with covariance matrix $(\theta\sigma)^2I_{k-1}$. As it is immediately seen, $\lambda^k$ has a probability density $p(\cdot)$ of the form $f(\|\cdot\|_2)$ with
nonincreasing $f$ and
whenever $e\in\bR^{k-1}$ is a unit vector, the probability density   $\gamma_k(\cdot)$ of the scalar random variable $e^T\lambda^k$ is
$$
\gamma_k=\gamma_\S\star\gamma_{\theta\sigma}.
$$
Here, as above, $\gamma_\S$ is the density of the  univariate Student's {$t_\nu$} distribution, and $\gamma_{\theta_k\sigma}$ is the density of $\N(0,(\theta\sigma)^2)$.
\par
Now let $V_1,V_2$ be two closed convex sets in the space $\bR^{d}$ of inputs such that the sets $Q_kV_\chi$, $\chi=1,2$, are closed, and one of these two sets is bounded, and let $(u_*^1,u_*^2)$ be an optimal solution to the convex optimization problem
\begin{equation}\label{meq48}
\delta=\min_{u^1\in V_1,u^2\in V_2}\half\|Q_k[u^1-u^2]\|_2.
\end{equation}
By Remark \ref{rem1}, no test based on observation (\ref{meq47}) can decide on two simple hypotheses $u=u^1_*$, $u=u^2_*$ with risk $<\int\limits_{\delta}^\infty\gamma_k(s)ds$, implying that the
same lower risk bound also holds true for all tests utilizing observations $y^k$ rather than ${w}^k$.
\par
Now let $V_1=\{0\}$ and $V_2=U_j(\rho)$, so that $\delta$ as defined in (\ref{meq48}) becomes a (clearly, continuous and nonincreasing when $\rho>0$) function $\delta_{kj}(\rho)$ of $\rho$. Let us define $\rho^*_{kj}$ as follows: if
$
\int_{\delta_{kj}(R)}^\infty\gamma_k(s)ds>\epsilon,
$
we set $\rho^*_{kj}=R$, otherwise $\rho^*_{kj}\in(0,R]$ is the smallest $\rho>0$ such that
$
\int_{\delta_{kj}(\rho)}^\infty\gamma_k(s)ds\leq\epsilon.
$
By construction, for every $\rho\in(0,\rho^*_{kj})$ there is no test which, given an observation $y^k$, would decide with risk $\leq\epsilon$ on the hypothesis ``{$u=0$}'' vs. the alternative ``$u$ is a signal {from} $U_j(\rho)$ with $\rho>0$.'' It is natural to quantify the conservatism of our decision rules $\T_k$  by the {\sl performance indexes} $\rho_{kj}/\rho^*_{kj}$; the less are these indexes, the less is the conservatism.\\

\subsubsection{Numerical results}
\begin{figure}
$$\begin{array}{cc}
\includegraphics[scale=0.64]{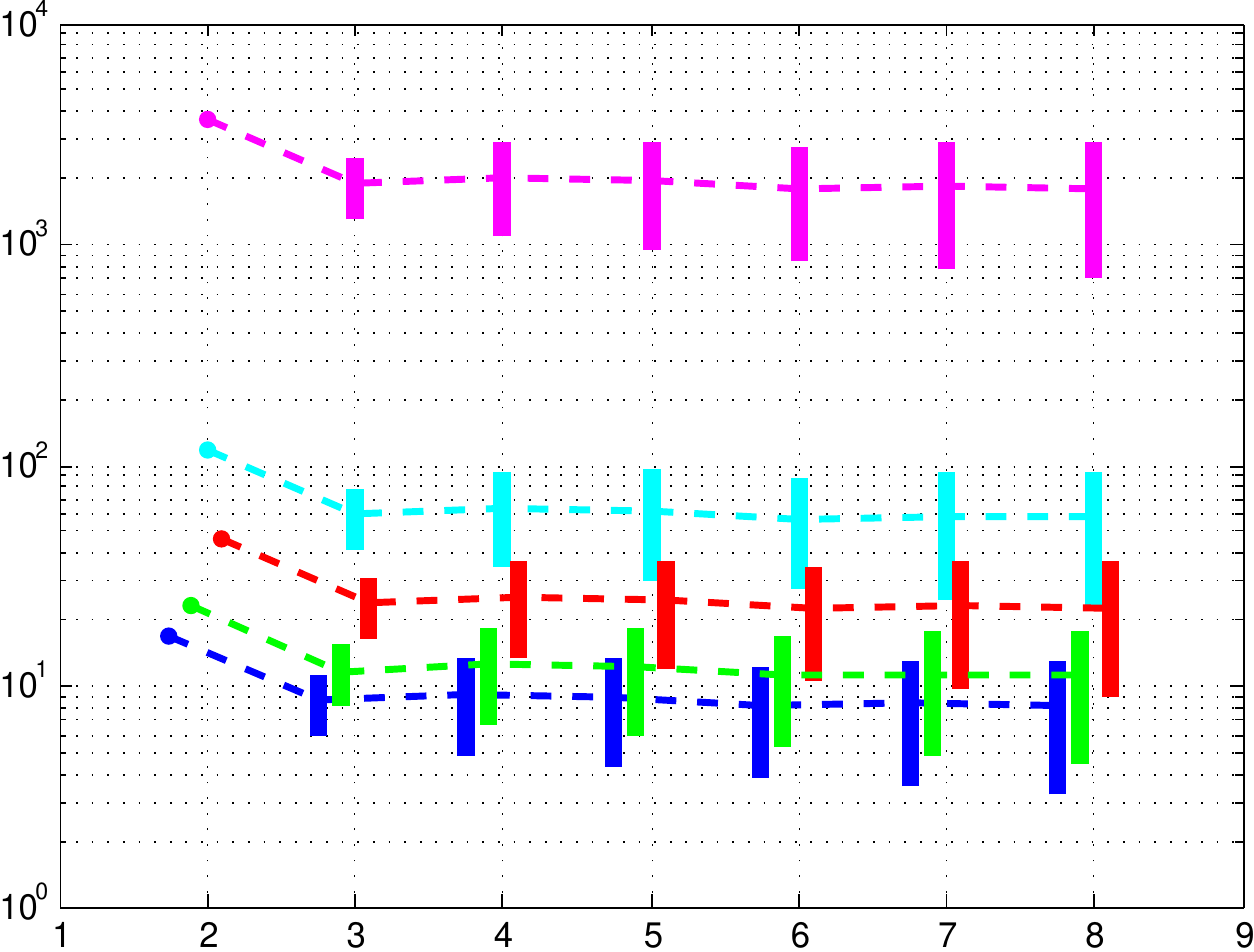}&
\includegraphics[scale=0.60]{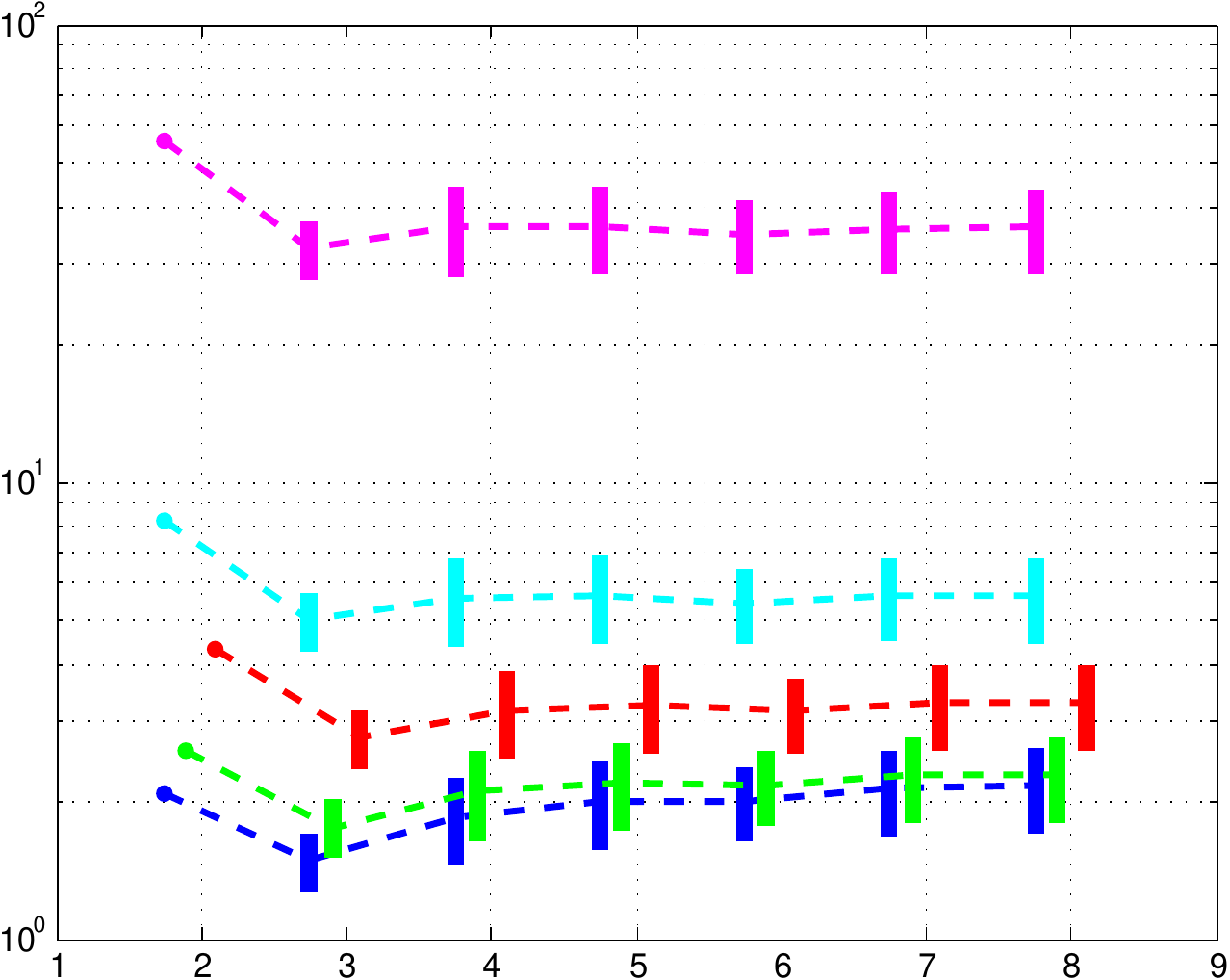}\\
 \multicolumn{2}{c}{\hbox{\small Signal input: step, $2\leq k\leq 8$, $1\leq i\leq k$}}\\
 & \\
\includegraphics[scale=0.64]{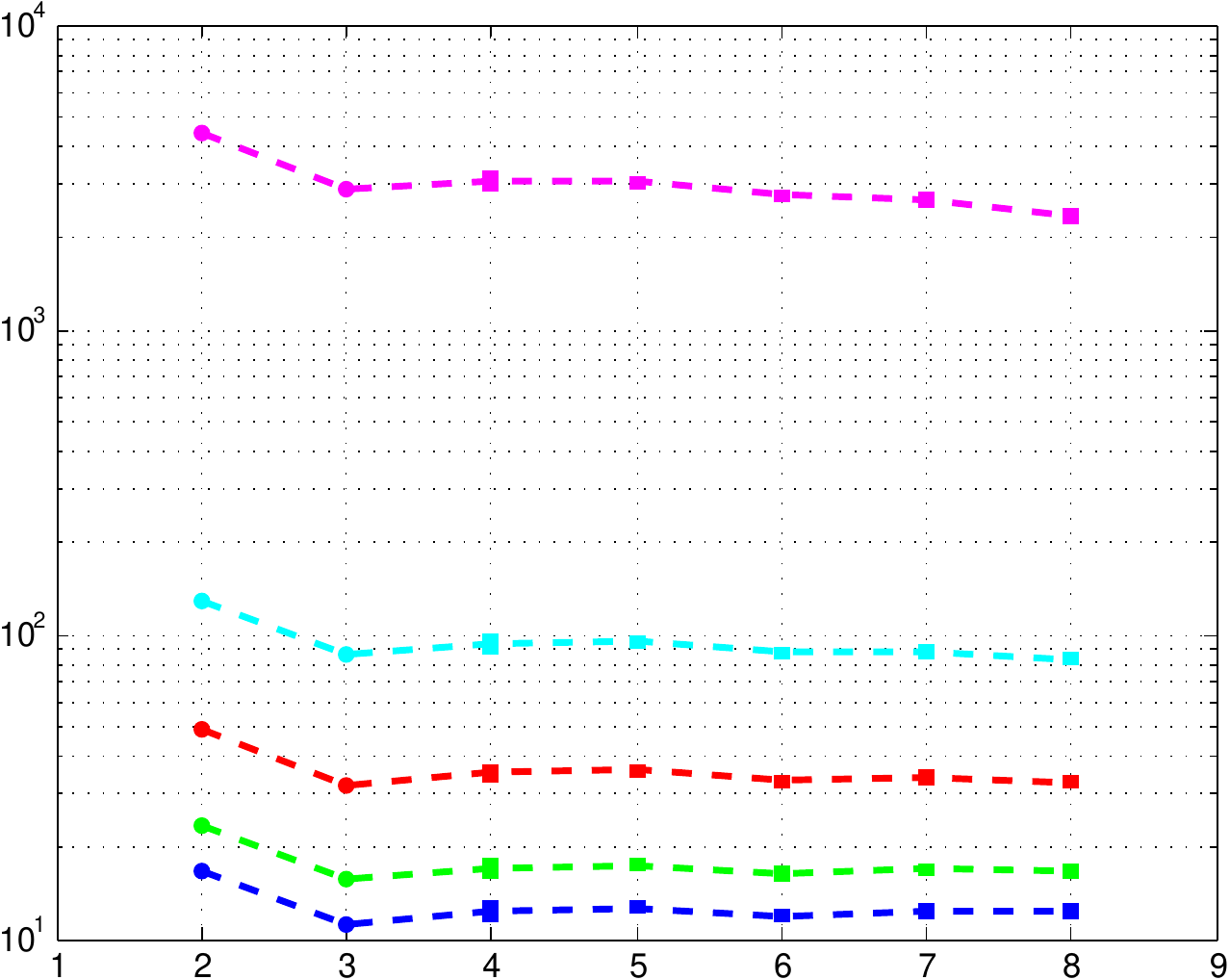}&
\includegraphics[scale=0.64]{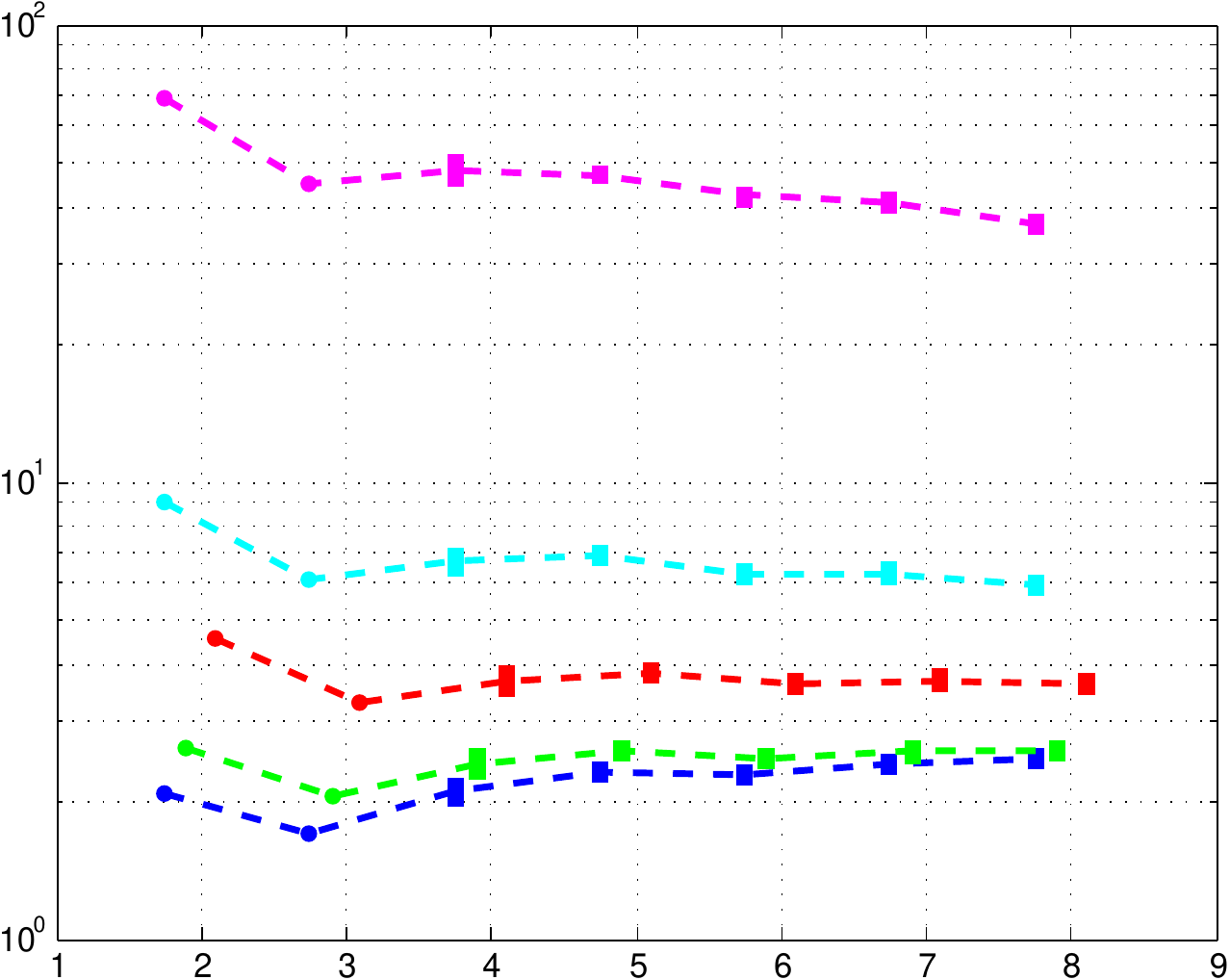}\\
\multicolumn{2}{c}{\hbox{\small Signal input: pulse, $2\leq k\leq 8$, $2\leq i\leq k$}}
\end{array}
$$
\caption{\label{figure1} {\small Detecting changes in the trend of a simple time series. Blue/green/red/cyan/magenta: $\nu=\infty/6/3/2/1$.  Left plots: ranges (vertical segments) of
$\rho_{k,2i-1}=\rho_{k,2i}$, vs. $k$. Right plots: ranges (vertical segments) of performance indexes $\rho_{k,2i-1}/\rho^*_{k,2i-1}=\rho_{k,2i}/\rho^*_{k,2i}$ vs. $k$. Ranges of $i$ and $k$ cover the domain where $\rho^*_{k,2i}=\rho^*_{k,2i-1}<R=10^4$. Charts are shifted horizontally to improve the plot readability.
On these plots both $\rho_{k,2i-1}$ and $\rho_{k,2i-1}/\rho^*_{k,2i-1}$ decrease with $\nu$.
}}
\end{figure}

We operate on time horizon ${d}=8$ and deal with $\sigma=1$ and with  5 values of the number $\nu$ of degrees of freedom of the Student distribution of $\eta$, specifically, the values $1,2,3,6,\infty$. In the experiments of this section we use parameter values $
R=10^4,\,\epsilon=0.01.$\footnote{We use  {\tt Mosek} and 
YALMIP Matlab toolbox \cite{mosek, Lofberg2004}  to solve corresponding optimization problems.}
The range of parameters $\rho_{k j}$ and of ratios $\rho_{k j}/\rho^*_{k,i}$
are presented in  Figure \ref{figure1}. Some comments are in order.
\begin{enumerate}
\item Relation $U_{2i-1}(\rho)=-U_{2i}(\rho)$ implies that $\rho_{k,2i-1}=\rho_{k,2i}$, $\rho^*_{k,2i-1}=\rho^*_{k,2i}$, $1\leq i\leq {d}$.
\item We display the range of quantities $\rho_{kj}$ and $\rho_{kj}/\rho^*_{kj}$ only for those values of $k,j$ for which $\rho^*_{kj}<R$, that is, ignore pairs $k,j$ for which already an optimistic lower bounds $\rho_{kj}^*$ on the magnitude of signal inputs of shape $j$ which can be detected, with the required risk, at time $k$ should be $\geq R$, which is forbidden by (\ref{meq2}). On a closest inspection, the ignored pairs $k,j$ are the pairs of the form $k,j=2i-1$ and $k,j=2i$ where
    \\
    $\quad$(a) $i>k$, or\\
    $\quad$(b) $k=1$, or \\
    $\quad$(c) [only for pulses!] $i=1<k$.\\
The reasons are clear: (a) stems from the fact that at time $k$ it is impossible to detect a whatever large signal input which starts at time $i>k$. (b) reflects the fact that the contribution of a whatever large signal input, if any, to the very first observation is fully masked by the initial condition $\alpha_0$, and in our model we do not impose any restrictions on this initial condition. (c) is of a similar origin: when the signal inputs are pulses, of a whatever magnitude,  at time 1, are fully masked by the initial conditions and thus cannot be detected at all.
\item As it could be guessed, when the signal inputs are steps, the quantities $\rho_{kj}$, for $j$ fixed, decrease as $k\geq j$ grows, since influence of step-change on our observations accumulates with time. In contrast, no such phenomenon is observed for pulse signal inputs where there is ``nothing to accumulate.''\footnote{Or, rather, that the noise in the states $\alpha_k$ of the model accumulates at the same rate, thus {cancelling} the effect of the growing observation sample.}
\item The conservatism of our decision rules, as presented in the tables, while unpleasant, seems to be not {too high} when $\nu\geq 3$, and becomes really {arresting} when $\nu=1$. The origin of this phenomenon is quite transparent. The conservatism  seems to stem primarily from systematic use in our constructions, for absence of something better, of the union bounds for probabilities. For example,  when computing $\rho^*_{kj}$, we allow for the probability of false alarm at time $k$ to be as large as $\epsilon$, while in our decision rules, we ``distribute'' this probability  between
${d}$ instants where we make our decisions. Similarly, when computing $\rho^*_{kj}$, we act {\sl as if} the only alternative to the nuisance hypothesis were a particular signal hypothesis $u\in U_j(\rho)$, while in fact we have several signal hypotheses to consider and should take into account the resulting ``accumulation of risk.'' As a result, we require from pairwise tests participating in $\T_k$   to have risk essentially smaller than $\epsilon$, which in the case of a ``heavy tail'' noise distribution allowed by our model requires an essentially larger magnitudes of detectable signal inputs than those allowing for detection when the shape of signal input is known in advance. And indeed, we see that the ratios $\rho_{kj}/\rho^*_{kj}$ rapidly increase as $\nu$ decreases.
    \end{enumerate}
\hide{
\begin{table}
\begin{center}
{\tiny\rm
\begin{tabular}{|c|}
\hline
\begin{tabular}{||c||c|c|c|c|c|c|c|c|c|c||} \cline{2-9}
\multicolumn{1}{c|}{}&$i=1$&$i=2$&$i=3$&$i=4$&$i=5$&$i=6$&$i=7$&$i=8$\\ \hline\hline
$k=1$&&&&&&&&\\ \hline
$k=2$&11.7/1.5&11.7/1.5&&&&&&\\ \hline
$k=3$&5.9/1.3&5.9/1.3&10.2/1.5&&&&&\\ \hline
$k=4$&4.8/1.5&4.8/1.5&6.3/1.6&11.0/2.0&&&&\\ \hline
$k=5$&4.2/1.7&4.2/1.7&5.0/1.7&6.4/1.8&11.0/2.2&&&\\ \hline
$k=6$&3.8/1.8&3.8/1.8&4.3/1.8&4.9/1.8&6.2/1.8&10.7/2.2&&\\ \hline
$k=7$&3.5/1.8&3.5/1.8&3.8/1.8&4.3/1.8&4.9/1.8&6.4/1.9&10.9/2.3&\\ \hline
$k=8$&3.2/1.8&3.2/1.8&3.5/1.8&3.8/1.8&4.3/1.8&5.0/1.8&6.4/1.9&10.9/2.3\\ \hline
\hline
 \end{tabular}
 \\
 {\scriptsize  $\nu=\infty$}\\
 \hline
\begin{tabular}{||c||c|c|c|c|c|c|c|c||} \cline{2-9}
\multicolumn{1}{c|}{}&$i=1$&$i=2$&$i=3$&$i=4$&$i=5$&$i=6$&$i=7$&$i=8$\\ \hline\hline
$k=1$&&&&&&&&\\ \hline
$k=2$&16.3/1.8&16.3/1.8&&&&&&\\ \hline
$k=3$&8.2/1.5&8.2/1.5&14.0/1.8&&&&&\\ \hline
$k=4$&6.7/1.7&6.7/1.7&8.8/1.8&15.3/2.2&&&&\\ \hline
$k=5$&5.8/1.8&5.8/1.8&6.9/1.8&8.9/1.9&15.1/2.3&&&\\ \hline
$k=6$&5.2/1.8&5.2/1.8&5.9/1.8&6.8/1.8&8.5/1.9&14.8/2.3&&\\ \hline
$k=7$&4.8/1.8&4.8/1.8&5.3/1.9&5.9/1.9&6.8/1.9&8.8/1.9&15.0/2.4&\\ \hline
$k=8$&4.4/1.8&4.4/1.8&4.8/1.9&5.3/1.9&5.9/1.9&6.8/1.9&8.8/1.9&15.0/2.4\\ \hline
\hline
 \end{tabular}
 \\
 {\scriptsize  $\nu=6$}\\
 \hline
\begin{tabular}{||c||c|c|c|c|c|c|c|c||} \cline{2-9}
\multicolumn{1}{c|}{}&$i=1$&$i=2$&$i=3$&$i=4$&$i=5$&$i=6$&$i=7$&$i=8$\\ \hline\hline
$k=1$&&&&&&&&\\ \hline
$k=2$&33.0/3.0&33.0/3.0&&&&&&\\ \hline
$k=3$&16.5/2.3&16.5/2.3&28.5/2.9&&&&&\\ \hline
$k=4$&13.6/2.5&13.6/2.5&17.9/2.7&30.8/3.3&&&&\\ \hline
$k=5$&11.9/2.6&11.9/2.6&14.0/2.6&18.2/2.8&30.8/3.3&&&\\ \hline
$k=6$&10.7/2.6&10.7/2.6&12.0/2.6&13.7/2.6&17.4/2.7&30.1/3.3&&\\ \hline
$k=7$&9.7/2.6&9.7/2.6&10.7/2.6&11.9/2.6&13.8/2.6&17.9/2.8&30.5/3.3&\\ \hline
$k=8$&8.9/2.6&8.9/2.6&9.8/2.6&10.7/2.6&11.9/2.6&13.9/2.6&17.9/2.8&30.5/3.3\\ \hline
\hline
 \end{tabular}
 \\
 {\scriptsize  $\nu=3$}\\
 \hline
\begin{tabular}{||c||c|c|c|c|c|c|c|c||} \cline{2-9}
\multicolumn{1}{c|}{}&$i=1$&$i=2$&$i=3$&$i=4$&$i=5$&$i=6$&$i=7$&$i=8$\\ \hline\hline
$k=1$&&&&&&&&\\ \hline
$k=2$&85.4/5.8&85.4/5.8&&&&&&\\ \hline
$k=3$&42.7/4.2&42.7/4.2&73.9/5.1&&&&&\\ \hline
$k=4$&35.1/4.3&35.1/4.3&46.1/4.6&79.3/5.6&&&&\\ \hline
$k=5$&30.8/4.4&30.8/4.4&36.0/4.5&46.7/4.7&79.3/5.7&&&\\ \hline
$k=6$&27.5/4.4&27.5/4.4&30.8/4.4&35.4/4.4&44.9/4.5&77.5/5.6&&\\ \hline
$k=7$&25.0/4.4&25.0/4.4&27.6/4.4&30.8/4.4&35.7/4.4&46.1/4.6&78.7/5.6&\\ \hline
$k=8$&23.0/4.4&23.0/4.4&25.2/4.4&27.6/4.4&30.8/4.4&36.0/4.5&46.1/4.6&78.7/5.6\\ \hline
\hline
 \end{tabular}
 \\
 {\scriptsize  $\nu=2$}\\
 \hline
\begin{tabular}{||c||c|c|c|c|c|c|c|c||} \cline{2-9}
\multicolumn{1}{c|}{}&$i=1$&$i=2$&$i=3$&$i=4$&$i=5$&$i=6$&$i=7$&$i=8$\\ \hline\hline
$k=1$&&&&&&&&\\ \hline
$k=2$&2246/35&2246/35&&&&&&\\ \hline
$k=3$&1123/25&1123/25&1953/30&&&&&\\ \hline
$k=4$& 928/25& 928/25&1211/27&2109/33&&&&\\ \hline
$k=5$& 806/25& 806/25& 952/26&1230/27&2090/33&&&\\ \hline
$k=6$& 723/25& 723/25& 811/25& 933/25&1182/26&2051/32&&\\ \hline
$k=7$& 659/25& 659/25& 728/25& 811/25& 938/25&1211/27&2090/33&\\ \hline
$k=8$&  610/25& 610/25& 664/26& 728/25& 811/25& 942/26&1221/27&\ 2090/33\ \\ \hline
\hline
 \end{tabular}
\\
 {\scriptsize  $\nu=1$}\\
 \hline
  \end{tabular}
}
\end{center}
\caption{\label{tab1} Signal inputs: steps. Empty cells $ki$: $\rho^*_{k,2i-1}=\rho^*_{k,2i}=R$. Nonempty cells $ki$: the first entry is $\rho_{k,2i-1}=\rho_{k,2i}$, the second is the performance index
$\rho_{k,2i-1}/\rho^*_{k,2i-1}$.}
\end{table}
\begin{table}
\begin{center}
{\tiny\rm
\begin{tabular}{|c|}
\hline
\begin{tabular}{||c||c|c|c|c|c|c|c|c||} \cline{2-9}
\multicolumn{1}{c|}{}&$i=1$&$i=2$&$i=3$&$i=4$&$i=5$&$i=6$&$i=7$&$i=8$\\ \hline\hline
$k=1$&&&&&&&&\\ \hline
$k=2$&&11.7/1.5&&&&&&\\ \hline
$k=3$&&10.2/1.5&10.2/1.5&&&&&\\ \hline
$k=4$&&11.0/2.0&10.4/1.9&11.0/2.0&&&&\\ \hline
$k=5$&&11.0/2.2&10.6/2.1&10.6/2.1&11.0/2.2&&&\\ \hline
$k=6$&&10.7/2.2&10.3/2.1&10.8/2.3&10.3/2.1&10.7/2.2&&\\ \hline
$k=7$&&10.9/2.3&10.5/2.2&10.8/2.3&10.8/2.3&10.5/2.2&10.9/2.3&\\ \hline
$k=8$&&10.9/2.3&10.5/2.2&10.8/2.3&10.7/2.3&10.8/2.3&10.5/2.2&10.9/2.3\\ \hline
\hline
 \end{tabular}
 \\
 {\scriptsize  $\nu=\infty$}\\
 \hline
\begin{tabular}{||c||c|c|c|c|c|c|c|c||} \cline{2-9}
\multicolumn{1}{c|}{}&$i=1$&$i=2$&$i=3$&$i=4$&$i=5$&$i=6$&$i=7$&$i=8$\\ \hline\hline
$k=1$&&&&&&&&\\ \hline
$k=2$&&16.3/1.8&&&&&&\\ \hline
$k=3$&&14.0/1.8&14.0/1.8&&&&&\\ \hline
$k=4$&&15.3/2.2&14.3/2.1&15.3/2.2&&&&\\ \hline
$k=5$&&15.1/2.3&14.6/2.3&14.6/2.3&15.1/2.3&&&\\ \hline
$k=6$&&14.8/2.3&14.3/2.2&15.0/2.3&14.3/2.2&14.8/2.3&&\\ \hline
$k=7$&&15.0/2.4&14.4/2.3&14.9/2.3&14.9/2.3&14.4/2.3&15.0/2.4&\\ \hline
$k=8$&&15.0/2.4&14.5/2.3&14.9/2.3&14.7/2.3&14.9/2.3&14.5/2.3&15.0/2.4\\ \hline
\hline
 \end{tabular}
 \\
 {\scriptsize  $\nu=6$}\\
 \hline
\begin{tabular}{||c||c|c|c|c|c|c|c|c||} \cline{2-9}
\multicolumn{1}{c|}{}&$i=1$&$i=2$&$i=3$&$i=4$&$i=5$&$i=6$&$i=7$&$i=8$\\ \hline\hline
$k=1$&&&&&&&&\\ \hline
$k=2$&&33.0/3.0&&&&&&\\ \hline
$k=3$&&28.5/2.9&28.5/2.9&&&&&\\ \hline
$k=4$&&30.8/3.3&29.0/3.1&30.8/3.3&&&&\\ \hline
$k=5$&&30.8/3.3&29.8/3.2&29.8/3.2&30.8/3.3&&&\\ \hline
$k=6$&&30.1/3.3&29.0/3.2&30.5/3.3&29.0/3.2&30.1/3.3&&\\ \hline
$k=7$&&30.5/3.3&29.3/3.2&30.2/3.3&30.2/3.3&29.3/3.2&30.5/3.3&\\ \hline
$k=8$&&30.5/3.3&29.4/3.2&30.2/3.3&29.9/3.3&30.2/3.3&29.4/3.2&30.5/3.3\\ \hline
\hline
 \end{tabular}
 \\
 {\scriptsize  $\nu=3$}\\
 \hline
\begin{tabular}{||c||c|c|c|c|c|c|c|c||} \cline{2-9}
\multicolumn{1}{c|}{}&$i=1$&$i=2$&$i=3$&$i=4$&$i=5$&$i=6$&$i=7$&$i=8$\\ \hline\hline
$k=1$&&&&&&&&\\ \hline
$k=2$&&85.4/5.8&&&&&&\\ \hline
$k=3$&&73.9/5.1&73.9/5.1&&&&&\\ \hline
$k=4$&&79.3/5.6&75.1/5.3&79.3/5.6&&&&\\ \hline
$k=5$&&79.3/5.7&76.9/5.5&76.9/5.5&79.3/5.7&&&\\ \hline
$k=6$&&77.5/5.6&75.1/5.4&78.7/5.6&75.1/5.4&77.5/5.6&&\\ \hline
$k=7$&&78.7/5.6&75.7/5.4&78.1/5.6&78.1/5.6&75.7/5.4&78.7/5.6&\\ \hline
$k=8$&&78.7/5.6&76.3/5.5&78.1/5.6&77.5/5.6&78.1/5.6&76.3/5.5&78.7/5.6\\ \hline
\hline
 \end{tabular}
 \\
 {\scriptsize  $\nu=2$}\\
 \hline
\begin{tabular}{||c||c|c|c|c|c|c|c|c||} \cline{2-9}
\multicolumn{1}{c|}{}&$i=1$&$i=2$&$i=3$&$i=4$&$i=5$&$i=6$&$i=7$&$i=8$\\ \hline\hline
$k=1$&&&&&&&&\\ \hline
$k=2$&&2246/35&&&&&&\\ \hline
$k=3$&&1953/30&1953/30&&&&&\\ \hline
$k=4$&&2109/33&1973/31&2109/33&&&&\\ \hline
$k=5$&&2090/33&2031/32&2031/32&2090/33&&&\\ \hline
$k=6$&&2051/32&1973/31&2070/32&1973/31&2051/32&&\\ \hline
$k=7$&&2090/33&1992/31&2051/32&2051/32&1992/31&2090/33&\\ \hline
$k=8$&&2090/33&2012/31&2051/32&2051/32&2051/32&2012/31&\ 2090/33\ \\ \hline
\hline
 \end{tabular}
\\
 {\scriptsize  $\nu=1$}\\
 \hline

\end{tabular}
}
\end{center}
\caption{\label{tab2} Signal inputs: pulses. Empty cells $ki$: $\rho^*_{k,2i-1}=\rho^*_{k,2i}=R$. Nonempty cells $ki$: the first entry is $\rho_{k,2i-1}=\rho_{k,2i}$, the second is the performance index
$\rho_{k,2i-1}/\rho^*_{k,2i-1}$.}
\end{table}
}

\begin{figure}
\centering
$$\begin{array}{cc}
\includegraphics[scale=0.64]{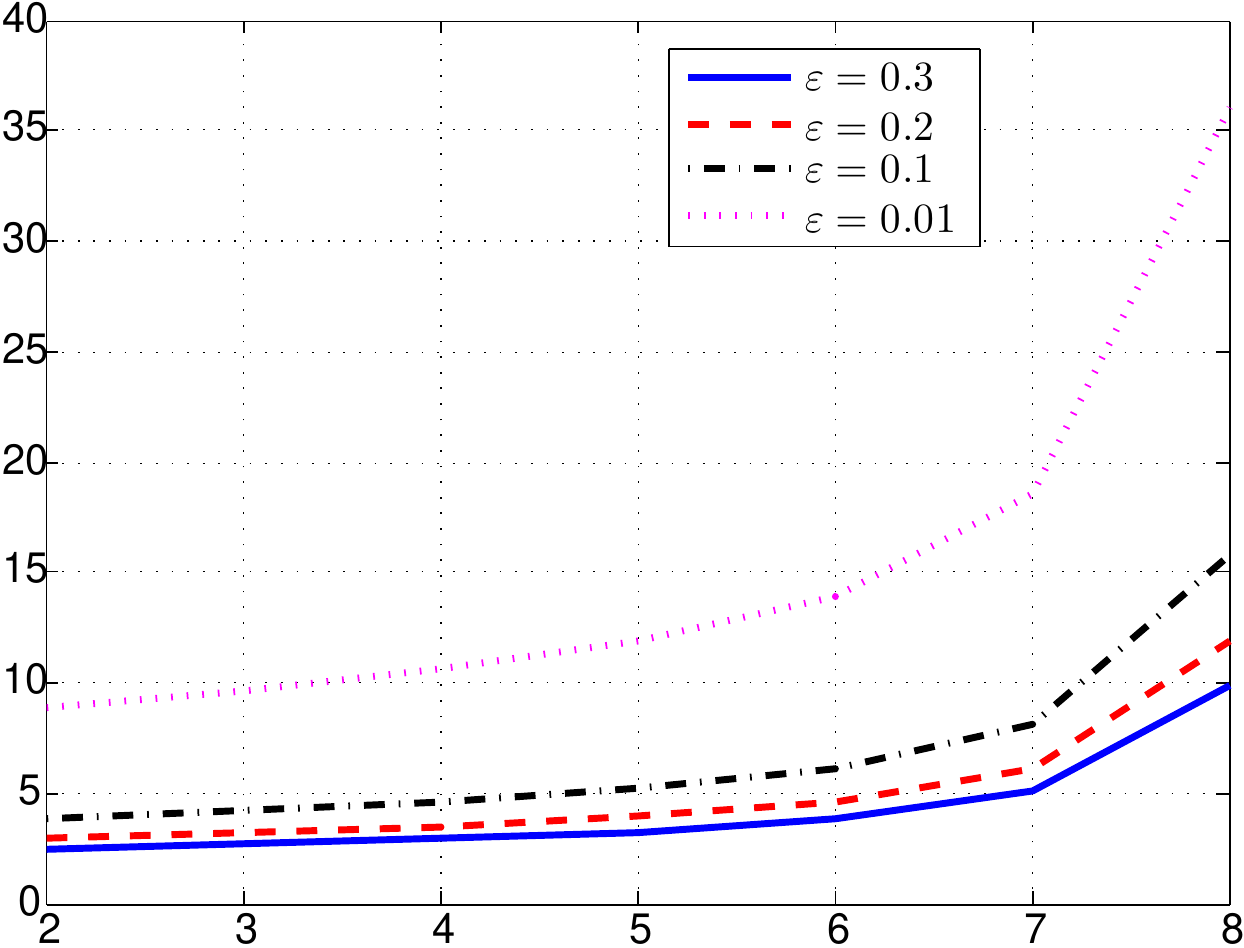}
&
\includegraphics[scale=0.60]{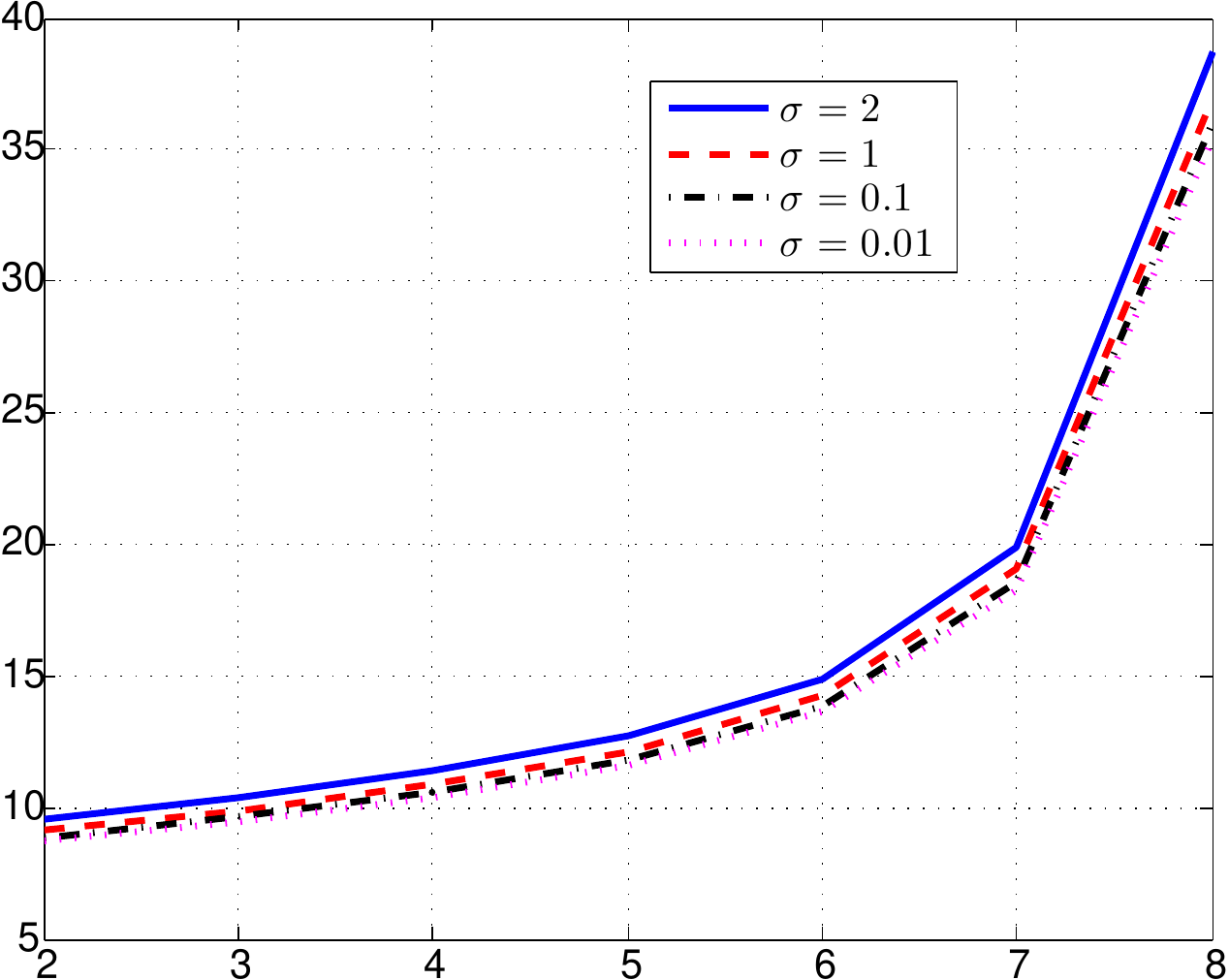}
\\
\multicolumn{2}{c}{\hbox{\small Signal input: step, $k=8$, $2\leq i\leq k$}}
\\
\includegraphics[scale=0.64]{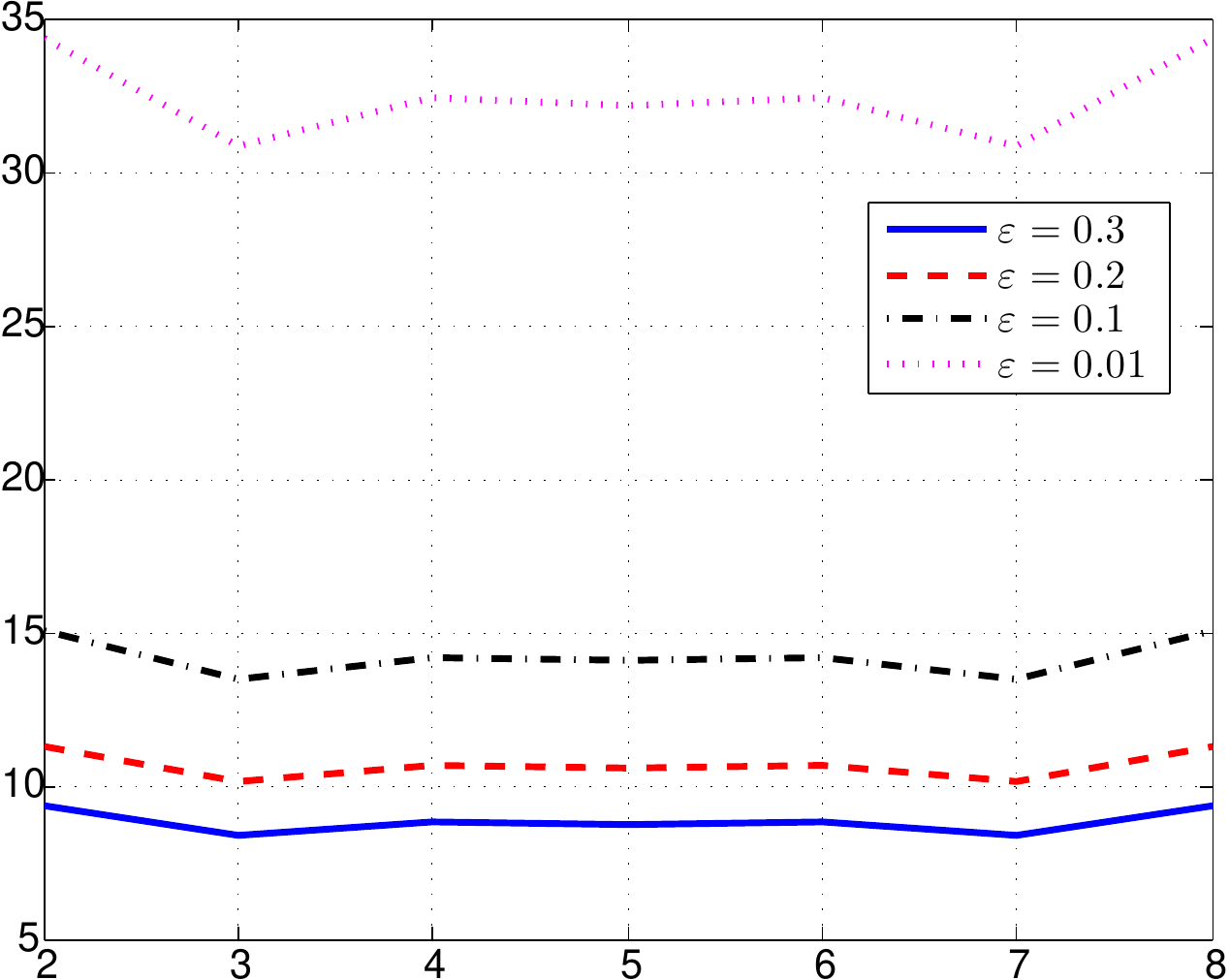}
&
\includegraphics[scale=0.64]{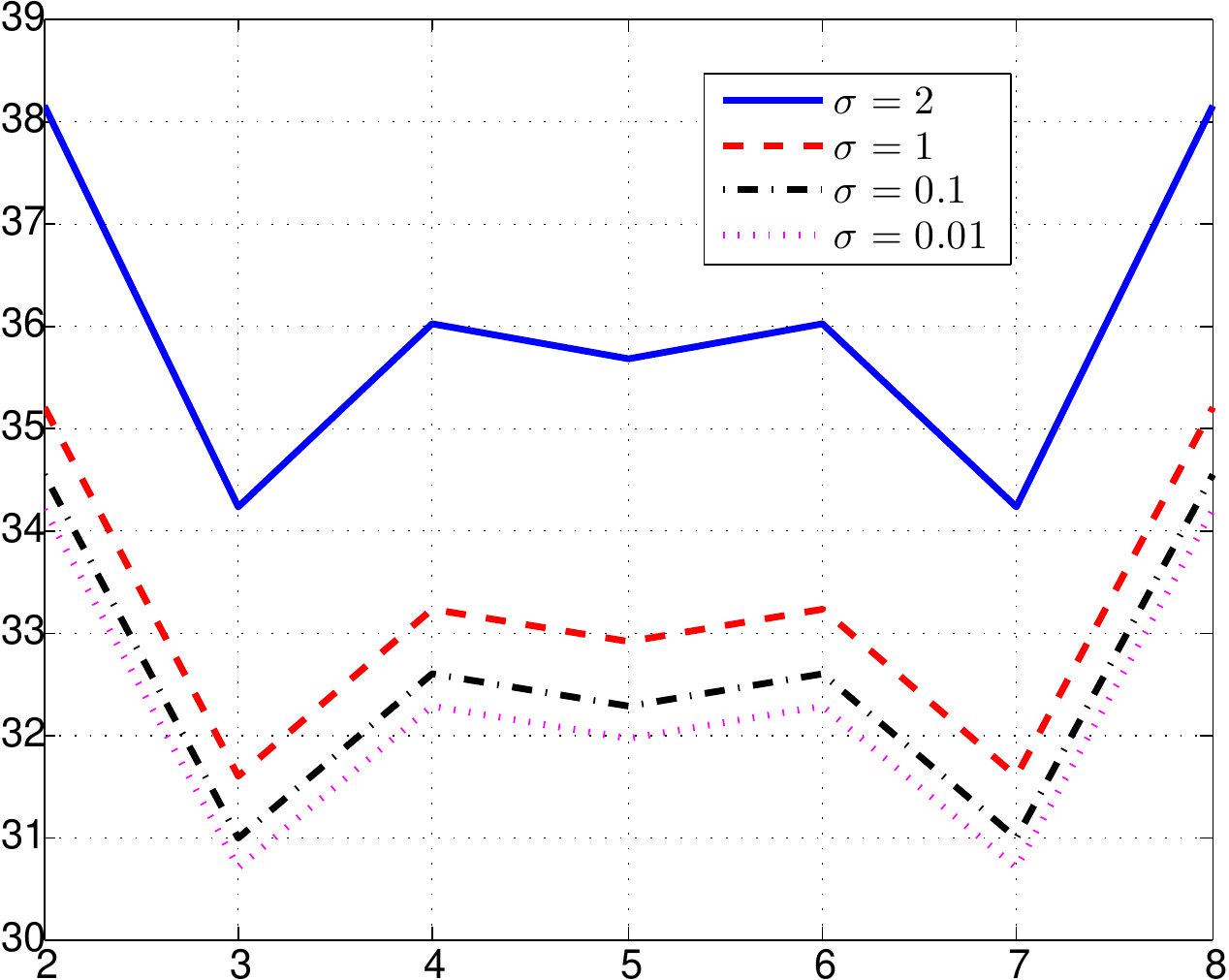}\\
\multicolumn{2}{c}{\hbox{\small Signal input: pulse, $k= 8$, $2\leq i\leq k$}}
\end{array}
$$
\caption{\label{RhosEpsSigDF} {\small Coefficients $\rho_{k,2i-1}=\rho_{k,2i}$ for $\nu=3$, $k=8$, and $i=2,\ldots,8$, as a function of the desired risk $\epsilon$ of the test (left pane)  and the standard deviation $\sigma$ of the Gaussian noise $\zeta_t$ (right pane).
}}
\end{figure}

We report on Figure \ref{RhosEpsSigDF} the evolution of parameters $\rho_{k j}$
for $k=8$ for step and pulse signals as a function of the risk of the test $\epsilon$ and
the standard deviation $\sigma$ of the Gaussian component $\zeta_t$ of the noise.
As expected, $\rho_{k j}$ increase with $\sigma$ and when $\epsilon$ decreases.\\

\subsubsection{{Change detection in linear dynamic system} revisited}\label{sectrefinement}
The methodology {developed in Section \ref{sectschemeII}} allows for {a straightforward refinement} which hopefully improves the resulting inference {performance}. {Note that the inference rules, as given in Sections \ref{sectschemeI} and \ref{sectschemeII} use very conservative bound for the probability of false alarm -- for multiple tests this probability is simply the sum of probabilities of false rejections of the nuisance hypothesis for each test.
This results in the increase of the testing thresholds $\rho_{kj}$, which  amounts to a ``logarithmic factor'' in the case of Gaussian observation noise, but becomes much more severe in the case of a heavy-tail noise distribution. One way to make the decision less cautious is to reduce the number of hypotheses to test by aggregating the alternatives.
 Here we illustrate the idea of the proposed modification on the simple numerical example in section \ref{numericalill}. } Specifically, {when building the detection procedure,} at time $k$ {we act} as follows:
\begin{itemize}
\item {We} compute the quantities $\rho^*_{kj}$, $j=1,...,2{d}$, and denote by ${J^*_k}$ the set of those $j$ for which $\rho^*_{kj}<R$. As it was explained, there is no reason to bother to detect at time $k$ signal inputs of shape $j\not\in {J^*_k}$.
\item Assume we have somehow associated {thresholds $\rho_{kj}\geq \rho^*_{kj}$ to} indexes $j\in{J^*_k}$; our goal, same as before, is to build a decision rule $\T_k$ which, given $\omega^k$, \\
    ---  with probability at least $1-\epsilon$, makes signal conclusion at time $k$, provided the input belongs to  $U_j(\rho_{kj})$ with some $j\in{J^*_k}$;\\
    ---  has false alarm probability at time $k$ (the probability to make signal conclusion when the input is a nuisance) $\leq\epsilon_k$, $\sum_k\epsilon_k=\epsilon$.
\end{itemize}
{Next, let us} color the sets $U_j(\rho_{kj})$, $j\in{J^*_k}$ (and their indexes) in a number $L$ of colors; let $I_\ell$ be the set of indexes $j\in {J^*_k}$ colored by color $\ell$. {We associate with each $\ell$ a
convex alternative $U^\ell$ -- the convex hull of ``alternatives of color $\ell$'':}
$$
U^\ell=\Conv\left(\bigcup\limits_{j\in I_\ell} U_j(\rho_{kj})\right),\,\ell=1,...,L,
$$
{and replace the original detection problem with the following one:}  given observation
$\omega^k$
we want to decide on the null hypothesis $H_0$ ``the input is nuisance'' (in our case, zero) vs. the alternative $H_1$ ``the input belongs to $\bigcup\limits_{\ell=1}^L U^\ell$.'' Same as before, our goal is to ensure probability of false alarm $\leq \epsilon_k$ and probability of miss $\leq\epsilon$. \par
Note that if we are able to do so, we meet our initial design specifications
-- with input from $U_j(\rho_{kj})$ for some $j\in {J^*_k}$, the probability of signal conclusion at time $k$ will be at least $1-\epsilon$.
To build the decision rule, let us use pairwise tests given by
the construction from Section \ref{sectschemeII}: we need $L$ tests, $\ell$-th of them deciding on $H_0$ vs.
the alternative $H_1^\ell: u\in U^\ell$, with risks $\epsilon^\ell_k$ of false alarm and $\epsilon$ of miss.
If we can ensure $\sum_{\ell=1}^L\epsilon^\ell_k \leq\epsilon_k$, we are done -- the decision rule which makes nuisance conclusion when all our $L$ tests ``vote'' for $H_0$, and makes the signal conclusion otherwise, is what we
are looking for.
Note that {the original} construction {in Section \ref{numericalill}} is of exactly this structure, with $L=\Card({J_k^*})$ (i.e., every $U_j(\rho_{kj})$,
$j\in {J_k^*}$, has its own color). It could make sense, however, to use {aggregated alternatives, thus reducing the} number of colors. When doing so, \par
--- on one hand, it {is} more difficult to decide on $H_0$ vs. the alternatives, because the image $A_kU^\ell$ of the ``aggregated'' set of signal inputs is closer to the image $\{0\}$ of the nuisance set than the images $A_kU_j(\rho_{kj})$, $j\in I_\ell$, of individual sets of signal inputs participating in the aggregation. Therefore,  to ensure the same risk, we now
need a somewhat larger separation of the images of the nuisance and the signal inputs  in the observation space;\par
--- on the other hand, we should now ``distribute'' $\epsilon_k$ {among} $L<\Card({J^*_k})$ miss probabilities $\epsilon^\ell_k$. {This would allow to} operate with larger  miss probabilities, thus reducing the {necessary} separation of the images of the nuisance and the signal inputs in the observation space.
\par
It is {hard} to tell in advance which  of these two opposite {effects will prevail}; an answer, however, could be provided by computation, and it makes sense to give to the outlined {modification} a try. To make things as simple as possible, let us act as follows.
\begin{itemize}
\item After the colors are assigned and the sets $U^\ell$, $\ell\leq L$, are built, we specify $\epsilon^\ell_k$ and the quantities $\alpha_{1k}$, $\alpha_{2k}$, $\delta_k$ to meet the requirements
\[
\epsilon^\ell_k={\epsilon_k\over L},\;1\leq\ell\leq L;\;\;\int\limits_{\alpha_{1k}}^\infty\gamma(s)ds={\epsilon_k\over L};\;\;\int\limits_{\alpha_{2k}}^\infty \gamma(s)ds=\epsilon;\;\delta_k={1\over 2}[\alpha_{1k}+\alpha_{2k}].
\]
\end{itemize}
{Let us suppose that} the relation
\begin{equation}\label{meq200}
\min_{u\in U^\ell}\half\|A_ku\|_2\geq\delta_k,\,\ell=1,...,L,
\end{equation}
is satisfied. We set
$$
h_{k\ell}=A_ku_{k\ell}/\|A_ku_{k\ell}\|_2,\;\;c_{k\ell}=\half h_{k\ell}^TA_ku_{k\ell},\;\;\phi_{k\ell}(\omega^k)=h_{k\ell}^T\omega_k-c_{k\ell},
$$
{where $u_{k\ell}$ are optimal solutions to the optimization problems in (\ref{meq200}). It is immediately seen that}  making at time $k$ the nuisance conclusion if and only if $\phi_{kj}(\omega^k)\geq \half[\alpha_{2k}-\alpha_{1k}]$, we ensure {simultaneously} the probability of false alarm at time $k$ at most $\epsilon_k$, and the probability of {miss} when the input belongs to $\bigcup\limits_{j\in{J^*_k}}U_j(\rho_{kj})$ at most $\epsilon$, thus meeting our design specifications.
\paragraph{{Specifying} $\rho_{kj}$.} The question we did not address so far is how to {choose} $\rho_{kj}$, $j\in{J^*_k}$. What we expect of these quantities is to ensure the validity of (\ref{meq200}), and the simplest way to achieve this goal is as follows. Setting $\rho_{kj}=\theta\rho_{kj}^*$, the left hand side in
(\ref{meq200}) is a nondecreasing function of $\theta$, and we can find by bisection the smallest $\theta=\theta_k\geq1$ for which (\ref{meq200}) takes place. After $\theta_k$ is found, we set $\rho_{kj}=\theta_k\rho^*_{kj}$, $j\in{J^*_k}$. Clearly, with this approach, the performance indexes $\rho_{kj}/\rho_{kj}^*$, $j\in {J^*_k}$, are all  equal to $\theta_k$.
\paragraph{How it works.} We applied the just outlined construction to the data underlying the numerical experiment described in Section \ref{numericalill}. Our implementation was the simplest possible: given $k$, we looked at all $j$'s such that $\rho^*_{kj}<R$; the set ${J^*_k}$ of these $j$'s with our data is nonempty only when $k\geq2$ and is either $\{1,...,2k\}$ (step signals), of $\{3,4,...,2k\}$ (pulse signals). A set $U_j(\rho)$ with odd index $j=2i-1$ (even index $j=2i$), $j\in{J^*_k}$, is comprised of {signals} which are zero before time $i$ and ``jump up/jump down'' at time $i$ depending on whether $j$ is even or odd. We color these sets in $L=2$ colors, depending on whether the corresponding indexes $j$ are odd or even, that is, we use at step $k$
$$
U^\ell=\bigcup\limits_{{j\in{J^*_k},\;j\,\hbox{\tiny\rm mod}\,2=\ell}} U_j(\theta_k\rho^*_{kj}),\,\ell=1,2,
$$
with $\theta_k$  as explained above.\par
\begin{figure}
$$\begin{array}{cc}
\includegraphics[scale=0.64]{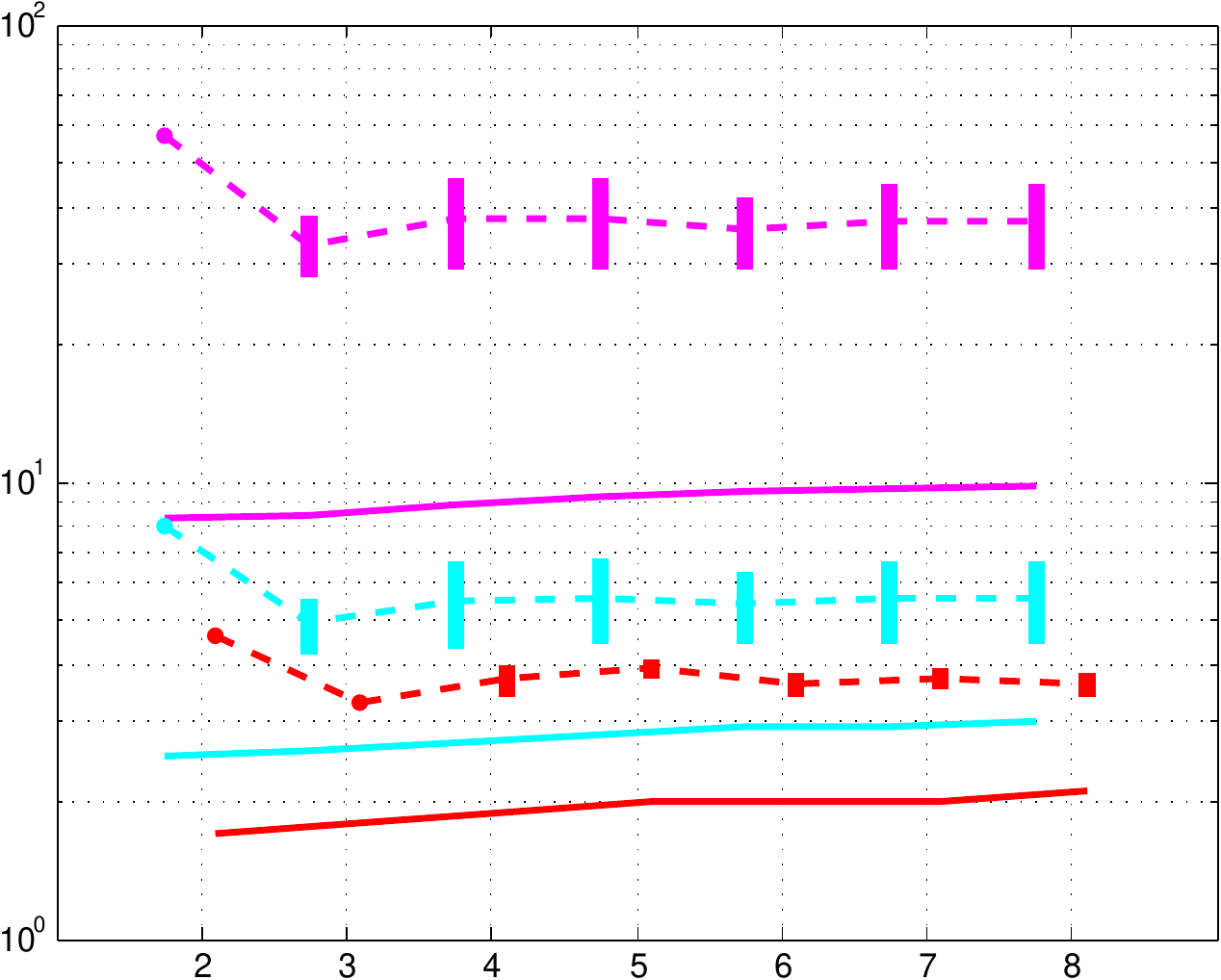} &
\includegraphics[scale=0.64]{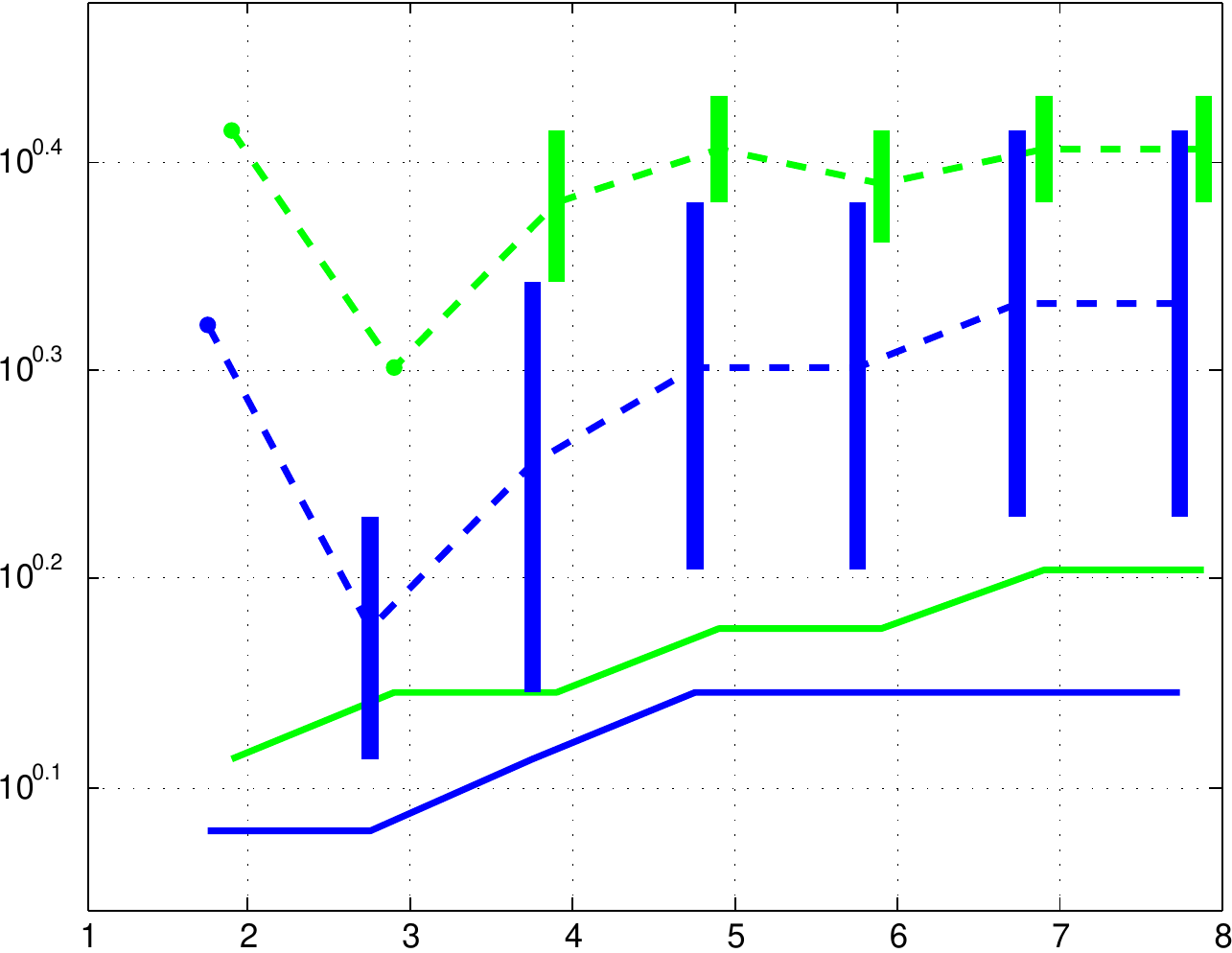} \\
 \multicolumn{2}{c}{\hbox{\small Signal input: step, $2\leq k\leq 8$}}\\
\includegraphics[scale=0.64]{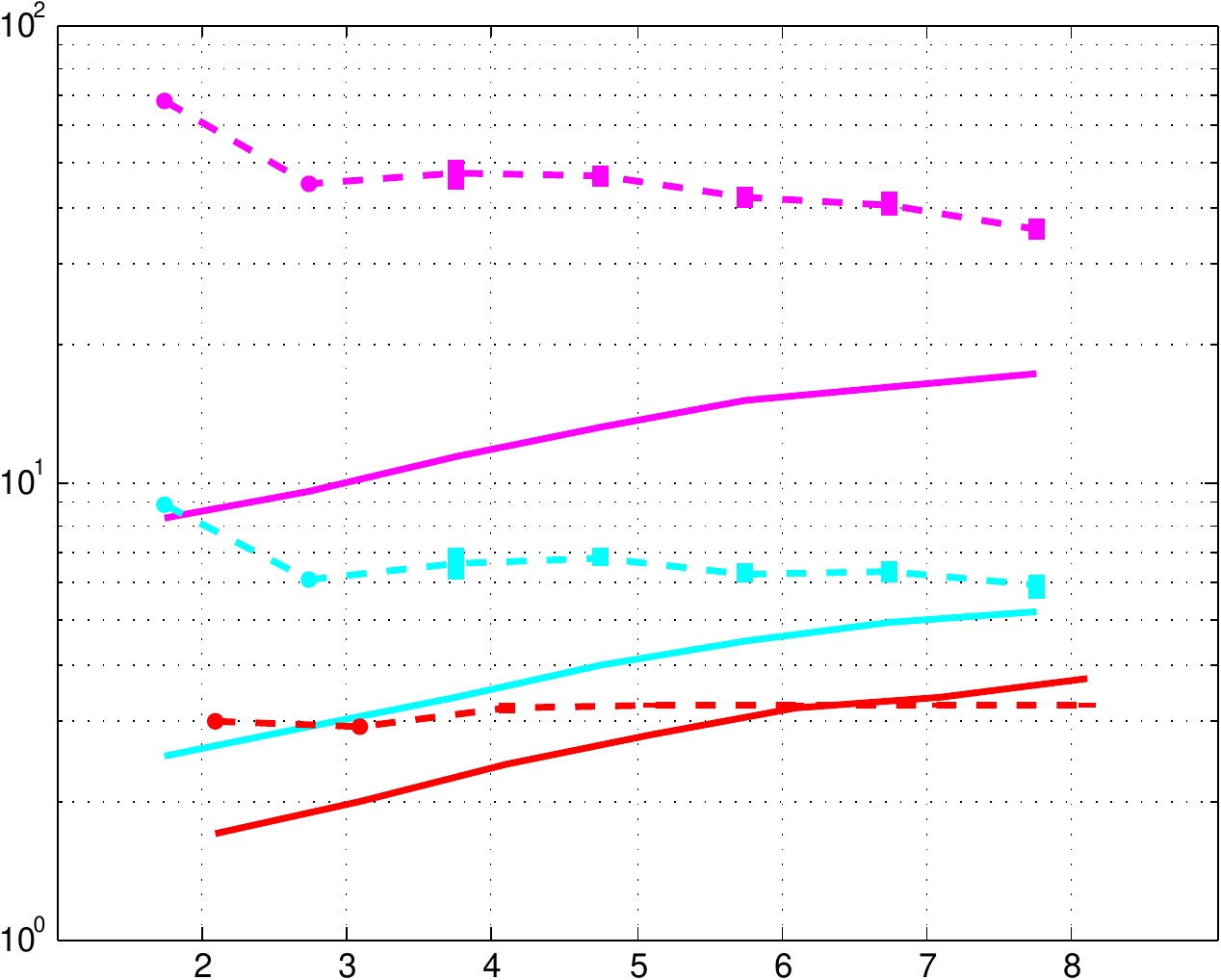} &
\includegraphics[scale=0.64]{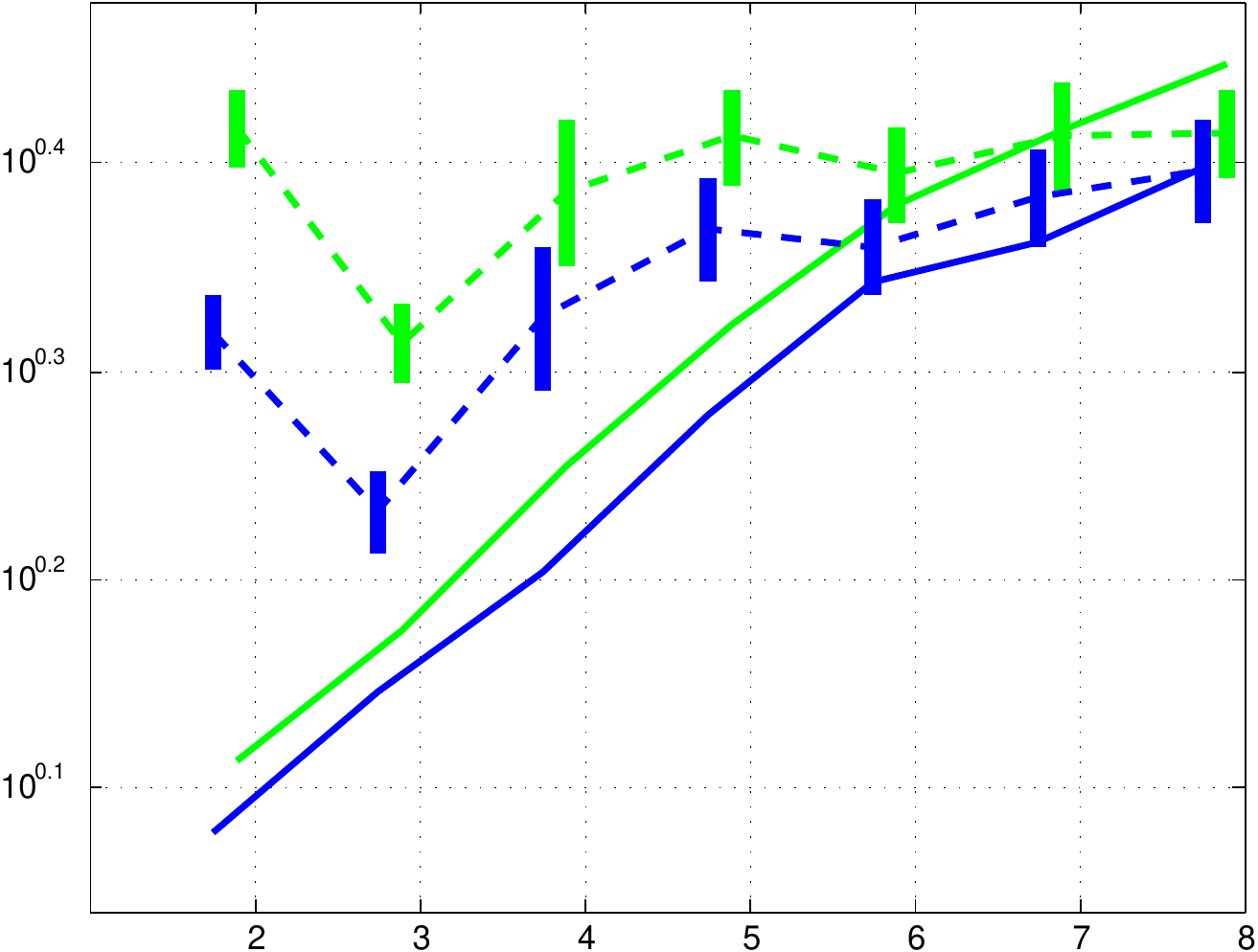} \\
\multicolumn{2}{c}{\hbox{\small Signal input: pulse, $2\leq k\leq 8$}}
\end{array}
$$
\caption{\label{figureconvagg} {\small Ranges of performance indexes $\rho_{kj}/\rho^*_{kj}$ (cf. right plots of Figure \ref{figure1}) as compared to the
performance index for the refined inference (solid lines). Red/cyan/magenta: $\nu=3/2/1$ (left plots), blue/green: $\nu=\infty/6$ (right plots). Charts are shifted horizontally to improve the plot readability.
On these plots $\rho_{k,j}/\rho^*_{k,j}$ decrease with $\nu$.
}}
\end{figure}

We present on Figure \ref{figureconvagg} the comparison of performance indexes $\rho_{kj}/\rho^*_{kj}$, $j\in {J^*_k}$, for {the original} inference routine (these indexes are presented on Figure \ref{figure1}) and the performance indexes of the just described {modified inference}. For our initial routine, the performance indexes slightly vary with $j\in{J^*_k}$, and we present their ranges; for the new routine, the performance indexes do not depend on $j$.
We observe that in the considered example the proposed straightforward aggregation of signal inputs typically {results in degradation of} the performance indexes, and improves these indexes  {significantly} for the ``heavy tailed'' {noise distributions} (the case of $\nu=1$).

\paragraph{Acknowledgement.} Research of the first author was partially supported by an
FGV grant, CNPq grants 307287/2013-0, and 401371/2014-0, FAPERJ grant E-26/201.599/2014.
Research of the second author was partially supported by the LabEx PERSYVAL-Lab (ANR-11-LABX-0025)
and CNPq grant 401371/2014-0.
Research of the third author was partially supported by NSF grants CCF-1523768, CCF-1415498, and CNPq grant 401371/2014-0.
\par
The authors are grateful to anonymous referees for their inspiring suggestions.
\bibliographystyle{abbrv}
\bibliography{ReferencesChP}

\begin{thebibliography}{10}

\bibitem{mosek}
E.~D. Andersen and K.~D. Andersen.
\newblock {\em The MOSEK optimization toolbox for MATLAB manual. Version 7.0},
  2013.
\newblock \url{http://docs.mosek.com/7.0/toolbox/}.

\bibitem{BN2001}
A.~Ben-Tal and A.~Nemirovski.
\newblock {\em Lectures on modern convex optimization: analysis, algorithms,
  and engineering applications}, volume~2.
\newblock Siam, 2001.

\bibitem{Burnashev1979}
M.~Burnashev.
\newblock On the minimax detection of an innacurately known signal in a white
  noise background.
\newblock {\em Theory Probab. Appl.}, 24:107--119, 1979.

\bibitem{Burnashev1982}
M.~Burnashev.
\newblock Discrimination of hypotheses for gaussian measures and a geometric
  characterization of the gaussian distribution.
\newblock {\em Math. Notes}, 32:757--761, 1982.

\bibitem{alltogether}
Y.~Cao, V.~Guigues, A.~Juditsky, A.~Nemirovski, and Y.~Xie.
\newblock Change detection via affine and quadratic detectors.
\newblock {\em arXiv preprint arXiv:1608.00524}, 2016.

\bibitem{chernoff1952}
H.~Chernoff.
\newblock A measure of asymptotic efficiency for tests of a hypothesis based on
  the sum of observations.
\newblock {\em The Annals of Mathematical Statistics}, pages 493--507, 1952.

\bibitem{eltoft2006multivariate}
T.~Eltoft, T.~Kim, and T.-W. Lee.
\newblock On the multivariate laplace distribution.
\newblock {\em IEEE Signal Processing Letters}, 13(5):300--303, 2006.

\bibitem{GJN2015}
A.~Goldenshluger, A.~Juditsky, and A.~Nemirovski.
\newblock Hypothesis testing by convex optimization.
\newblock {\em Electronic Journal of Statistics}, 9(2):1645--1712, 2015.

\bibitem{cvx2014}
M.~Grant and S.~Boyd.
\newblock {\em The {\tt CVX} Users’ Guide. {Release} 2.1}, 2014.
\newblock \url{http://web.cvxr.com/cvx/doc/CVX.pdf}.

\bibitem{Ingster2002}
Y.~Ingster and I.~A. Suslina.
\newblock {\em Nonparametric goodness-of-fit testing under Gaussian models},
  volume 169 of {\em Lecture Notes in Statistics}.
\newblock Springer, 2002.

\bibitem{Seq2015}
A.~Juditski and A.~Nemirovski.
\newblock On sequential hypotheses testing via convex optimization.
\newblock {\em Automation and Remote Control}, 76:809--825, 2015.

\bibitem{AffineDetectors}
A.~Juditsky and A.~Nemirovski.
\newblock Hypothesis testing via affine detectors.
\newblock {\em Electronic Journal of Statistics}, 10:2204--2242, 2016.

\bibitem{kotz2004multivariate}
S.~Kotz and S.~Nadarajah.
\newblock {\em Multivariate t-distributions and their applications}.
\newblock Cambridge University Press, 2004.

\bibitem{Lehmann2006}
E.~L. Lehmann and J.~P. Romano.
\newblock {\em Testing statistical hypotheses}.
\newblock Springer Science \& Business Media, 2006.

\bibitem{Lofberg2004}
J.~L{\"{o}}fberg.
\newblock Yalmip : A toolbox for modeling and optimization in matlab.
\newblock {\em In Proceedings of the IEEE CACSD Conference}, 2004.

\end{thebibliography}
\appendix
\section{Proofs}
\subsection{Proof of Proposition \ref{propnew}}\label{propnewproof}
In what follows, for functions $f, g: \bf{R} \rightarrow \bf{R}$, we say that $f$ dominates $g$ (notation: $f\succeq g$, or, equivalently, $g\preceq f$), if
$$
\int\limits_\delta^\infty f(s)ds\geq\int\limits_\delta^\infty g(s)ds, \;\forall \delta \geq 0.
$$Let
\begin{itemize}
\item $\E$ be the family of even probability densities on the real axis,
\item $\N$ be the family  of nice functions on the axis.
\end{itemize}
 {\bf 1$^o.$}
Note that $\succeq$ clearly is transitive: if $f\succeq g$ and $g\succeq h$, then $f\succeq h$. Furthermore,
 $\N\subset \E$ (by definition of $\N$), and $\E$ is closed with respect to taking convolution (evident). We also need the following technical facts.
\\{\bf 1$^o.$a} $\N$ is closed with respect to taking convolution.
\par Indeed, let $f,g\in \N$. The fact that $f\star g$ is even and continuous on the real axis is evident. Therefore in order to show that $f\star g\in \N$ it suffices to verify that $h=f\star g$ is nonincreasing on the nonnegative ray.
For $0< z \leq f(0)$, denoting
$
{\bar y}(z) = \min \left\{y\geq 0 : f(y)=z \right\},
$
we have for every $x \geq 0$,
$$
h(x)= \int\limits_{-\infty}^\infty f(x-y) g(y) dy = \int\limits_{-\infty}^\infty \big(  \int\limits_{0}^{f(x-y)} dz \big) g(y) dy
= \int\limits_{0}^{f(0)} H(x,z) dz
$$
where $H(x,z) =\int\limits_{x- {\bar y}(z)}^{x+ {\bar y}(z)} g(y) dy$ for $0 < z \leq f(0)$ and $H(x,0)=0$.
To conclude we observe that for every fixed $f(0) \geq z \geq 0$, $H(\cdot,z)$ is a differentiable, nonnegative, and nonincreasing function on the nonnegative ray.
Indeed, for $0 <z \leq f(0)$, the derivative $H_x'(x,z)$ of this function at $x$ is
$H_x'(x,z)=g(x+ {\bar y}(z)  ) - g(x- {\bar y}(z) )$
and it is clear that this quantity is $\leq 0$ since $g \in \N$.
We have checked that for every fixed $z \geq 0$, $H(\cdot,z)$ is nonincreasing on the nonnegative ray.
It follows that $h$ is nonincreasing on the nonnegative ray.
\\{\bf 1$^o.$b}
Let $\bar{f},f\in\E$, $g \in\N$, and let $\bar{f}\succeq f$. Then $\bar{f}\star g\succeq f\star g$.
\par  Let us verify that for all $x\geq 0$,
$\bar{H}(x)\leq H(x)$, where $\bar{H}$ and $H$ are cumulative distribution functions (c.d.f.) of the densities $\bar{f}\star g$ and ${f}\star g$, respectively. Observe that
\[
H(x)=
\int_{-\infty}^\infty g(s) F(x-s)ds=\int_{-\infty}^\infty g(x-t)F(t)dt,
\]
same as
\[
\bar{H}(x)=\int_{-\infty}^\infty g(x-t)\bar{F}(t)dt,
\]
where $F$ and $\bar{F}$ are the c.d.f.'s of the densities $f$ and $\bar{f}$, respectively. Thus, when setting  $\Delta(s)=F(s)-\bar{F}(s)$, we get
\[\begin{array}{rcl}
H(x)-\bar{H}(x)&=&\int_{-\infty}^\infty g(x-t)\Delta(t)dt
= \int_{0}^\infty g(x-t)\Delta(t)dt+\int_{0}^\infty g(x+t)\Delta(-t)dt\\
\mbox{[because $\Delta(-t)=-\Delta(t)$]}&=& \int_{0}^\infty \Delta(t)[g(x-t)-g(x+t)]dt\geq 0,\end{array}
\]
where the final ``$\geq$'' is due to $\Delta(t)\geq 0$ and $g(x-t)\geq g(x+t)$ for $x,t\geq 0$ (recall that $g\in \N$) for $t\geq 0$.
\\{\bf 1$^o.$c.}
Let $f\in\E$, $\rho\in(0,1)$, and let $f_\rho(s)=\rho^{-1}f(\rho^{-1}s)$. Then $f_\rho\in\E$ and $f_\rho\preceq f$.
\par
Indeed, the inclusion $f_\rho\in\E$ is obvious. On the other hand, if $\xi\sim f$, one has $P(\rho \xi\geq x)=P(\xi\geq x/\rho)\leq P(\xi\geq x)$ for $0<\rho<1$ and $x\geq 0$, what is exactly $f_\rho\preceq f$.
\\{\bf 1$^o.$d.} Let $f\in\N$ and $g\in \E$. Then $f\star g\succeq f$.
\par
Note that the c.d.f. $H$ of $f\star g$ satisfies
\[
H(x)=\int_{-\infty}^\infty g(s) F(x-s)ds=\int_{0}^\infty g(s) [F(x-s)+F(x+s)]ds,
\]
where $F$ is the c.d.f. of $f$ (recall that $g\in \E$), and therefore for $x\geq0$ it holds
\[
H(x)-F(x)=\int_{0}^\infty g(s) [F(x-s)+F(x+s)-2F(x)]ds\leq 0
\]
due to the concavity of $F$ on $\bR_+$.
\paragraph{2$^0$.} Now let $p(\cdot)\in\P^n$, $q(\cdot)\in\P^m$, and let $e\in\bR^r$, $\|e\|_2=1$, be given. Let us set $e_1=A^T\Theta^{-1}e$, $e_2=B^T\Theta^{-1}e$, let $\eta\sim p$ and $\zeta\sim q$ be independent, and let $\xi=\Theta^{-1}[A\eta+B\zeta]$. We have
$$
\omega:=e^T\xi=[A^T\Theta^{-1}e]^T\eta+[B^T\Theta^{-1}e]^T\zeta=e_1^T\eta+e_2^T\zeta.
$$
Observe that $e_1^Te_1=e^T\Theta^{-1}AA^T\Theta^{-1}e\leq 1$ due to $\Theta^{-1}AA^T\Theta^{-1}\preceq I_r$, and for similar reasons $e_2^Te_2\leq 1$. We denote
$\rho_1=\|e_1\|_2 \mbox{ and }\rho_2=\|e_2\|_2,$
so that $\rho_1,\rho_2\in[0,1]$. When $\rho_\chi>0$, $\chi=1,2$, we set
$
 \bar{e}_\chi=\rho_\chi^{-1}e_\chi.$
Now, let $f_1\in\E$ (respectively, $f_2\in\E$) be the density of the scalar random variable $e_1^T\eta$ (respectively, $e_2^T\zeta$). Note that when $\rho_1>0$ ($\rho_2>0$) random variable $e_1^T\eta$ ($e_2^T\zeta$) indeed has a density, and this density is even. Let also $\bar{f}_1\in\E$ ($\bar{f}_2\in\E$) be the density of $\bar{e}_1^T\eta$ ($\bar{e}_2^T\zeta$).
\paragraph{2$^0.$a.}
Assume for a moment that $\rho_1>0$ and $\rho_2>0$.   Then
\begin{itemize}
\item $p\in\P^n_\mu$, whence $f_1\in\E$, $\bar{f}_1\in\E$ and $f_1(s)=\rho_1^{-1}\bar{f}_1(\rho_1^{-1} s)$, whence $f_1\preceq \bar{f}_1$ by {\bf 1$^o.$c}.
Since $\bar{f}_1\preceq \mu$ and $\preceq$ is transitive, we have also $f_1\preceq \mu$.
\item $q\in \P^m$, and $\P^m$ is a completely monotone subfamily of $\P^m_\nu$, whence $f_2\in\N$, $\bar{f}_2\in\N$ and, same as above, $f_2\preceq\bar{f}_2\preceq \nu$, whence also $f_2\preceq \nu$.
\item The density of $\omega$ is $f_1\star f_2$. We have
$$
\begin{array}{lll}
&f_1\star f_2 \preceq \mu\star f_2&\hbox{\ [by {\bf 1$^o.$b} in view of $f_1,\mu\in \E$, $f_2\in\N$ and $f_1\preceq\mu$],}\\
&f_2\star\mu\preceq \mu\star \nu  &\hbox{\ [by {\bf 1$^o.$b} in view of $f_2,\nu\in\E$, $\mu\in\N$ and $f_2\preceq \nu$].}\\
\end{array}
$$
Hence, $f_1\star f_2\preceq \gamma:=\mu\star\nu$, such that $\gamma\in \E$, and (\ref{aeq5}) follows. Besides this, in the case in question $\omega$ has an even density.
\end{itemize}
\paragraph{2$^0$.b.}
Now let $\rho_1=0$. Then $\rho_2>0$ due to $AA^T+BB^T\succ0$, and the probability density of $\omega$ is ${f}_2$. We have $f_2\preceq\bar{f}_2\preceq\nu\preceq\mu\star\nu$ (the concluding $\preceq$ is due to {\bf 1$^o.$d}), and (\ref{aeq5}) follows.  Similarly, when $\rho_2=0$, we have $\rho_1>0$, and the probability density of $\omega$ is $f_1$. Similarly to the
case where $\rho_1=0$, we have $f_1\preceq\bar{f}_1\preceq \mu\preceq\mu\star\nu$, and (\ref{aeq5}) follows. Furthermore, as we have seen, $\omega$ always has an even probability density. Finally, $\gamma=\mu\star \nu$ is nice by {\bf 1$^o.$a}. (i) is proved.
\paragraph{3$^0$.}
To prove (ii), note that in the notation of item 2$^0$ and under the premise of (ii), $\rho_\chi>0$ implies that $f_\chi\in\N$. Hence, due to {\bf 1$^o.$a}, the distribution of $\omega$ has density from $\N$ (this density is $f_1\star f_2$ when $\rho_1,\rho_2>0$, $f_2$ when $\rho_1=0$, and  $f_1$ when $\rho_2=0$), which combines with (i) to imply (ii).\qed
\subsection{Proof of Proposition \ref{prop1}}\label{prop1proof}
{\bf Proof.} When $x\in X_1$, we have $h_*^Tx\geq h_*^Tx^1_*=c_*+\delta$, thus
\[\begin{array}{rcl}
\{{\xi:\,}s_*(x+\xi)<\half(\alpha_2-\alpha_1)\}&=&\{\xi:h_*^T(x+\xi)<c_*+\half(\alpha_2-\alpha_1)\}\\
&=&
\{{\xi:\,}h_*^T\xi<c_*-h_*^Tx+\half(\alpha_2-\alpha_1)\}\\
&\subseteq& \{{\xi:\,}h_*^T\xi<c_*-h_*^T x^1_*+\half(\alpha_2-\alpha_1)\}\\
&=&\{{\xi:\,}h_*^T\xi<-\delta+\half(\alpha_2-\alpha_1)\}
\subseteq
\{h_*^T\xi<-\alpha_1\},
\end{array}
\]
where the last inclusion is due to  $\alpha_1+\alpha_2\leq 2\delta$.
{Therefore,} for $p\in \P$ it holds
\[
\Prob_{\xi\sim p}\{s_*(x+\xi)<\half(\alpha_2-\alpha_1)\}\leq\int\limits_{h_*^T\xi<-\alpha_1} p(\xi)d\xi\underbrace{=}_{(a)}\int\limits_{h_*^T\xi>\alpha_1} p(\xi)d\xi\underbrace{\leq}_{(b)} P_\gamma(\alpha_1)
\]
where $(a)$ is due to the fact that $p(\cdot)$ is even, {and $(b)$} is due to $p\in\P_\gamma$ and $\alpha_1\geq0$.
{When} $x\in X_2$, we have $h_*^Tx\leq h_*^Tx^2_*=c_*-\delta$, and { using a completely similar argument we conclude that} for $p\in\P$ it holds

\[\Prob_{\xi\sim p}\{s_*(x+\xi)\geq\half(\alpha_2-\alpha_1)\}\leq P_\gamma(\alpha_2).\eqno{\mbox{\qed}}\]

\subsection{Proofs for Section \ref{sectNearOpt}}\label{ProofOptMajTest}
\subsubsection{Proof of Proposition \ref{OptMajTest}}
Let $\T^\maj_K$ be $K$-observation majority test for problem $(\S_K)$ associated with $X_1$, $X_2$, and $\P_\gamma$.
{Observe that (\ref{weq14}),  (\ref{eqDDb}), and (\ref{eqDDc}) give $\epsilon_* \leq  {1\over2}-\beta\delta$},
so that (\ref{jeq20}) implies
\be
\risk_{\S}(\T^\maj_K|\P,X_1,X_2)&\leq& \sum_{K\geq k\geq K/2} \left({K\atop k}\right)(1/2-\beta\delta)^k(1/2+\beta\delta)^{K-k}\nn
&\leq&
\sum_{K\geq k\geq K/2} 2^{-K}\left({K\atop k}\right)(1-4\beta^2\delta^2)^k\leq (1-4\beta^2\delta^2)^{K/2}\leq\exp\{-2K\beta^2\delta^2\}.
\ee{riskbin}
In particular, when $\epsilon\in(0,1)$, we have
$$
K\geq {\ln(1/\epsilon)\over 2\beta^2\delta^2} \Rightarrow \risk_{\S}(\T^\maj_K|\P,X_1,X_2)\leq\epsilon,
$$
as claimed in (\ref{eqDDd}).
\par
To prove (\ref{eqDDg}), {assume that the risk of a $K$-observation test $\T_K$ for $(\S_K)$ is} $\leq\epsilon$. Let $x^1_*$, $x^2_*$ form an optimal solution to (\ref{jeq3}), so that $\|x^1_*-x^2_*\|_2=2\delta$. Consider two simple hypotheses stating that the observations $\omega_1,...,\omega_K$ are iid drawn from the distribution of $x^1_*+\xi$, resp., $x^2_*+\xi$, with $\xi\sim q(\cdot)$. Test $\T_K$
decides on these hypotheses with risk $\leq\epsilon$; consequently, assuming w.l.o.g. that $x^1_*+x_*^2=0$, so that $x^1_*=e$, $x^2_*=-e$ with $\|e\|_2=\delta$, we have by Neyman-Pearson lemma
\begin{equation}\label{eqDDe}
\int \min\Big[\underbrace{\prod_{k=1}^Kq(\xi_k+e)}_{q_+(\xi^K)},\underbrace{\prod_{k=1}^Kq(\xi_k-e)}_{q_-(\xi^K)}\Big]\underbrace{d\xi_1...d\xi_K}_{d\xi^K}\leq2\epsilon.
\end{equation}
On the other hand, we have
\begin{equation}\label{eqDDf}
\begin{array}{l}
\left[\int_{\bR^n}\sqrt{q(\xi+e)q(\xi-e)}d\xi\right]^K=\int\sqrt{q_+(\xi^K)q_-(\xi^K)}d\xi^K\\
=
\int\sqrt{\min[q_+(\xi^K),q_-(\xi^K)]\dot\max[q_+(\xi^K),q_-(\xi^K)]}d\xi^K\\
\leq \left[\int\min[q_+(\xi^K),q_-(\xi^K)]d\xi^K\right]^{1/2}\left[\int\max[q_+(\xi^K),q_-(\xi^K)]d\xi^K\right]^{1/2}\\
=\left[\int\min[q_+(\xi^K),q_-(\xi^K)]d\xi^K\right]^{1/2}\left[\int\left[q_+(\xi^K)+q_-(\xi^K)-\min[q_+(\xi^K),q_-(\xi^K)]\right]d\xi^K\right]^{1/2}\\
\leq 2\sqrt{\epsilon(1-\epsilon)},\\
\end{array}
\end{equation}
where the concluding $\leq$ is given by (\ref{eqDDe}) combined with $\epsilon<1/2$.
Since $\|e\|_2=\delta\leq\bar{d}$, (\ref{eqDDa}) combines with (\ref{eqDDf}) to imply that
$
\exp\{-K\alpha\delta^2\}\leq 2\sqrt{\epsilon(1-\epsilon)}\leq \sqrt{4\epsilon},
$
implying  (\ref{eqDDg}) when $\epsilon<\four$. \qed
\subsubsection{Justifying Illustration}
All we need to verify is that in the situation in question,  denoting by $q(\cdot)$ the density of $\N(0,\half I_n)$, we have $q\in\P_\gamma$. Indeed, taking this fact for granted
we ensure the validity of (\ref{eqDDa}) with $\alpha=1$ and $\bar{d}=\infty$. Next, it is immediately seen that the function $\gamma$, see (\ref{eqDDz}), satisfies the relation
$$
\forall (s,0\leq s\leq 1):\gamma(s)\geq \gamma_1(1)={1\over 2\pi},
$$
implying that (\ref{eqDDb}) holds true with $\bar{d}=1$ and $\beta={1\over 2\pi}$, as claimed.
\par
It remains to verify that $q\in\P_\gamma$, which reduces to verifying that the marginal univariate density $\bar{\gamma}(s)={1\over\sqrt{\pi}}\exp\{-s^2\}$ of $q(\cdot)$ satisfies the relation $\int_{\delta}^\infty\bar{\gamma}(s)ds\leq \int_{\delta}^\infty\gamma(s)ds$, $\delta\geq0$, or, which is the same since both $\gamma$ and $\bar{\gamma}$ are even probability densities on the axis, that
$$
\forall(\delta\geq 0): \int_0^\delta \bar{\gamma}(s)ds\geq \int_0^\delta \gamma(s)ds.
$$
The latter relation is an immediate consequence of the fact that the ratio $\bar{\gamma}(s)/\gamma(s)$ is a strictly decreasing function of $s\geq0$ combined with $\int_0^\infty(\bar{\gamma}(s)-\gamma(s))ds=0$.
\subsection{Proof of Proposition \ref{prop0}}\label{prop0proof}
When $x_k\in X_1^k$ and $p_k\in\P^k$ for all $k\leq K$, due to the origin of $h_{k}$, ${\delta_k}$ and $c_{k}$, we have \[
h_k^T(x_k+\xi_k)-c_k\geq \delta_k+h_k^T\xi_k,
\]
{and, because $\eta_k(\cdot)$ is nondecreasing,
\[
\int{\rm e}^{-{\eta_k(h_k^T[x_k+\xi_k]-c_k)}}p_k(\xi_k)d\xi_k\leq \int{\rm e}^{-\eta_k(\delta_k+h_k^T\xi_k)}p_k(\xi_k)d\xi_k\leq \Risk_{\delta_k}(\eta_k|\P^k),\;\;k=1,...,K,
\]
where the concluding ``$\leq$'' is due to (\ref{zeq10}.$a$). Hence,
\[
\begin{array}{r}
\int{\rm e}^{-\phi^{(K)}(x_1+\xi_1,...,x_K+\xi_K)}\prod_{k=1}^K[p_k(\xi_k)d\xi_k]=\prod\limits_{k=1}^K\left[\int{\rm e}^{-\phi_k(x_k+\xi_k)}p_k(\xi_k)d\xi_k\right]\\=\prod\limits_{k=1}^K\left[\int{\rm e}^{-{\eta_k(h_k^T[x_k+\xi_k]-c_k)}}p_k(\xi_k)d\xi_k\right]
\leq\prod_{k=1}^K \left[\Risk_{\delta_k}(\eta_k|\P^k)\right].
\end{array}
\]
}

When $x_k\in X_2^k$ and $p_k\in {\P^k}$ for all $k\leq K$, {in a completely similar way we obtain
\[
\int{\rm e}^{\phi^{(K)}(x_1+\xi_1,...,x_K+\xi_K)}\prod_{k=1}^K[p_k(\xi_k)d\xi_k]
\leq \prod_{k=1}^K\left[\Risk_{\delta_k}(\eta_k|\P^k)\right]. \eqno\hbox{\qed}\\
\]}

\subsection{Proof of Proposition \ref{prop2}}\label{prop2proof}
Let $p(\cdot)\in\P$, let $e\in\bR^n$ be a unit vector, and let $q(\cdot)$ be the probability density of the scalar random variable $e^T\xi$ induced by the density $p(\cdot)$ of $\xi$.
We start with the following well known observation:
\begin{lemma}\label{lemsimp} Let $f$, $g$ be two probability densities on $\bR$ such that
\begin{equation}\label{peq1}
\begin{array}{ll}
(a)&\int\limits_0^\infty f(s)ds=\int\limits_0^\infty g(s)ds,\\
(b)&\int\limits_r^\infty f(s)ds\geq \int\limits_r^\infty g(s)ds,\;\;\forall r\geq0,
\end{array}
\end{equation}
and let $h(s)$ be a nondecreasing real-valued function on the nonnegative ray such that $\int\limits_0^\infty h(s)f(s)ds<\infty$. Then
\begin{equation}\label{peq2}
\int\limits_0^\infty h(s)f(s)ds\geq\int\limits_0^\infty h(s)g(s)ds.
\end{equation}
\end{lemma}
To make the presentation self-contained, here is the proof of the lemma:
\begin{quote}
In view of (\ref{peq1}.$a$), we can assume w.l.o.g. that $h(0)=0$. Let us extend $h(s)$ from the nonnegative ray to the entire real axis by setting $h(s)=0$, $s<0$, thus arriving at a monotone on the axis nonnegative function. Let $\eta\sim f$ and $\zeta\sim g$.
When denoting $H$ the c.d.f. of $h(\eta)$, under the premise of the lemma we clearly have
\[
\bE_{\eta\sim f}\{h(\eta)\}=\int_{0}^\infty t dH(t)=\int_{0}^\infty (1-H(t)) dt=\int_{0}^\infty \Prob\{h(\eta)> t\} dt.
\]
On the other hand, for $t\geq0$ the set $\{s:\;h(s)>t\}$ is a ray either of the form $[a_t,+\infty)$ or $(a_t,+\infty )$ with $a_t\geq0$, so that  (\ref{peq1}.b) implies that
\[
\forall t\geq 0\;\;\Prob\{h(\eta)> t\} 
\geq
\Prob\{h(\zeta)> t\}.
\]
As a result,
\[\bE_{\eta\sim f}\{h(\eta)\}=\int_{0}^\infty \Prob\{h(\eta)> t\}dt\geq \int_{0}^\infty \Prob\{h(\zeta)> t\}dt.
\]
We conclude that $\bE_{\zeta\sim g}\{h(\zeta)\}$ is finite and satisfies
\[
\bE_{\zeta\sim g}\{h(\zeta)\}=\int_{0}^\infty \Prob\{h(\zeta)> t\}dt\leq \bE_{\eta\sim f}\{h(\eta)\}.\eqno{\mbox{\qed}}
\]
\end{quote}

\paragraph{1$^0$.}
Note that $q(\cdot)$ is an even probability density on the axis and  \[
\int\limits_s^\infty [\gamma(r)-q(r)]dr={P_\gamma}(s)-\Prob_{\xi\sim p}\{\xi:\,e^T\xi\geq s\}\left\{\begin{array}{ll}\geq0,&s\geq 0,\\
=0,&s=0.
\end{array}\right.
\]
We have
$$
\begin{array}{l}
\int\limits_{\bR^n}{\rm e}^{-\eta(\delta+e^T\xi)}p(\xi)d\xi= \int\limits_{-\infty}^\infty{\rm e}^{-\eta(\delta+s)}q(s)ds
\underbrace{=}_{(a)}\int\limits_0^\infty H_{{\delta \eta}}(s)q(s)ds
\underbrace{\leq}_{(b)}\int\limits_0^\infty H_{{\delta \eta}}(s)\gamma(s)ds\underbrace{=}_{(c)}\epsilon_\delta(\eta|\gamma),
\\
\end{array}
$$
where $(a)$ is due to the fact that $q$ is even, $(b)$ is a result of applying Lemma \ref{lemsimp} to densities $ \gamma$, $q$ and nondecreasing $H_{{\delta \eta}}$, and $(c)$  is due to the definition (\ref{jeq8}) of the $\delta$-index.

\paragraph{2$^0$.} We have
$$
\begin{array}{rl}
&\int\limits_{\bR^n}{\rm e}^{\eta(-\delta+e^T\xi)}p(\xi)d\xi= {\int\limits_{\bR^n}}{\rm e}^{-\eta(\delta-e^T\xi)}p(\xi)d\xi\hbox{\ [since $\eta(\cdot)$ is odd]}\\
=&\int\limits_{\bR^n}{\rm e}^{-\eta(\delta+e^T\xi)}p(\xi)d\xi \hbox{\ [since $p(\cdot)$ is even]},
\end{array}
$$
and we have already seen in {\bf 1$^0$} that the concluding quantity is $\leq \epsilon_\delta(\eta|\gamma)$. The bottom line is that inequalities (\ref{zeq10}) hold true with $\epsilon=\epsilon_\delta(\eta|\gamma)$, and (\ref{jeq9}) follows. \qed

\subsection{Proof of Proposition \ref{prop33}}\label{sec:proofprop33}
Let $e\in\bR^n$ be a unit vector, $\delta\geq0$, and $p\in \P_\sG$. We have
$$
\int\limits_{\bR^n}{\rm e}^{-\eta(\delta+e^T\xi)}p(\xi)d\xi=\int\limits_{\bR^n}{\rm e}^{-\delta^2-\delta e^T\xi}p(\xi)d\xi\leq {\rm e}^{-\delta^2+\delta^2/2}={\rm e}^{-\delta^2/2},
$$
where the concluding $\leq$ is due to $p\in\P_\sG$ and $\|e\|_2=1$. Similarly,
$$
\int\limits_{\bR^n}{\rm e}^{\eta(e^T\xi-\delta)}p(\xi)d\xi=\int\limits_{\bR^n}{\rm e}^{-\delta^2+\delta e^T\xi}p(\xi)d\xi\leq {\rm e}^{-\delta^2+\delta^2/2}={\rm e}^{-\delta^2/2}.
$$
The resulting inequalities hold true for all unit vectors $e$ and all $p\in\P_\sG$, implying (\ref{zeq20}). \qed
\subsection{Proof of Proposition \ref{prop133}}\label{sec:proofprop133}
1$^o.$ Let $p\in \P$ and $x\in X_1$ be fixed. Due to the monotonicity of $\eta$, and by the definition of $c_*$, we get from (\ref{eq:evar}.a):
\[
\bE_{\xi\sim p}\{\eta(h_*^T(x+\xi)+c)\}\geq \bE_{\xi\sim p}\{\eta(h_*^T(x^1_* + \xi)+c)\} =  c_*+{\e^\kappa\varrho\over 1+\e^\kappa},
\]
and so
\[
\bE_{\xi^K\sim p\times ...\times p}\psi_j([x+\xi_1;...;x+\xi_K])=
\bE_{\xi\sim p}\{\eta(h^T(x+\xi)+c)\}\geq c_*+ {\e^\kappa\varrho\over 1+\e^\kappa}.
\]
On the other hand, by (\ref{eq:evar}.b) we have
\[
\Var_{\xi^K\sim p\times ...\times p}(\psi_j(\omega^K))=m^{-1}\Var_{\xi\sim p}\{\eta(h^T(x+\xi)+c)\}\leq m^{-1}.
\]
Now, by the Chebyshev inequality,
\[
\begin{array}{rcl}\Prob_{\xi^K\sim p\times ...\times p}\left\{\psi_j(\omega^K)<c_*\right\}&\leq&
\Prob_{\xi^K\sim p\times ...\times p}\left\{\psi_j(\omega^K)-\bE_{\xi^K\sim p\times ...\times p}\{\psi_j(\omega^K)
\} <-{\e^\kappa\varrho\over 1+\e^\kappa}\right\}\\
&\leq &
{(1+\e^\kappa)^2\over m\e^{2\kappa}\varrho^2}\leq {1\over 4\e}.\end{array}
\]
As a result, the risk $\risk_{1\S}$ of the test $\T^\mm_K$ satisfies the bound
\[
\risk_{1\S}(\T^\mm_K|\P,X_1,X_2)\leq \sum_{J/2\leq j\le J} \left(\begin{array}{c} J\\j\end{array}\right)(4\e)^{-j}\left(1-{1\over 4\e}\right)^{J-j}\leq 2^J(4\e)^{-J/2}=e^{-J/2}\leq \epsilon_1,
\]
where the final inequality is due to \rf{eq:Kdef}.
\par\noindent
2$^o.$
The same argument, as applied to $x\in X_2$ results in
\[
\bE_{\xi^K\sim p\times ...\times p}\psi_j([x+\xi_1;...;x+\xi_K])=
\bE_{\xi\sim p}\{\eta(h_*^T(x+\xi)+c)\}\leq c_*-{\varrho\over 1+\e^\kappa},
\]
and
\[
\begin{array}{rcl}\Prob_{\xi^K\sim p\times ...\times p}\left\{\psi_j(\omega^K)\geq c_*\right\}&\leq&
\Prob_{\xi^K\sim p\times ...\times p}\left\{\psi_j(\omega^K)-\bE_{\xi^K\sim p\times ...\times p}\{\psi_j(\omega^K)
\} \geq {\varrho\over 1+\e^\kappa}\right\}\\
&\leq &
{(1+\e^\kappa)^2\over m\varrho^2}\leq \four \e^{2\kappa-1}.\end{array}
\]
Same as above, we conclude that
\[
\risk_{2 \S}(\T^\mm_K|\P,X_1,X_2)\leq \e^{J\left(\kappa-\half\right)}=\exp\left(-{J\over 2}{\ln \epsilon_2^{-1}\over \ln \epsilon_1^{-1}}\right)\leq \epsilon_2,
\]
where the concluding inequality is due to \rf{eq:Kdef}.\qed

\end{document}